\documentclass[11pt,english]{article}
\usepackage{a4wide}
\usepackage{amssymb,amsthm}
\usepackage{amsmath}
\usepackage{amsfonts}
\usepackage{latexsym}
\usepackage{amsxtra}
\usepackage{graphicx}
\usepackage{wrapfig}
\usepackage{mathrsfs}
\usepackage{epsfig}
\usepackage{times}
\usepackage{color}
\usepackage{framed}
\usepackage[titles]{tocloft}

\newcommand{\detail}[1]{}
\newcommand{\heikodetail}[1]{}
\newcommand{\xx}{\mbox{$\clubsuit$}}
\def\yy{\mbox{$\spadesuit$}}
\newcommand{\invisible}[1]{%
    \par ($\spadesuit$ \emph{hidden comments in \TeX{} file\/})}
\newcommand{\heikoinvisible}[1]{%
    \par ($\clubsuit$ \emph{hidden comments in \TeX{} file\/})}

\renewcommand{\proof}{\noindent{\sc Proof:}\hskip 1.1em}
\renewcommand{\qed}{\hfill\mbox{$\Box$}\\}

\newcommand{\INT}{{\textnormal{\,int\,}}} 

\newcommand{\s}{\sigma}        %

\newtheorem{lemma}{\bf Lemma}[section]
\newtheorem{theorem}[lemma]{\bf Theorem}
\newtheorem{proposition}[lemma]{\bf Proposition}
\newtheorem{corollary}[lemma]{\bf Corollary}
\newtheorem{definition}[lemma]{\bf Definition}
\newtheorem{REMARK}[lemma]{\bf Remark}
\newtheorem{example}[lemma]{\bf Example}
\newcommand{\Fo}{\,\,\,\text{for }\,\,}
\newcommand{\Foa}{\,\,\,\text{for all }\,\,}

\newcommand{\AND}{\,\,\,\text{and }\,\,}
\newcommand{\OR}{\,\,\,\text{or }\,\,}

\newcommand{\As}{\,\,\,\text{as }\,\,}

\newcommand\Reals{{\mathbb R}}
\newcommand\R{{\mathbb R}}

\newcommand\bbbn{{\mathbb N}}

\newcommand\Sphere{{\mathbb S}}

\newcommand\N{{\mathbb N}}
\renewcommand\S{{\mathbb S}}

\newcommand{\bbbr}{\Reals}

\newcommand{\rd}{\Reals^3}

\newcommand{\stwo}{{\Sphere^2}}

\newcommand\dist{\mathop{\rm dist}\nolimits}
\newcommand\diam{\mathop{\rm diam}\nolimits}

\newcommand\ang{\mathop{\mbox{$<\!\!\!)$}}\nolimits}

\newcommand\Id{{{\rm Id}}}

\newcommand{\F}{\mathcal{F}}

\newcommand{\A}{\mathscr{A}}
\newcommand{\M}{\mathscr{M}}

\newcommand{\Spac}{\mathscr{S}}
\newcommand{\V}{\mathscr{V}}
\newcommand{\B}{\mathscr{B}}
\newcommand{\K}{\mathcal{K}}
\newcommand{\eps}{\varepsilon}
\renewcommand{\H}{\mathscr{H}}

\begin{document}

\title{Integral Menger curvature for surfaces}

\author{Pawe\l{} Strzelecki, Heiko von der Mosel}

\maketitle

\frenchspacing

\begin{abstract}

\detail{\tt\yy Pawel compiles the final version -> bold headings
for theorems}

We develop the concept of integral Menger curvature for  a large
class of nonsmooth surfaces. We prove uniform Ahlfors regularity
and a $C^{1,\lambda}$-a-priori bound for surfaces for which this
functional is finite. In fact, it turns out that there is an
explicit length scale $R>0$ which depends only on an upper bound
$E$ for the integral Menger curvature $\M_p(\Sigma)$ and the
integrability exponent  $p$, and \emph{not\/} on the surface
$\Sigma$ itself; below that scale, each surface with energy
smaller than $E$ looks like a nearly flat disc with the amount of
bending controlled by the (local) $\M_p$-energy. Moreover,
integral Menger curvature can be defined a priori for surfaces
with self-intersections or branch points; we prove that a
posteriori all such singularities are excluded for surfaces with
finite integral Menger curvature.  By means of slicing and
iterative arguments we bootstrap the H\"older exponent $\lambda$ up
to the optimal one, $\lambda=1-(8/p)$, thus establishing a new
geometric `Morrey-Sobolev' imbedding theorem.

As two of the various possible variational applications  we prove
the existence of surfaces in given isotopy classes minimizing
integral Menger curvature with a uniform bound on area, and of
area minimizing surfaces subjected to a uniform bound on integral
Menger curvature.

\vspace{5mm}


\detail{

\bigskip

28A75 =  Length,
area, volume, other geometric measure theory, 49Q10 = optimization
of shapes other than minimal surfaces, 49Q20 = Variational
problems in a geometric measure-theoretic setting, 53A05 =
surfaces in Euclidean space, 53C21 = Methods of Riemannian
geometry, including PDE methods; curvature restrictions, 46E35 =
Sobolev spaces and other spaces of ``smooth'' functions, embedding
theorems, trace theorems; cf. the classification in MathRev of
David and Semmes on quasiminimal sets

\bigskip

}

\end{abstract}

\tableofcontents

\clearpage

\renewcommand\theequation{{\thesection{}.\arabic{equation}}}
\def\setnumbers{\setcounter{equation}{0}}


\section{Introduction}\label{sec:1}

For three different non-collinear points $x,y,z\in\R^n$ the
expression
\begin{equation}\label{one}
R(x,y,z):=\frac{|x-y||x-z||y-z|}{4A(x,y,z)},
\end{equation}
where $A(x,y,z)$ is the area of the triangle with vertices at $x$,
$y$ and $z$, provides the radius of the uniquely defined
circumcircle through $x,y,$ and $z.$ This gives rise to a family
of {\it integral Menger curvatures}\footnote{Coined after K.
Menger who generalized expression \eqref{one} to metric spaces as
a foundation of a metric coordinate free geometry; see
\cite{menger}, \cite{BlM70}.}, that is, geometric curvature
energies of the form
\begin{equation}\label{two}
\M_p(E):=\int_E\int_E\int_E\frac{1}{R^p(x,y,z)}\,d\H^1(x)\,
d\H^1(y)\, d\H^1(z),\quad p\ge 1,
\end{equation}
defined on one-dimensional Borel sets $E\subset \R^n.$ According
to a remarkable result of J.C.~L\'eger \cite{leger} such sets $E$
with Hausdorff measure $\H^1(E)\in (0,\infty)$ and with finite
integral Men\-ger curvature $\M_p(E)$ for some $p\ge 2$, are
$1$-re\-cti\-fi\-able in the sense of geometric measure theory. To
be precise, $\H^1$-almost all of $E$ is contained in a countable
union of Lipschitz graphs. Ahlfors-regular\footnote{A set $E$ of
Hausdorff dimension $1$ is said to be {\it Ahlfors-regular} if and
only if there is a constant $C_E\ge 1$ such that
$C_E^{-1}R\le\H^1(E\cap B(x,R))\le C_ER$ for every $x\in E$ and
$R>0,$ where $B(x,R)$ denotes a closed ball of radius $R$.}
one-dimensional Borel sets $E\subset\R^2$ satisfying the local
condition
\begin{equation}\label{three}
\M_2(E\cap B(\xi,r))\le Cr\quad\Foa \xi\in\R^2,\, r\in (0,r_0]
\end{equation}
turn out to be {\it uniformly rectifiable}, i.e., they are
contained in the graph of {\it one} bi-Lipschitz map
$f:\R\to\R^2$; see \cite[Theorem 39]{P02} referring to work of P.
Jones. M. Melnikov and J. Verdera \cite{Mel95}, \cite{MelV95} realized
that $\M_2$ is a crucial quantity in harmonic analysis to
characterize removable sets for bounded analytic functions; see
e.g. the surveys \cite{Ma98}, \cite{Ma04}, \cite{V01}.

If one considers the $\M_p$-energy on rectifiable closed curves
$E=\gamma(\S^1)\subset\R^3$ the following {\it geometric
Morrey-Sobolev imbedding theorem} was proven in \cite[Theorem
1.2]{ssvdm-triple}, and this result may be viewed as a counterpart
to J.C.~L\'eger's regularity result on a higher regularity level:

\medskip

\noindent {\it If $\M_p(\gamma)$ is finite for some $p\in
(3,\infty]$ and if the arclength parametrization $\Gamma$ of the
curve $\gamma$ is a local homeomorphism then $\gamma(\S^1)$ is
diffeomorphic to the unit circle $\S^1$, and $\Gamma$ is a finite
covering of $\gamma(\S^1)$ of class $C^{1,1-(3/p)}.$ }

\medskip

In fact, even the stronger local version holds true \cite[Theorem
1.3]{ssvdm-triple}, which may be viewed as a {\it geometric
Morrey-space imbedding} and whose superlinear growth assumption
\eqref{three_a} is the counterpart of \eqref{three}:

\medskip

\noindent {\it If the arclength parametrization $\Gamma$ is a
local homeomorphism, and if
\begin{equation}\label{three_a}
\int_{B(\tau_1,r)} \int_{B(\tau_2,r)} \int_{B(\tau_3,r)}\frac{ds\,
dt\, d\sigma}{R^2(\Gamma(s),\Gamma(t),\Gamma(\sigma))}\le
Cr^{1+2\beta}
\end{equation}
holds true for all $r\in (0,r_0]$ and all arclength parameters
$\tau_1,\tau_2,\tau_3$, then $\Gamma $ is a $C^{1,\beta}$-covering
of the image $\gamma(\S^1)$ which itself is diffeomorphic to the
unit circle. }

\medskip

From the results on one-dimensional sets and in particular on
curves it becomes apparent that integral Menger curvature $\M_p$
exhibits regularizing and self-avoidance effects (as already
suggested in \cite{GM99} and \cite{banavargmm}). These effects
become stronger with increasing $p$, in fact, one has
$$
\lim_{p\to\infty}\left(\M_p(\gamma)\right)^{1/p}=\frac{1}{\inf\limits_{\sigma\not=
s\not= t\not= \sigma}R(\Gamma(s),\Gamma(t),\Gamma(\sigma))}=:
\frac{1}{\triangle[\gamma]},
$$
where $\triangle[\gamma]$ is the notion of {\it thickness} of
$\gamma$ introduced by O. Gonzalez and J.H. Maddocks \cite{GM99} who were
motivated by analytical and computational issues arising in the
natural sciences such as the modelling of knotted DNA molecules.
In fact, it was shown in \cite{GMSvdM02} and \cite{SvdM03a} that
closed curves with finite energy $1/\triangle[\gamma]$, i.e. with
positive thickness, are {\it exactly} those embeddings with a
$C^{1,1}$-arclength parametrization, which lead to variational
applications for nonlinearly elastic curves and rods with positive
thickness; see also \cite{SvdM03b}, \cite{SvdM04}. We generalized
this concept of thickness in \cite{StvdM05} and \cite{svdm-global}
to a fairly general class of nonsmooth surfaces
$\Sigma\subset\R^n$ with the central result {\it that surfaces
with positive thickness $\triangle[\Sigma]$ are in fact
$C^{1,1}$-manifolds with a uniform control on the size of the
local $C^{1,1}$-graph patches depending only on the value of
$\triangle[\Sigma]$.} Uniform estimates on sequences then allow
for the treatment of various energy minimization problems in the
context of thick (and therefore embedded) surfaces of prescribed
genus or isotopy class; see \cite[Theorem~7.1]{svdm-global}.

In the present situation we ask the question:

\begin{quote}
\noindent {\it Is it possible to extend the definition of
integral Menger curvature $\M_p$ for $p<\infty$ to surfaces with
similar regularizing and self-avoidance effects as in the curve
case?}
\end{quote}

The most natural generalization of $\M_p$ to two-dimensional
closed surfaces $\Sigma\subset\R^3$ would be to replace the
circumcircle radius $R(x,y,z) $ of three points $x,y,z$ in
\eqref{two} by the circumsphere radius $R(\xi,x,y,z)$ of the
tetrahedron $T:=(\xi,x,y,z)$ spanned by the four non-coplanar
points $\xi,x,y,z$. This radius is given by
\begin{equation}
\label{1/R-sphere} \frac{1}{2R(T)} = \frac{\bigl|\langle z_3,
z_1\times z_2 \rangle\bigr|}{\bigl| \, |z_1|^2 z_2 \times z_3 +
|z_2|^2 z_3 \times z_1 + |z_3|^2 z_1 \times z_2\, \bigr| },
\end{equation}
where $z_1=\xi-z$, $z_2=x-z$, $z_3=y-z$. This would lead to the
geometric curvature energy
\begin{equation}\label{five}
\int_\Sigma\int_\Sigma\int_\Sigma\int_{\Sigma}\frac{d\H^2(\xi)
d\H^2(x) d\H^2(y) d\H^2(z)}{R^p(\xi,x,y,z)},
\end{equation}
which in principle would serve our purpose: all our results stated
below extend to this energy. But -- although the integrand is
trivially constant if $\Sigma$ happens to be a round sphere --
there are smooth surfaces with straight no\-dal lines (such as the
graph of the function $f(x,y):=xy$) where the integrand is not
pointwise bounded; see Appendix B. This is a problem since we want
to consider arbitrarily large $p$, and we envision a whole family
of integral Menger curvatures that are finite on {\it any} closed
smooth surface for {\it any\/} value of $p$.

Rewriting \eqref{one} as
$$
R(x,y,z)=\frac{|x-z||y-z|}{2\dist(z,L_{xy})},
$$
where $L_{xy}$ denotes the straight line through $x$ and $y$, one
is naively tempted to consider {\it $4$-point-integrands} of the
form
\begin{equation}\label{six}
\left(\frac{\dist(\xi,\langle x,y,z\rangle)}{M(|\xi-x|,|\xi-y|,
|\xi-z|)^\alpha}\right)^p,
\end{equation}
where $\langle x,y,z\rangle $ denotes the affine $2$-plane through
generic non-collinear points $x,y,z\in\R^3$. Here, $\alpha\ge 1$
is a power and the function $M:\R_+\times\R_+\times \R_+\to\R_+$
is a {\it mean}, i.e. $M$ is monotonically increasing with respect
to each variable and satisfies the inequality
$$
\min\{a,b,c\}\le M(a,b,c)\le\max\{a,b,c\}.
$$
Again, all our results that will be stated below would hold if we
worked with integrands as in \eqref{six} for $\alpha=2$. This is
very similar to a suggestion of J.C.~L\'eger \cite[p. 833]{leger}
who proposes a general integrand for $d$-dimensional sets; for
$d=2$ his choice boils down to \eqref{six} with $M$ being the
geometric mean and $\alpha=3$.  However, the situation for such
integrands, due to the lack of symmetry w.r.t.  permutations of
the 4 points, is even worse than for inverse powers of the
circumsphere radius: for any choice of $\alpha>1$ there are
sufficiently large $p=p(\alpha)$ such that  even a round sphere
has infinite energy. This singular behaviour is  caused by small
tetrahedra for which the plane through $(x,y,z)$ is almost
perpendicular to the surface.  See Appendix B for more details.

\detail{

\bigskip

Incidentally: these are tetrahedra for which the ratio of
${\dist(\xi,\langle x,y,z\rangle)}$ (one of the heights) to
${M(|\xi-x|,|\xi-y|, |\xi-z|)^2}$ (which, roughly, measures the
average area of a face) ceases to be comparable with, say, volume
of $T$ divided by $(\diam T)^4$.

\bigskip

} Roughly speaking, the trouble with \eqref{1/R-sphere} or
\eqref{six} for surfaces comes from the fact that various
`obviously equivalent' formulae for $1/R$ for triangles (relying
on the sine theorem) are no longer equivalent for tetrahedra; to
obtain a whole scale of surface integrands which penalize
wrinkling, folding, appearance of narrow tentacles,
self-intersections etc. but stay bounded on smooth surfaces, one
should make a choice here.  In their pioneering work
\cite{LW08a,LW08b} dealing with $d$-rectifiability and least
square approximation of $d$-regular measures, G. Lerman and J.T.
Whitehouse suggest a whole series of high-dimensional counterparts
of the one-dimensional Menger curvature. Their ingenious discrete
curvatures are based, roughly speaking, on the so-called polar
sine function scaled by some power of the diameter of the simplex,
and can be used to obtain powerful and very general
characterizations of rectifiability of measures. (In
\cite[Sec.~1.5 and~6]{LW08b} they also note that the integrand
suggested by L\'{e}ger does not fit into their setting.)

Motivated by this and by the explicit formula for the circumsphere
radius, we are led to consider another $4$-point integrand with
symmetry and with fewer cancellations in the denominator. For a
tetrahedron $T$ consider the function
\begin{equation}\label{seven}
\K(T):=\begin{cases}\displaystyle \frac{V(T)}{A(T)(\diam T)^2} &
\textnormal{if the vertices of $T$ are not
co-planar},\\
0 & \textnormal{otherwise,}
\end{cases}
\end{equation}
where $V(T)$ denotes  the volume of $T$ and $A(T)$ stands for the
total area, i.e., the sum of the areas of all four triangular
faces of $T$.  Thus, up to a constant factor $\K$ is the ratio of
the minimal height of $T$ to the square of its diameter, which is
similar but not identical to the numerous curvatures considered by
Lerman and Whitehouse in \cite{LW08b}. The difference is that our
$\K$ scales like the inverse of length whereas their
$d$-dimensional curvatures, cf. e.g. the definition of
$c_{\text{MT}}$ in \cite[p. 327]{LW08b}, for $d=2$ scale like the
inverse of the \emph{cube\/} of length. Such scaling enforces too
much singularity for our purposes; we explain that in
Remark~\ref{badscaling} in Section 5.

This leads us to the {\it integral Menger curvature for
two-dimensional surfaces $\Sigma\subset\R^3$},
\begin{equation}\label{eight}
\M_p(\Sigma):=\int_\Sigma\int_\Sigma\int_\Sigma\int_\Sigma\K^p(T)
\,d\H^2\otimes d\H^2\otimes d\H^2\otimes d\H^2(T),
\end{equation}
which is finite for any $C^{2} $-surface for any finite $p$, since
$\K(T)$ is bounded on the set of all nondegenerate tetrahedra with
vertices on such a surface; see Appendix A.

To keep a clear-cut situation in the introduction we state our
results here for closed Lipschitz surfaces only and refer the
reader to Definition \ref{umbrella} in Section \ref{adm} for the
considerably more general class $\A$ of admissible surfaces, and
to Sections \ref{sec:3}, \ref{sec:5}, and \ref{sec:6} for the
corresponding theorems in full generality. Let us just remark,
however, that our admissibility class  $\A$ contains surfaces that
are not even topological submanifolds of $\bbbr^3$: e.g. a sphere
with the north and south pole glued together. The finiteness of
$\M_p(\Sigma)$ has therefore topological, measure-theoretic and
analytical consequences.

\begin{theorem}[Uniform Ahlfors regularity and a diameter bound]
\label{thm:1.1} There exists an absolute constant $\alpha >0$ such
that for any $p>8$, every $E>0$, and for every closed compact and
connected Lipschitz surface $\Sigma\subset\R^3$ with
$\M_p(\Sigma)\le E $ the  following estimates hold:
\begin{gather}
\diam
\Sigma\ge\left(\frac{\alpha^{5p}}{E}\right)^{\frac{1}{p-8}},\nonumber \\
\label{eleven} \H^2(\Sigma\cap B(x,R))\ge\frac{\pi}{2}R^2\quad\Foa
x\in\Sigma\, \mbox{ and }\, R\in (0,(\alpha^{5p}/E)^{1/(p-8)}]
\end{gather}
\end{theorem}
General Lipschitz surfaces may have conical singularities with a
very small opening angle, but finite $\M_p$-energy controls
uniformly the lower density quotient. These quantitative lower
estimates for diameter and density quotient resemble L. Simon's
results \cite[Lemma 1.1 and Corollary 1.3]{Si93} for smooth embedded
two-dimensional surfaces of finite Willmore energy, derived by
means of the first variation formulas. Here, in contrast, we set
up an intricate algorithm (see Theorem \ref{goodtetra} and its
proof in Section~\ref{sec:4}), starting with a growing double cone
and continuing with an increasingly  complicated growing set
centrally symmetric to a surface point, to scan the possibly highly
complex exterior and interior domain bounded by $\Sigma$ in
search for three more complementing surface points to produce a
``nice'' tetrahedron whose size is controlled in terms of the
energy. Along the way, the construction allows for large
projections onto affine $2$-planes which leads to the uniform estimate
\eqref{eleven}.

The case $p=8$ yields a result which may be interpreted as a
two-dimensional variant of Fenchel's theorem on the total
curvature of closed curves \cite{fenchel}:
\begin{theorem}[Fenchel for surfaces]\label{thm:1.2}
There is an absolute constant $\gamma_0>0$ such that
$\M_8(\Sigma)\ge\gamma_0$ for any closed compact connected
Lipschitz surface $\Sigma\subset\R^3$.
\end{theorem}

The exponent $p=8$ is a limiting one here: $\M_8$ is scale
invariant. Invoking scaling arguments, it is easy to see that any
cone over a smooth curve must have infinite $\M_p$-energy for
every $p\ge 8$.

Uniform control over the lower Ahlfors regularity constant as in
Theorem \ref{thm:1.1} permits us to prove the existence of a field
of tangent planes for finite energy surfaces $\Sigma$ (coinciding
with the classical tangent planes at points of differentiability
of $\Sigma$), and quantitatively control its oscillation:
\begin{theorem}[Oscillation of the  tangent planes]\label{thm:1.3}
For any closed compact and connected Lipschitz surface
$\Sigma\subset\R^3$ with $\M_p(\Sigma)\le E$ for some $p>8$ the
tangent plane $T_x\Sigma$ is defined everywhere and depends
continuously on $x$: there are constants $\delta_1=\delta_1(p)>0$
and $A=A(p)\ge 0$ such that
\begin{equation}\label{twelve}
\ang (T_x\Sigma, T_y\Sigma) \le A
E^{\frac{1}{p+16}}|x-y|^{\frac{p-8}{p+16}}
\end{equation}
whenever $|x-y|\le \delta_1(p)E^{-1/(p-8)}$.

\end{theorem}

\detail{\yy\xx In particular, $\Sigma$ is an orientable
$C^{1,\kappa}$-manifold for $\kappa=(p-8)/(p+16)$. \yy}

One might compare this theorem with Allard's famous regularity
theorem \cite[Theorem 8.19]{allard} for varifolds: Supercritical
integrability assumptions (with exponent $p>$ dimension) on the
generalized mean curvature are replaced here by integrability
assumptions on our four-point Menger curvature integrand $\K$ for
$p>8=4\cdot $ dimension -- with possible extension to metric
spaces, since our integrand may be expressed in terms of distances
only. To prove Theorem~\ref{thm:1.3} (see Section~\ref{sec:5} for
all details), we start with a technical lemma, ascertaining that
the so-called P. Jones' $\beta$-numbers of $\Sigma$, measuring the
distance from $\Sigma$ to the best `approximating plane' and
defined by
\begin{equation}
\label{Jones-intro} \beta_\Sigma(x,r) := \inf\left\{ \sup_{y\in
\Sigma\cap B(x,r)} \frac{\dist (y,F)}{r} \quad \colon \quad
\mbox{$F$ is an affine $2$-plane through $x$} \right\},
\end{equation}
can be estimated by $\mathrm{const}\cdot r^{(p-8)/(p+16)}$ at
small scales. This estimate is uniform, i.e. depends only on $p$
and on the energy bounds, due to Theorem~\ref{thm:1.1}. For a wide
class of \emph{Reifenberg flat sets with vanishing constant\/},
see G. David, C. Kenig and T. Toro \cite{davidkenigtoro} or D.
Preiss, X. Tolsa and Toro \cite{ptt}, this would be enough to
guarantee the desired result. However, at this stage we cannot
ensure that the surface we consider is Reifenberg flat with
vanishing constant; it might be just a Lipschitz surface with some
folds or conical singularities which are not explicitly excluded
in Theorem~\ref{thm:1.1}. Reifenberg flatness, introduced by
 E.R. Reifenberg \cite{reifenberg} in his famous paper on the Plateau
problem in high dimensions, requires not only some control of
$\beta$'s, but also a stronger fact: one needs to know that the
Hausdorff distance between the approximating planes and $\Sigma$
is small at small scales. To get such control, we use some
elements of the proof of Theorem~\ref{thm:1.1} to guarantee the
existence of large projections of $\Sigma$ onto planes, and,
proceeding iteratively, combine this with the decay of $\beta$'s
to reach the desired conclusion. The proof is presented in
Section~\ref{sec:5}; it is self-contained and independent of
\cite{davidkenigtoro} and \cite{ptt}.

Once Theorem~\ref{thm:1.3} is established, we know that in a small
scale, depending solely on $p$ and on the energy bound, the
surface is a graph of a $C^{1,\kappa}$ function. Slicing arguments
similar to, but technically more intricate than those in the proof
of optimal H\"older regularity for curves in \cite[Theorem
1.3]{ssvdm-triple}, are employed in Section \ref{sec:6} to
bootstrap the H\"{o}lder exponent from $\kappa= (p-8)/(p+16)$ to
$(p-8)/p$ and prove
\begin{theorem}[Optimal H\"older exponent]\label{thm:1.4}
Any closed, compact and connected Lipschitz surface $\Sigma$ in
$\R^3$ with $\M_p(\Sigma)\le E<\infty$ for some $p>8$ is an orientable
$C^{1,1-(8/p)}$-manifold with local graph representations whose
domain size is controlled solely in terms of $E$ and $p$.
\end{theorem}
We expect that $1-8/p$  is the optimal exponent,  like the
corresponding optimal exponent $1-3/p$ in the curve case in
\cite[Theorem 1.3]{ssvdm-triple}; see the example  for
curves in \cite[Section 4.2]{marta}.

The last section deals with sequences of surfaces with a uniform
bound on their $\M_p$-energy. Using a combination of Blaschke's
selection theorem and Vitali covering arguments with balls on the
scale of uniformly controlled local graph representations we can
establish the following compactness result.
\begin{theorem}[Compactness for surfaces with equibounded ${\M_p}$-energy]
\label{thm:1.5} Let $\{\Sigma_j\}$ be a sequence of closed,
compact and connected Lipschitz surfaces containing $0\in\R^3$
with
$$
\M_p(\Sigma_j)\le E \Foa
j\in\N\quad\And\quad\sup_{j\in\N}\H^2(\Sigma_j)\le A,
$$
for some $p>8$. Then there is a compact $C^{1,1-8/p}$-manifold
$\Sigma$ without boundary embedded in $\R^3$, and a subsequence
$j'$, such that $\Sigma_{j'}$ converges to $\Sigma$ in $C^1$, and
such that
$$
\M_p(\Sigma)\le \liminf_{j\to\infty}\M_p(\Sigma_j)
$$
\end{theorem}
Instead of the uniform area bound one could also assume a uniform
diameter bound.

Using this compactness result and the self-avoidance effects of
integral Menger curvature we will prove that one can minimize area
in the class of closed, compact and connected Lipschitz surfaces
of fixed genus under the constraint of equibounded energy. For
given $g\in\N$ let $M_g$ be a closed, compact and connected
reference surface of genus $g$ that is smoothly embedded in
$\R^3$, and consider the class $\mathscr{C}_E(M_g)$ of closed,
compact and connected Lipschitz surfaces $\Sigma \subset\R^3 $
ambiently isotopic to $M_g$ with $\M_p(\Sigma)\le E$ for all
$\Sigma\in\mathscr{C}_E(M_g)$.
\begin{theorem}[Area minimizers in a given isotopy class]
\label{thm:1.6} For each $g\in\N$, $E>0$ and each fixed reference
surface $M_g$ the class $\mathscr{C}_E(M_g)$ contains a surface of
least area.
\end{theorem}
We can also minimize the integral Menger curvature $\M_p$ itself
in a given isotopy class with a uniform area bound, i.e. in the
class $\mathscr{C}_A(M_g)$ of closed, compact and connected
Lipschitz surfaces $\Sigma\subset\R^3$ ambiently isotopic to $M_g$
with $\H^2(\Sigma)\le A<\infty.$
\begin{theorem}[$\M_p$-minimizers in a given isotopy class]
\label{thm:1.6a} For each $g\in\N$, $A>0$, there exists a surface $\Sigma
\in\mathscr{C}_A(M_g)$ with
$$
\M_p(\Sigma)=\inf_{\mathscr{C}_A(M_g)}\M_p.
$$
\end{theorem}

The proofs of Theorems~\ref{thm:1.5}--\ref{thm:1.6a} are given in
Section~7.
\begin{REMARK}\rm
It can be checked that our Theorems \ref{thm:1.1} and
\ref{thm:1.2} can be proven for a large class of integrands
including the two-dimensional $c_{\text{MT}}$ and other curvatures
of Lerman and Whitehouse, and even the one suggested by L\'eger.
(One just has to check what is the critical scaling-invariant
exponent, and work above this exponent.) However, Theorems
\ref{thm:1.3} and \ref{thm:1.4}, and consequently also Theorems
\ref{thm:1.5}, \ref{thm:1.6}, and \ref{thm:1.6a}  seem to fail for
any choice of integrand $\K_s(T)$ which scales like the inverse of
length to some power $1+s$, $s>0$. Such a choice enforces too much
singularity for large $p$, and the methods we employ to prove
H\"{o}l\-der regularity of the unit normal show that the only
surface with $\int \K_s^p \, d\mu$ finite for all $p$ would be (a
piece of) the flat plane. See Remark~\ref{badscaling} in
Section~5.
\end{REMARK}

\begin{REMARK} \rm Our work is related to the theory of
uniformly rectifiable sets of G. David and S. Semmes, see their
monograph \cite{davidsemmes}. Numerous equivalent definitions of
these sets involve subtle conditions stating how well, in an
average sense, the set can be approximated by planes. One of the
deep ideas behind this is to try and search for the analogies
between classes of sets and function spaces. It turns out then
that various approximation or imbedding theorems for function
spaces have geometric counterparts for sets, see e.g. the
introductory chapter of \cite{davidsemmes}. Speaking naively and
vaguely, David and Semmes work in the realm which corresponds to
the subcritical case of the Sobolev imbedding theorem: there is no
smoothness but subtle tools are available to give nontrivial
control of the rate of approximation of a function by linear
functions (or rather: a set by planes). Here, we are in the
supercritical realm. For exponents larger than the critical $p=8$
related to scale invariance, excluding conical singularities,
finiteness of our curvature integrands gives continuity of tangent
planes, with precise local control of the oscillation. Note that
the exponent $1-8/p$ in Theorem~\ref{thm:1.6} is computed
according to Sobolev's recipe: the domain of integration has
dimension 8 and we are dealing with the $p$'th power of
`curvature'.
\end{REMARK}

\medskip

\noindent {\bf Acknowledgement.}\, The  authors  would like to
thank the Deutsche Forschungsgemeinschaft, Polish Ministry of
Science and Higher Education, and the Centro di Ricerca Matematica
Ennio De Giorgi at the Scuola Normale Superiore in Pisa, in
particular Professor Mariano Giaquinta, for generously supporting
this research. We are also grateful to Professor Gilad Lerman for
his comments on an earlier version of this paper which helped us
to improve the presentation.


\section{Notation. The class of admissible surfaces}

\setnumbers

\subsection{Basic notation}

\label{basic}

\textbf{Balls, planes and slabs.} $B(a,r)$ denotes always the
\emph{closed\/} ball of radius $r$, centered at $a$. When $a=0\in
\bbbr^3$, we often write just $B_r$ instead of $B(0,r)$.

For non-collinear points $x,y,z\in\R^3$ we denote by $\langle
x,y,z\rangle$ the affine $2$-plane through $x,y, $ and $z.$ If $H$
is a 2-plane in $\rd$, then $\pi_H$ denotes the orthogonal
projection onto $H$. For an affine plane $F\subset \rd$ such that
$0\not \in F$, we write $\sigma_F$ to denote the central
projection from $0$ onto $F$.

If $F$ is an affine plane in $\rd$ and $d>0$, then we
denote the infinite slab about $F$ by
\[
U_d(F) := \{y\in \rd\colon \dist (y,F) \le d\}.
\]

\smallskip\noindent\textbf{Cones.}  Let $\varphi \in (0, \frac \pi 2)$ and $w\in
\Sphere^2$. We set
\[
C(\varphi,w) := \{y\in \rd \, \colon \, |y\cdot w | \ge
|y|\cos \varphi\}
\]
describing the infinite double-sided cone of opening angle $2\varphi$ whose axis
is determined by $w$,
and we define
$
C_r(\varphi,w):= B(0,r)\cap C(\varphi,w).
$
We also distinguish between the two conical halves
\[
C^+(\varphi,w) := \{y\in \rd \, \colon \,  y\cdot w  \ge
|y|\cos \varphi\}, \qquad C^-(\varphi,w):=\{y\in \rd \, \colon \,
- y\cdot w  \ge |y|\cos \varphi\},
\]
and set
$
C_r^\pm(\varphi,w):= B(0,r)\cap C^\pm(\varphi,w).
$

\smallskip\noindent\textbf{Rotations in $\rd$.} Throughout, we fix an
orientation of $\rd$. Assume that $u,v\in \Sphere^2$ are
orthogonal and $u\times v = w\in\stwo$. We write $R(\varphi,w)$ to
denote the rotation which, in the orthonormal basis $(u,v,w)$, is
represented by the matrix
\[
\left(
\begin{array}{ccc}
\cos \varphi & -\sin \varphi & 0 \\
\sin \varphi & \cos \varphi & 0 \\
0 & 0 & 1
\end{array}
\right)\, .
\]
Note that this formula gives in fact a linear map which does not
depend on the choice of orthonormal vectors $u,v$ with $u\times v=
w$.

\smallskip\noindent\textbf{Segments.} Whenever
$z\in \rd$, $s > 0$ and $w\in \stwo$, we set
\[
I_{s,w}(z):=\{z+ tw \, \colon \, |t|\le s\}
\]
(this is the segment of length $2s$, centered at $z$ and parallel
to $w$).

\smallskip\noindent\textbf{Tetrahedra.} Since we deal with
an inte\-grand defined on quadruples of points in $\rd$, and in
various places we need to estimate that integrand on specific
quadruples satisfying some additional conditions, we introduce
some notation now to shorten the statements of several results in
Sec. 3--6.

By a \emph{tetrahedron\/} $T$ we mean a quadruple of points,
$T=(x_0,x_1,x_2,x_3)$ with $x_i\in \rd$ for $i=0,1,2,3.$
By a triangle $\Delta$ we
mean a triple of points, $\Delta=(x_0,x_1,x_2)$, $x_i\in \rd$. We
say that $\Delta=\Delta(T)$ is the base of $T$ iff
$\Delta=(x,y,z)$ and $T=(x,y,z,w)$ for some $x,y,z,w\in \rd$.

For $T=(x_0,x_1,x_2,x_3)$ and $T=(x_0',x_1',x_2',x_3')$ we set
\[
\|T-T'\| := \min_{\s\in S_4}\left[\max_{0\le i\le 3} |x_i-x_i'|\right],
\]
where $|x_i-x_i'|$ denotes the Euclidean norm and $S_4$ is the symmetric
group  of all
permutations of sets with four elements. We write
$\B_r(T):=\{T'\colon \|T-T'\|\le r\}$.

To investigate the local and global behaviour of a surface, we
often estimate its $\M_p$-energy on $\B_\eps(T)\cap \Sigma$ where
either $T$ resembles, roughly speaking, a regular tetrahedron or
at least its base $\Delta(T)$ resembles, again roughly, a regular
triangle. Here are the appropriate definitions.

\begin{definition}
\label{Def-Ved} Let $\theta\in (0,1)$ and $d>0$. We say that
$T=(x_0,x_1,x_2,x_3)$ is $(\theta,d)$--voluminous, and write $T\in
\V(\theta,d)$, if and only if
\begin{enumerate}
\renewcommand{\labelenumi}{{\rm (\roman{enumi})}}

\item $x_i\in B(x_0,2d)$ for all $i=1,2,3$;

\item
$\theta d \le |x_i-x_j|$ \, for all $i\not =j$, where $i,j=0,1,2,3$;

\item $\ang(x_1-x_0,x_2-x_0) \in [\theta,\pi-\theta]$;

\item $\dist \bigl(x_3, \langle x_0,x_1,x_2\rangle\bigr)\ge \theta
d$.
\end{enumerate}
\end{definition}

\begin{definition} Let $\theta\in (0,1)$ and $d>0$.
We say that $\Delta=(x_0,x_1,x_2)$ is $(\theta,d)$--wide, and
write $\Delta \in \Spac (\theta,d)$, if and only if
\begin{enumerate}
\renewcommand{\labelenumi}{{\rm (\roman{enumi})}}

\item $x_i\in B(x_0,2d)$ for $i=1,2$;

\item
$\theta d \le |x_i-x_j|$ \, for $i\not =j$, where $i,j=0,1,2$;

\item $\ang(x_1-x_0,x_2-x_0) \in [\theta,\pi-\theta]$.
\end{enumerate}
\end{definition}

\noindent
\begin{REMARK}\rm  Similar classes of simplices have been used by
Lerman and Whitehouse, see \cite[Sec.~3]{LW08b}. The class of $T$
with $\Delta(T)\in \Spac(\theta,d)$ differs from their class of
$2$-separated tetrahedra as the minimal face area of $T$ with
$\Delta(T)\in \Spac(\theta,d)$ does not have to be comparable to
the square of $\diam T$. This plays a role in Section~5 and
Section~6.
\end{REMARK}

\subsection{The class of admissible surfaces}

\label{adm}

Throughout  the paper we consider only compact and closed
surfaces.

\begin{definition} We say that a compact connected subset $\Sigma\subset \rd$
such that $\Sigma = \partial U$ for some bounded domain $U\subset
\rd$ is \emph{an admissible surface\/}, and write $\Sigma\in \A$,
if the following two conditions are satisfied: \label{umbrella}
\begin{enumerate}
\renewcommand{\labelenumi}{(\roman{enumi})}

\item[\rm (i)] There exist a constant $K=K(\Sigma)$ such that
\[
\infty > \H^2(\Sigma\cap B(x,r)) \ge K^{-1} r^2 \qquad\mbox{for
all $x\in \Sigma$ and all $0<r\le \diam\Sigma$;}
\]

\item[\rm (ii)]
There exists a dense subset $\Sigma^\ast \subset \Sigma$
with the following property: for each $x\in \Sigma^\ast$ there
exists a vector $v=v(x)\in \S^2$ and a radius
$\delta_0=\delta_0(x)>0$ such that
\[
B(x,\delta_0) \cap \bigl(x+C^+(\pi/4,v)\bigr) \subset U \cup
\{x\}, \qquad B(x,\delta_0) \cap \bigl(x+C^-(\pi/4,v)\bigr)
\subset (\rd\setminus\overline{U})\cup \{x\}\, .
\]
\end{enumerate}
\end{definition}

\medskip

Condition (ii) seems to be rather rigid because of the symmetry
requirement. We could have used some smaller angle $\varphi_0$
instead of $\varphi_0=\pi/4$ with the only effect that the
absolute constants in Theorems~\ref{Thm:3.1}--\ref{goodtetra},
\ref{thm:5.3}, and \ref{improved} would change, but we stick to
$\varphi_0=\pi/4$ for the sake of simplicity.

Condition (i) excludes sharp cusps around an isolated point of
$\Sigma$ but allows for isolated conical singularities and various
cuspidal folds along arcs.

Note that this is a large class of surfaces, and if $\Sigma\in
\A$, then $\Sigma$ does not have to be an embedded topological
manifold. Consider for example a sphere on which two distinct
points have been identified, or, more generally, a sphere with
$2N$ distinct smooth arcs  and identify $N$ pairs of these arcs.

Here are further examples.

\begin{example}[$C^1$ surfaces]\label{C1-surf} \rm
If $\Sigma$ is a $C^1$ manifold which bounds a domain $U$, then
$\Sigma\in \A$. One can take $\Sigma^\ast \equiv \Sigma$; by
definition of differentiability, for each point $x\in \Sigma$
condition (ii) is satisfied for $v(x)={}$ the inner normal to
$\Sigma$ at $x$, and one can choose a uniform lower bound for the
numbers $\delta_0(x)$, i.e. we can always pick a $\delta_0(x)\ge
\delta_0=\delta_0(\Sigma)>0$.
\end{example}

\begin{example}[Lipschitz surfaces] \rm
If $\Sigma=\partial U$ is a Lipschitz manifold, then $\Sigma\in
\A$. We can take $\Sigma^\ast={}$ the set of all points where
$\Sigma$ has a classically defined tangent plane. By Rademacher's
theorem, $\Sigma^\ast$ is a set of full surface measure, hence it
is dense. Obviously, $\delta_0(x)$ does depend on $x\in
\Sigma^\ast$ now. It is an easy exercise to check (with a covering
argument using compactness of $\Sigma$) that condition (i) is also
satisfied. \label{lip-ex}
\end{example}

\begin{example}[$W^{2,2}$ surfaces]\rm
If $\Sigma=\partial U$ is locally a graph of a $W^{2,2}$ function
and condition (i) is satisfied, then $\Sigma\in \A$. This follows
from Toro's \cite{Toro-w22} theorem on the existence of
bi-Lipschitz parame\-tri\-zations for such surfaces.
\end{example}

\begin{example} \rm If a compact, connected surface $\Sigma=\partial U$ is locally a
graph of an $AC^2$-function (see J. Mal\'{y}'s paper \cite{maly}
for a definition of absolutely continuous functions of several
variables)  and if (i) is satisfied -- which is a necessary
assumption as graphs of $AC^2$ functions may have cusps -- then
$\Sigma$ is admissible. ($AC^2$ functions are differentiable a.e.
and this implies condition (ii) of Definition~\ref{umbrella}.)
\end{example}

\subsection{The energy and two simple estimates of the inte\-grand}

As mentioned in the introduction, we consider the energy
\begin{equation}
\M_p(\Sigma) := \int_{\Sigma^4} \K^p(T) \, d\mu(T), \qquad
\Sigma \in \A \, ,
\end{equation}
where
$$
\K(T):=\begin{cases} {\displaystyle \frac{V(T)}{A(T)(\diam T)^2}}
& \textnormal{if the vertices of $T$ on $\Sigma$
are not co-planar}\\
0 & \textnormal{otherwise}.
\end{cases}
$$
Here $V(T)$ denotes the volume of $T$ and $A(T)$ the total area,
i.e. the sum of the areas of all four triangular faces of $T$. For
the sake of brevity we write
\begin{equation}
\label{mu} d\mu(\xi,x,y,z) := d\H^2(\xi)\, d\H^2(x)\, d\H^2(y)\,
d\H^2(z).
\end{equation}
If $T=(x_0,x_1,x_2,x_3)$ and one sets $z_i=x_i-x_0$ for $i=1,2,3$,
then  we have
\begin{equation}
\label{1/R} \K(T) = \frac{1}{3}\cdot\frac{|z_3\cdot (z_1\times
z_2)|}{ \Big[|z_1\times z_2|+|z_2\times z_3|+|z_1\times
z_3|+|(z_2-z_1)\times (z_3- z_2)|\Big](\diam T)^2}\,.
\end{equation}
We  will mostly not work with \eqref{1/R} directly. In almost all  proofs in
Sections 3--6, we use iteratively two simple estimates of $\K $
on appropriate classes of tetrahedra.

\begin{lemma}
\label{R-Ved} If $\, T\in \V(\theta,d)$, then
\[
\K(T) > \frac 1{50^2} \theta^4 d^{-1}.
\]
\end{lemma}

\begin{lemma}
\label{R-flat} If $\, T=(x_0,x_1,x_2,x_3)$ is such that
$\Delta(T)=(x_0,x_1,x_2)\in \Spac (\theta,d)$, $x_3 \in B(x_0,2d)$
and $\dist (x_3, \langle x_0,x_1,x_2\rangle ) \ge \kappa d$, then
\[
\K(T) > \frac 1{50^2} \theta^3 \kappa d^{-1}.
\]
\end{lemma}

\medskip\noindent\textbf{Proof of Lemma~\ref{R-Ved}.} Let $T=(x_0,x_1,x_2,x_3)$,
$z_i:=x_i-x_0$ for $i=1,2,3$. Using conditions (ii) and (iii) of
Definition~\ref{Def-Ved}, we obtain $|z_1\times z_2| \ge \theta^2
d^2\sin\theta \ge \frac 2\pi \theta^3 d^2$ and by (iv)
\begin{equation}
\label{num} \left| \frac{z_3\cdot (z_1\times
z_2)}{|z_1\times z_2|} \right| = \dist (x_3, \langle
x_0,x_1,x_2 \rangle) \overset{\textnormal{Def. \ref{Def-Ved}(iv)}}{\ge} \theta d.
\end{equation}
Therefore we can estimate
\begin{eqnarray*}
\K(T) & \overset{\eqref{num}}{\ge} & \frac{1}{3(\diam T)^2}\cdot\frac{\theta d}{
1+
\frac{|z_2\times z_3|}{|z_1\times z_2|}+
\frac{|z_1\times z_3|}{|z_1\times z_2|}+
\frac{|(z_2-z_1)\times (z_3-z_2)|}{|z_1\times z_2|}}\\
& \ge & \frac{1}{3(4d)^2}\cdot\frac{\theta d}{
1+
2\cdot\frac{(2d)^2}{(2/\pi)\theta^3d^2}
+
\frac{(4d)^2}{(2/\pi)\theta^3d^2}}\\
& = &
\frac{\theta^4}{48d[\theta^3+12\pi]}> \frac{\theta^4}{50^2d}.
\end{eqnarray*}
\qed

%
%

\smallskip

The proof of Lemma~\ref{R-flat} is identical. One just replaces
\eqref{num} by $\dist (x_3, \langle x_0,x_1,x_2\rangle ) \ge
\kappa d$.

\section{From energy bounds to uniform Ahl\-fors regularity}\label{sec:3}
\setnumbers

The main result of this section is the following.

\begin{theorem}[Energy bounds imply uniform Ahlfors
regularity] There exists an absolute constant $\alpha>0$ such that
for every $p>8$, every $E>0$ and every $\Sigma\in\A$ with
$\M_p(\Sigma)\le E$ the following holds:

Whenever $x\in \Sigma$, then
\[
\H^2 (B(x,R) \cap \Sigma) \ge \pi R^2/2
\]
for all radii
\begin{equation}
\label{R0Ep} R\,\,\le\,\, R_0(E,p)\,\,:=\,\,
\biggl(\frac{\alpha^{5p}}{E}\biggr)^{\frac{1}{p-8}}\, .
\end{equation}
\label{Thm:3.1}
\end{theorem}

Note that the value of $R_0(E,p)$ depends on $E$ and $p$, but
\emph{not on\/} $\Sigma$ itself, which is by no means obvious.
Even infinitely smooth surfaces can have long `fingers' which
contribute a lot to the diameter but very little to the area. The
point is that fixing an energy bound $E$ we can be sure that
`fingers' cannot appear on $\Sigma$ at a scale smaller than
$R_0(E,p)$. Moreover, a general sequence $\Sigma_i$ of
$C^\infty$-surfaces could in principle gradually form a tip
approaching a cusp singularity as $i\to\infty$ (in fact, it is not
difficult to produce examples of sequences of smooth surfaces with
uniformly bounded area and infinitely many cusp or hair-like
singularities in the limit), whereas this cannot happen according
to Theorem \ref{Thm:3.1} for a sequence of smooth admissible
surfaces with equibounded $\M_p$-energy.

This fact plays a crucial role later on, in the derivation of
uniform estimates for the oscillation of the tangent in
Section~\ref{sec:5}. These estimates in turn allow us to prove in
Section \ref{sec:7} compactness for sequences of surfaces having
equibounded energy.

The scale-invariant limiting case $p=8$ leads to the following
result which can be viewed as a naive counterpart of the
Gau{\ss}--Bonnet theorem for closed surfaces, or the Fenchel
theorem for closed curves: one needs a fixed amount of energy to
`close' the surface. Our estimate of this necessary energy quantum
is by no means sharp; it would be interesting to know the optimal
value of that constant.

\begin{theorem} There exists an absolute constant $\gamma_0>0$
such that $\M_8(\Sigma) > \gamma_0$ for every surface $\Sigma\in
\A$. \label{Thm:3.2}
\end{theorem}

The proof of both theorems relies on a preparatory technical
result which might be of interest on its own, since it allows us
to find for any given admissible surface (no matter how ``crooked''
its shape might look) a good tetrahedron with vertices on the  surface,
i.e. a voluminous tetrahedron in the sense of Definition \ref{Def-Ved}.
This result is completely independent of Menger curvature, but in our
context it will allow us to prove $\M_p$-energy estimates from below.

\begin{theorem}[Good tetrahedra with vertices on
$\Sigma$]
There exist two absolute constants $\alpha,\eta\in
(0,1)$ such that
\begin{equation}
\label{e-a} 1 > 2\eta > 40 \alpha > 0
\end{equation}
with the following property: For every surface $\Sigma\in \A$ and
every $x_0\in \Sigma^\ast$ one can find a positive \emph{stopping
distance} $d_s(x_0) \in (\delta_0 (x_0), \diam \Sigma]$ and a
triple of points $(x_1,x_2,x_3)\in \Sigma\times\Sigma\times\Sigma$
such that

\begin{enumerate}
\renewcommand{\labelenumi}{{\rm (\roman{enumi})}}

\item  $T=(x_0,x_1,x_2,x_3)\in \V(\eta,d_s(x_0))$,
\item whenever $\|T'-T\|\le \alpha d_s(x_0)$,
we have $T'\in \V(\frac \eta 2, \frac 32 d_s(x_0))$.
\end{enumerate}
Moreover, for each $r\in (0,d_s(x_0)]$ there is an affine plane
$H=H(r)$ passing through $x_0$ such that
\begin{equation}
\pi_H (\Sigma\cap B(x_0,r)) \supset H\cap B(x_0,r/\sqrt{2})
\label{bigpi}
\end{equation}
and therefore we have
\begin{equation}
\label{ahlfors-x0} \H^2 (\Sigma \cap B(x_0,r)) \ge \frac \pi 2 r^2
\qquad\mbox{for all $r\in (0,d_s(x_0)]$}.
\end{equation}

\label{goodtetra}
\end{theorem}

The proof of this result is elementary but tedious. We give it in
the next section. We also state one direct corollary of that proof
for sake of further reference.

\begin{proposition}[Large projections and forbidden conical sectors]\label{sectors}
Let $p>8$, $E>0$, and $\partial U = \Sigma\in\A$ with
$\M_p(\Sigma)\le E$. Assume that $R_0=R_0(E,p)$ is given by
\eqref{R0Ep}. For each $x\in \Sigma$ and $r<R_0$ there exists a
plane $H$ passing through $x$ and a unit vector $v\in \Sphere^2$,
$v\perp H$, such that
\begin{eqnarray}
 D:= H\cap B(x,r/\sqrt{2}) & \subset & \pi_H (\Sigma\cap B(x,r))
\label{bigpi-1} \\
\mathrm{int}\, C^+_r(\varphi_0,v) \setminus
B(x,r/2) & \subset & U,
\label{C+in}\\
\mathrm{int}\, C^-_r(\varphi_0,v) \setminus B(x,r/2) & \subset &
\bbbr^3\setminus \overline{U}, \label{C-out}
\end{eqnarray}
where $\varphi_0=\pi/4$.
\end{proposition}

In the remaining part of this section we show how to derive
Theorems~\ref{Thm:3.1} and~\ref{Thm:3.2} from
Theorem~\ref{goodtetra}. We begin with an auxiliary result which
gives an estimate for the infimum of stopping distances considered
in Theorem~\ref{goodtetra}. Note that for $\Sigma$ of class $C^1$,
compact and closed, property (i) below is obvious: we have
$d_s(x_0)> \delta_0(x_0)$, and, as mentioned in
Example~\ref{C1-surf}, in this case one can in fact choose a
positive $\delta_0$ independent of $x_0$.

\begin{proposition}\label{prop:3.5}
Assume that $p > 8$, $\Sigma\in \A$ and $\M_p(\Sigma) < \infty$.
Then
\begin{enumerate}
\renewcommand{\labelenumi}{{\rm (\roman{enumi})}}

\item The stopping distances $d_s(x_0)$
given by Theorem~\ref{goodtetra} have a positive greatest lower
bound,
\[
d(\Sigma):= \inf_{x_0\in \Sigma^\ast} d_s(x_0)>0\, .
\]
\item We have
\begin{equation}
\label{bound-Mp} \M_p(\Sigma) \ge \alpha^{5p} d(\Sigma)^{8-p}\, .
\end{equation}
\end{enumerate}

\end{proposition}
\proof To prove (i), we argue by contradiction. Assume that
$d(\Sigma)=0$ and set
\begin{equation}
\label{eps-contra} \eps:=  \frac 12
\left(\frac{\alpha^{5p}}{K(\Sigma)^4
\M_p(\Sigma)}\right)^{1/(p-8)},
\end{equation}
where $K(\Sigma)$ is the constant from
Definition~\ref{umbrella}~(i). Select $x_0\in \Sigma^\ast$ with
$d_s(x_0)=:d_0<\eps$. Pick $x_1,x_2,x_3$ whose existence is
guaranteed by Theorem~\ref{goodtetra}. Perturbing these points
slightly, by at most $\alpha d_0/2$, we may assume that
\begin{eqnarray}
& & x_i \in \Sigma^\ast, \qquad i=0,1,2,3\, ; \label{c1} \\
& & T=(x_0,x_1,x_2,x_3) \in \V (\eta/2, \frac 32 d_0)\, ; \label{c2} \\
& & \|T'-T\|< \alpha d_0/2 \quad \Rightarrow \quad T'\in \V
(\eta/2, \frac 32 d_0)\, . \label{c3}
\end{eqnarray}
Integrating over all $T'$ close to $T$, we now estimate the energy
as follows:
\begin{eqnarray}
\M_p(\Sigma) & \ge & \int_{\Sigma^4 \cap \B_{\alpha d_0/2}(T)}
\K^p(T')\, d\mu (T') \nonumber \\
& > & \frac 1{K(\Sigma)^4} \left(\frac{\alpha d_0}{2}\right)^8
\left[ \frac 1{50^2} \left(\frac \eta 2\right)^4
(3d_0/2)^{-1}\right]^p \nonumber \\
& > & \frac 1{K(\Sigma)^4} d_0^{8-p}
\left(\frac{\alpha\eta^4}{50^2\cdot 2^6}\right)^p \nonumber \\
& \ge & \frac{\alpha^{5p}}{K(\Sigma)^4} d_0^{8-p}\qquad \mbox{as
$\eta/{20}\ge \alpha\, $.} \label{energy1}
\end{eqnarray}
(We have used Definition~\ref{umbrella}~(i) and Lemma~\ref{R-Ved}
in the second inequality.)
\heikodetail{\xx

\bigskip

In fact, sets of $\H^1\otimes\H^1\otimes\H^1\otimes\H^1$-measure
zero could be neglected for Def. \ref{umbrella}~(i) which might
increase our admissibility  set of surfaces. On the other hand,
this part of the definition is used later on again.

\bigskip

For the very first inequality we note that we only need the
energy locally near $T$, which could be a hint that local
Morrey-space-type theorems could also be true in the surface
case as in the curve case.\xx

\bigskip

}

This gives a contradiction with
\eqref{eps-contra} and the choice of $d_0$, as \eqref{energy1}
implies $d_0>2\eps$.

\medskip

\noindent
(ii)\,
Now we shall show that $d(\Sigma)$ is not only strictly positive,
but has a lower bound depending only on the energy. Fix $\eps>0$
small and pick $x_0\in \Sigma^\ast$ with $d_0:= d_s(x_0) <
(1+\eps)d(\Sigma)$. As in the first part of the proof, take
$x_1,x_2,x_3$ given by Theorem~\ref{goodtetra}. Perturbing these
points slightly, we may assume that \eqref{c1}--\eqref{c3} are
satisfied. Moreover, by \eqref{e-a}
\[
\frac{\alpha d_0} 2 < \frac{d_0}{80} < d(\Sigma) \le d_s(x_i)
\qquad \mbox{for $i=1,2,3$,}
\]
so that by \eqref{ahlfors-x0}
\[
\H^2(\Sigma\cap B(x_i,\alpha d_0/2)) \ge \frac \pi 2
\left(\frac{\alpha d_0}{2}\right)^2, \qquad i=0,1,2,3\, .
\]
Using this information, we again estimate the energy as in
\eqref{energy1}, replacing now the constant $1/K(\Sigma)$ by an
absolute one, $\frac \pi 2$. This yields
%
 \begin{eqnarray*}
 \M_p(\Sigma) & > & \left( \frac{\pi}{2}\right)^4
 \alpha^{5p} d_0^{8-p} \\
 & \ge & \alpha^{5p} (1+\eps)^{8-p} d(\Sigma)^{8-p}\, .
 \end{eqnarray*}
Upon letting $\eps\to 0$, we conclude the whole proof. \hfill
$\Box$

\medskip
\noindent\textbf{Proof of Theorem~\ref{Thm:3.1}:} Inequality
\eqref{bound-Mp} implies that
\[
d(\Sigma) \ge
\biggl(\frac{\alpha^{5p}}{\M_p(\Sigma)}\biggr)^{\frac{1}{p-8}}\, .
\]
Combining this estimate with \eqref{ahlfors-x0}, we see that
\begin{equation}
\label{end-3} \H^2(\Sigma\cap B(x,r))\ge \frac \pi 2 r^2, \qquad r\in
(0, d(\Sigma)]
\end{equation}
holds for all $x\in \Sigma^\ast$. Since $\Sigma^\ast$ is dense in
$\Sigma$, \eqref{end-3} must in fact hold for \emph{all\/} $x\in
\Sigma$.
\heikodetail{

\bigskip

Use density in the following way: for all $\eps >0$
there is $x_\eps\in\Sigma^\ast$ such that
$B(x,r)\supset B(x_\eps,r-\eps)$ which implies
$$
\H^2(\Sigma\cap B(x,r))\ge \H^2(\Sigma\cap B(x_\eps,r-\eps))
\ge \frac{\pi}{2}(r-\eps)^2
$$
and thus the desired estimate upon $\eps\to 0$.

\bigskip

}

\medskip
\noindent\textbf{Proof of Theorem~\ref{Thm:3.2}:} We shall
construct inductively a (possibly finite) sequence  of tetrahedra
with vertices in $\Sigma^\ast$.

Initially, we pick an arbitrary point $x_0=x_0^{1}\in
\Sigma^\ast$. Let $d_1:=d_s(x_0^{(1)})>0$. Use
Theorem~\ref{goodtetra} and density of $\Sigma^\ast$ to find a
tetrahedron
\begin{equation}
T_1 = (x_0^{(1)}, x_1^{(1)}, x_2^{(1)}, x_3^{(1)}) \in \V (\eta/2,
3 d_1 /2) \cap (\Sigma^\ast)^4
\end{equation}
such that
\begin{equation}
\|T'-T_1\| \le \frac{\alpha d_1}{2} \quad \Rightarrow \quad T'\in
\V (\eta/2, 3 d_1 /2).
\end{equation}
Assume that $T_1,T_2, \ldots, T_k$ have been already defined, $T_j
= (x_0^{(j)}, x_1^{(j)}, x_2^{(j)}, x_3^{(j)})$ for $j=1,\ldots,
k$, so that the following properties are satisfied:
\begin{eqnarray}
& & d_{j}:=d_s(x_0^{(j)})< \frac{\alpha d_{j-1}}4, \qquad j=2,\ldots,
k; \label{3.13} \\
& & T_j \in \V (\eta/2, 3 d_j /2) \cap (\Sigma^\ast)^4; \label{3.14} \\
& & \|T-T_j\| \le \frac{\alpha d_j}{2} \quad \Rightarrow \quad
T\in \V (\eta/2, 3 d_j /2); \label{3.15} \\
& & x_0^{(j)} = x_{i(j)}^{(j-1)} \qquad\mbox{for some $i(j)\in
\{1,2,3\}$.\label{3.16} }
\end{eqnarray}
(The last property simply means that $T_j$ and $T_{j-1}$ have one
vertex in common.) Now for $y\in\Sigma$, let
\[
R_\ast(y) = \sup \{r>0 : \H^2 (\Sigma \cap B(y,\varrho))\ge \pi
\varrho^2 /2 \quad\mbox{for all $\varrho\in (0,r]$}\}\, .
\]
We consider the following stopping condition:
\begin{equation}
\label{stop-R} R_\ast(x_i^{(k)}) \ge \frac{\alpha d_k}{4}=:r_k
\qquad\mbox{for all $i\in \{1,2,3\}$}.
\end{equation}
For a fixed value of $k$, there are two cases possible.

\smallskip
\noindent\emph{Case 1. Condition \eqref{stop-R} does hold.} We
then estimate the energy, integrating over small balls centered at
vertices of $T_k$. This yields
\begin{eqnarray*}
\M_8(\Sigma) & \ge & \int_{\Sigma^4 \cap \B_{r_k}(T_k)}
 \K^8(T) \,  d\mu(T) \\
  & > & \left(\frac \pi 2\right)^4 r_k^8 \left[ \frac
 1{50^2} \left(\frac \eta 2 \right)^4 (3d_k/2)^{-1} \right]^8
 \qquad\mbox{by Lemma~\ref{R-Ved}} \\
 & =: & \gamma_0\  >\ 0,
\end{eqnarray*}
where the constant $\gamma_0$ depends \emph{only\/} on the choice
of $\alpha$ and $\eta$ (note that the ratio $r_k/d_k = \alpha/4$
is constant). This is the desired estimate of $\M_8(\Sigma)$.

\smallskip
\noindent\emph{Case 2. Condition \eqref{stop-R} fails.} Then we choose
$i(k)\in \{1,2,3\}$ such that by \eqref{e-a}
\[
R_\ast \bigl(x^{(k)}_{i(k)}\bigr) < r_k=\frac{\alpha d_k}4 < \frac
1{160} d_k.
\]
We set $x_0^{(k+1)}:= x^{(k)}_{i(k)}$ and
$d_{k+1}:=d_s(x_0^{(k+1)})$. The choice of $i(k)$ gives
\begin{equation}\label{d_k-relation}
d_{k+1}
<\alpha d_k/4 < d_k/160.
\end{equation}
Again, we use Theorem~\ref{goodtetra} and
density of $\Sigma^\ast$ to find the next tetrahedron
\[
T_{k+1} = (x_0^{(k+1)}, \, x_1^{(k+1)}, \, x_2^{(k+1)}, \,
x_3^{(k+1)}) \in \V(\eta/2, 3d_{k+1}/2) \cap (\Sigma^\ast)^4
\]
such that \eqref{3.15} is satisfied for $j=k+1$. Thus, we have
increased the length of sequence of tetrahedra satisfying
\eqref{3.13}--\eqref{3.16}.

\medskip

Note that if the stopping condition \eqref{stop-R} is satisfied
for some $k=1,2,\ldots $, then we are done. The only possibility
left to consider is that \eqref{stop-R} fails for each $k$. We
then have an infinite sequence of tetrahedra satisfying
\eqref{3.13}--\eqref{3.16}. To prove that this also gives the
desired result, we shall show later that
\begin{equation}
\textnormal{the sets}\quad\Sigma^4 \cap \B_{r_k}(T_k), \quad k=1,2,\ldots,
\qquad\mbox{are pairwise disjoint.} \label{S-k}
\end{equation}
Assuming \eqref{S-k} for the moment,
we have by Definition~\ref{umbrella}~(i) and
Lemma~\ref{R-Ved}
\begin{eqnarray*}
\M_8(\Sigma) & \ge & \sum_{k=1}^\infty \int_{\Sigma^4 \cap \B_{r_k}(T_k)}
\K^8(T)\, d\mu(T) \\
& > & \sum_{k=1}^\infty \frac{1}{K(\Sigma)^4} r_k^8 \left[ \frac
1{50^2} \left(\frac \eta 2 \right)^4 (3d_k/2)^{-1} \right]^8 \\
&  = & \frac{\gamma_1}{K(\Sigma)^4} \sum_{k=1}^\infty 1 \\
&  = & +\infty,
\end{eqnarray*}
where $\gamma_1$ denotes some constant depending \emph{only\/} on
the choice of $\alpha$ and $\eta$ (again, note that $r_k/d_k =
\alpha/4$ for each $k$).

It remains to prove \eqref{S-k}. Since $T_k \in \V(\eta/2, 3 d_k
/2)$ for each $k$, we have by virtue of Part (i) of
Definition \ref{Def-Ved}
\[
|x_0^{(k+1)} - x_0^{(k)}| = |x_{i(k)}^{(k)} - x_0^{(k)}| \le 3d_k,
\]
so that \eqref{d_k-relation} implies for each $m>k$
\begin{eqnarray}
|x_0^{(m)} - x_0^{(k)}| & \le & 3d_k + 3d_{k+1} + \cdots \nonumber \\
& < & 3d_k (1 + 160^{-1} + 160^{-2} + \cdots) \nonumber \\
& < & 4d_k\, .\label{x-mk}
\end{eqnarray}

For $m=k$ \eqref{x-mk} holds trivially. Also by definition of
$\V(\eta/2, 3 d_{k-1} /2)$ we have
\begin{eqnarray}
\label{xkk-1} |x_0^{(k)} - x_0^{(k-1)}| & = & |x_{i(k-1)}^{(k-1)}
- x_0^{(k-1)}| \\
& \ge &
 \frac \eta 2\frac{3d_{k-1}}2 \nonumber \\
 & > & 15\alpha d_{k-1} \qquad\mbox{as
 $\eta > 20\alpha$.} \nonumber
 \end{eqnarray}

Using \eqref{x-mk}, \eqref{xkk-1}, and the condition $4d_k <
\alpha d_{k-1}$, we obtain
\begin{eqnarray*}
|x_0^{(m)} - x_0^{(k-1)}| & \ge & |x_0^{(k)} - x_0^{(k-1)}| -
|x_0^{(m)} - x_0^{(k)}| \\
& > &  15\alpha d_{k-1} - 4d_k \\
& > & 14\alpha d_{k-1}
\end{eqnarray*}
for each $m\ge k.$
The last inequality readily implies that $B_{r_m}
\bigl(x_0^{(m)}\bigr)$ and $B_{r_{k-1}} \bigl(x_0^{(k-1)}\bigr)$
are disjoint for all $m\ge k$, as
\[
r_m + r_{k-1} = \frac\alpha{4} (d_{m} + d_{k-1}) < \frac{\alpha
d_{k-1}}2.
\]
Thus, the sets $\B_{r_m}(T_m)$ and $\B_{r_{k-1}}(T_{k-1})$ are
disjoint in $\bigl(\rd\bigr)^4$, which proves \eqref{S-k}.

\bigskip

The whole proof of Theorem~\ref{Thm:3.2} is complete now. \hfill
$\Box$

\section{Good tetrahedra: Proof of Theorem~\ref{goodtetra}}
\setnumbers \label{sec:4}

The proof of Theorem~\ref{goodtetra} is lengthy but elementary. It is of
algorithmic nature and, at each of finitely many steps,
requires a case inspection which from a geometric point of view is
not so complicated but nevertheless includes three different cases
(and one of them has to be divided into three further subcases).
The crucial task is to find a triple $(x_1,x_2,x_3)$ such that the
$x_i$'s ($i=0,1,2,3$) satisfy conditions (i) and \eqref{bigpi} of
the theorem. Condition (ii) follows then from simple estimates
based on elementary linear algebra; for sake of completeness, we
present the details of that part in Section~\ref{Sec:3.4}.

\smallskip

Here are a few informal words about the main idea of the proof.

\smallskip

Assume for a while that $\Sigma=\partial U$ is of class $C^1$. To
find a candidate for $x_1$, we look at the surface
$M_\rho=\partial B_\rho \cap C$, where $\rho>0$, $B_\rho$ is
centered at $x_0$, and $C$ is a double cone with vertex $x_0$,
fixed opening angle, and axis given by $n(x_0)$, the normal to
$\Sigma$ at $x_0$. For small $\rho>0$, $x_0$ is the only point of
$\Sigma$ in $C_\rho:=B_\rho \cap C$. (If $\Sigma\in \A$ is not
$C^1$, then the existence of an appropriate cone follows from Part
(ii) of Definition~\ref{umbrella}.)

It is clear that as $\rho$ increases, the growing cone $C_\rho$
must hit $\Sigma$ for some (possibly large) $\rho=\rho_1>0$, at
some $x_1\in \Sigma\setminus\{x_0\}$, $x_1\in M_{\rho_1}$. If the
point of the first hit, $x_1$, is close to the center of one of
the two ``lids'' $\M_{\rho_1}$ of the cone $C_{\rho_1}$, then we
can use the fact that the two components $U^+,U^-$ of
$\mathrm{int\, } C_{\rho_1}$ are on two different sides of
$\Sigma$ to select a voluminous tetrahedron with two of its
vertices at $x_0$ and $x_1$, and all edges $\approx \rho_1$. To
convince yourself that this is indeed plausible, note that there
are many segments perpendicular to $T_{x_0} \Sigma$, with one
endpoint in $U^+$ and the other in $U^-$; each such segment
intersects $\Sigma$ and therefore contains a candidate for one of
the remaining vertices. And, as we shall check later, many of
those candidates are good enough for our purposes.

\begin{wrapfigure}[12]{l}[0cm]{5cm}
\includegraphics*[totalheight=4.6cm]{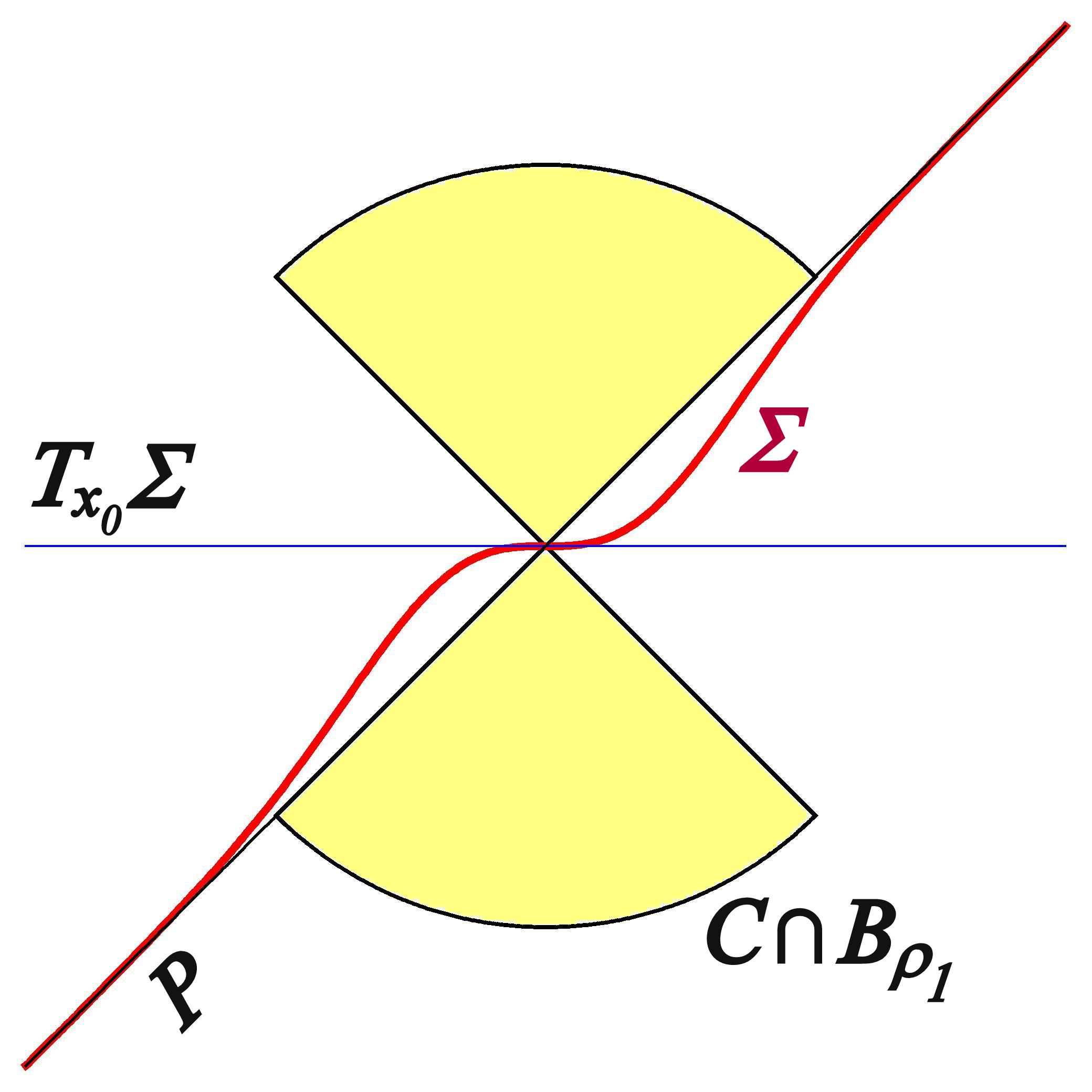}

{\footnotesize

\textbf{Fig. 1.} A little kink determines $T_{x_0}\Sigma$.

}

\end{wrapfigure}

However, it might happen that for this particular intermediate
value of $\rho_1>0$ --- somewhere between $\diam \Sigma$ and the
infinitesimal scale where a smooth $\Sigma$ is very close to the
tangent plane --- most points of $\Sigma\cap B_{\rho_1}$ are very
close to a fixed plane $P$ which might be completely different
from $T_{x_0}\Sigma$, due to a little kink of $\Sigma$ near to
$x_0$. In fact, such a plane might be tangent to $\partial C$, and
$\Sigma\cap B_{\rho_1}$ would look pretty flat at all length
scales $\approx \rho_1$.

If this were the case, then $x_1$ would be located close to the
rim of $C\cap B_{\rho_1}$, and one could not expect to find a good
tetrahedron with vertices $x_i\in \Sigma \cap B_{\rho_1}$ and
edges $\approx \rho_1$. But then, one might rotate $C$ around an
axis contained in $T_{x_0}\Sigma$, away from such a plane $P$, to
a new position $C'$ chosen so that two connected components of
$C'\cap (B_{\rho_1} \setminus B_{\rho_1/2})$ are still on two
different sides of $\Sigma$. One could look for possible vertices
of a good tetrahedron in $C'\cap B_\rho$ for $\rho>\rho_1$,
enlarging the radius $\rho$ until $C'\cap (B_\rho\setminus
B_{\rho_1})$ hits the surface again. This would happen for some
radius $\rho_2>\rho_1$.

It might turn out again that at scales comparable to $\rho_2$
large portions of $\Sigma$ are almost flat, close to a single
fixed plane $P'$ which is tangent to $\partial C'$ so that it is
not at all evident how to indicate a voluminous tetrahedron with
vertices $x_i\in \Sigma \cap B_{\rho_2}$ and edges
$\approx{}\rho_2$. One could try then to iterate the reasoning,
rotating portions of the cones if necessary.

Several steps like that might be needed if, for example, $x_0$
were at the end of a long tip that spirals many times --- in such
cases the points of $\Sigma$ that we hit, enlarging the
consecutive cones, might not convey enough information about the
shape of the surface. We make all this precise (including a
stopping mechanism, a procedure which allows one to select
appropriate rotations at each step of the iteration, and a bound
on the number of steps) in subsection~\ref{meat}, using
Definition~\ref{umbrella}~(ii) to construct the desired cones for
small radii. Before, in subsection~\ref{slanted}, we state two
elementary geometric lemmata which are then used to obtain (i) and
\eqref{bigpi} for various quadruples $(x_0,x_1,x_2,x_3)$.

\medskip

Without loss of generality we suppose throughout
Section~\ref{sec:4} that $x_0=(0,0,0)\in \rd$.

\subsection{Slanted planes and good vertical segments}

\label{slanted}

Suppose that we have a fixed a cone $C=C(\varphi_0,v)$ in $\rd$,
where $v \in \stwo$ and $\varphi_0\in (0,\frac \pi 4]$. We also
fix an auxiliary angle $\varphi_1\in (0,\frac \pi 2]$.

Throughout this subsection, we say that a segment $I$ is
\emph{vertical} (with respect to the cone $C$)
if $I$ is parallel to $v$, i.e., $I=I_{s,v}(z)$
for some $s>0$ and $z\in \rd$. Any plane $P=\langle 0,y_1,y_2
\rangle$ whose unit normal $n$ satisfies $0<|
n\cdot v|<1$ is called \emph{slanted\/}. We say that $I$ is \emph{good
(for $P$)} iff $\dist (I,P)\approx \diam I$, up to constants
depending \emph{only\/} on the angles $\varphi_i$.

We state and prove two  elementary lemmata which give quantitative
estimates of the distance between good vertical segments $I$ and
slanted planes spanned by $0$ and two other points $y_1,y_2$. In
the first lemma both $y_i$ have to be in $C\cap F$, on the same
affine plane $F$ whose normal equals the cone axis of $C$,
i.e. with unit normal $n_F=v$. In the second lemma we keep
one of the $y_i$'s in $C$ and allow the other one to belong to a
portion of $C'$, where $C'$ is a cone congruent to $C$ but rotated
by an angle $\gamma \in (0,\varphi_0/2]$.

To fix the whole setting,  pick a radius $\rho > 0$.  Set
$h=\rho\cos\varphi_0$ and $r=\rho\sin\varphi_0$. Moreover, set
$H:=v^\perp\subset \rd$, and  let $\pi_H\colon \rd\to H$ be the
orthonormal projection onto $H$. Let $\sigma_F$ denote the central
projection from $0$ to the affine plane $F:=H+hv$.

\begin{lemma}[{Slanted planes and good vertical segments, I}]
\label{slanted1} Suppose that $P=\langle 0,y_1,y_2\rangle \subset
\rd$ is span\-ned by $0$ and two other points $ y_1\not= y_2 \,
\in\,  F \cap C_\rho(\varphi_0,v) $ such that there is an angle
$\varphi_1\in (0,\pi)$ such that
\[
\pi> \ang \bigl(\pi_H (y_1), \pi_H(y_2)\bigr) \ge \varphi_1
\qquad\mbox{and}\qquad \pi_H(y_i)\not= 0 \quad\mbox{for $i=1,2$.}
\]
Then, there exists a point $z\in H \cap \partial B_{r}$ such that
\begin{equation}
\dist(I_{h,v}(z), P) \ge c_0 \rho,
\end{equation}
where the constant
\begin{equation}
c_0:=c_0(\varphi_0,\varphi_1) = \frac 12 \left(1-\cos\frac{\varphi_1}2\right)
\sin 2\varphi_0
 > 0\, . \label{def-c0}
\end{equation}
\end{lemma}

\begin{wrapfigure}[15]{l}[0cm]{7.5cm}
\includegraphics*[totalheight=6.3cm]{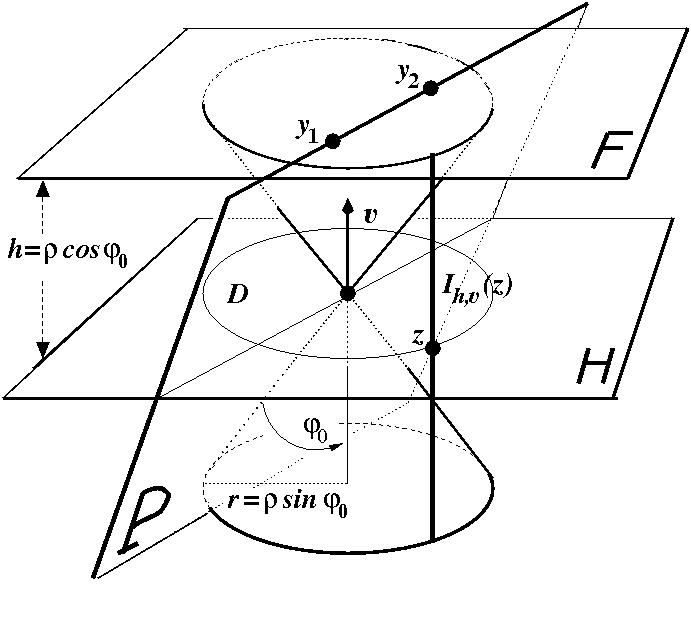}

\vspace*{-2mm}

{\footnotesize

\textbf{Fig. 2.} The setting in Lemma~\ref{slanted1}: a double
cone and three planes $H$, $F$, $P$.

}

\end{wrapfigure}

\medskip\noindent\textbf{Proof.} Let $z_i:=\pi_H(y_i)$ for $i=1,2$.
Consider the 2-dimensional disk
\[
D := H \cap B_r \ \ni\ z_1, z_2.
\]
Let $\gamma:= H\cap \partial B_r$ be the boundary of $D$ in $H$.
We select $z\in \gamma$ such that $z\perp z_2-z_1$ and the segment
$[0,z]$ has a common point with the straight line $l$ which passes
through $z_1$ and $z_2$. By elementary planar geometry, we have
\begin{eqnarray}
 d & := & \dist(z,l)  \nonumber\\
& \ge &  r \left(1-\cos\frac{\varphi_1}2\right) \label{d-est}  \\
&  = & \rho \sin\varphi_0 \left(1-\cos\frac{\varphi_1}2 \right) \,
. \nonumber
\end{eqnarray}
Now, let $\psi $ denote the angle between $v$ and $P$. It is easy
to see that we have $0<\psi<\varphi_0$ since $\varphi_1\in (0,\pi/2]$
and $y_1\not= y_2\in F\cap C_\rho(\varphi_0,v)$. Thus,
\begin{eqnarray*}
\dist(I_{h,v}(z),P) & = & d \cos \psi  \,\, \ge\,\,   d \cos \varphi_0 \\
& \stackrel{\eqref{d-est}}{\ge} & \rho \cos \varphi_0\,
\sin\varphi_0
\left(1-\cos\frac{\varphi_1}2\right) \\
& = & c_0 \rho,
\end{eqnarray*}
where the constant $c_0$ is given by \eqref{def-c0}.
\hfill $\Box$

\begin{lemma}[{Slanted planes and good vertical segments, II}]
\label{slanted2} Let $y_1 \in F\cap B_\rho$, assume $\pi_H
(y_1)\not =0$ and set $u= \pi_H (y_1)/|\pi_H (y_1)|$. Let
$w:=u\times v$ and consider the family of rotations $R_s :=
R(s\varphi_0,w)$, where $s \in [0,\frac 12].$ Then, for any point
\begin{equation}
\label{pos-y2} y_2\in \bigcup_{s\in [0,1/2]} R_s
\biggl(C_\rho(\varphi_0,v)\setminus \mathrm{int}\,
B_{\rho/2}\biggr)\qquad\mbox{such that $  y_2\cdot u < 0
<  y_2\cdot v $}
\end{equation}
there exists a point $z\in H\cap \partial B_r$ such that
$\dist(I_{h,v}(z), \langle  0,y_1,y_2 \rangle ) \ge c_1 \rho.$ One
can take $c_1\equiv c_1(\varphi_0) = \frac{1}{16}\sin 2\varphi_0
>0$.

\end{lemma}

\medskip\noindent\textbf{Proof.} Consider the two-dimensional disk $
D:=F \cap B_\rho$ and its boundary circle $\gamma = F\cap
\partial B_\rho$. Note that the radius of $D$ equals
$r=\rho\sin\varphi_0$. The key point is to observe what the union
of all the central projections $\sigma_F(R_s(D))$, $s\in [0,1/2]$,
looks like. The rest will follow from the previous lemma.

Without loss of generality we assume that $v=(0,0,1)\in \stwo$ and
$y_1=(a,0,h)\in \rd$ for some $a\in (0,r]$. Then $u= \pi_H
(y_1)/|\pi_H (y_1)| = (1,0,0)$ and $w=u\times v = (0,-1,0)$. In
the standard basis of $\rd$ --- which is $(u,-w,v)$ --- the
rotations $R_s=R(s\varphi_0,w)$ are given by
\[
R_s =\left(
\begin{array}{ccc}
\cos s\varphi_0 & 0 & - \sin s\varphi_0 \\
0 & 1 & 0 \\
\sin s\varphi_0 & 0 & \cos s\varphi_0
\end{array}
\right)\, .
\]
\begin{wrapfigure}[16]{l}[0cm]{8.5cm}
\vspace*{-2.5mm}

\includegraphics*[totalheight=6.2cm]{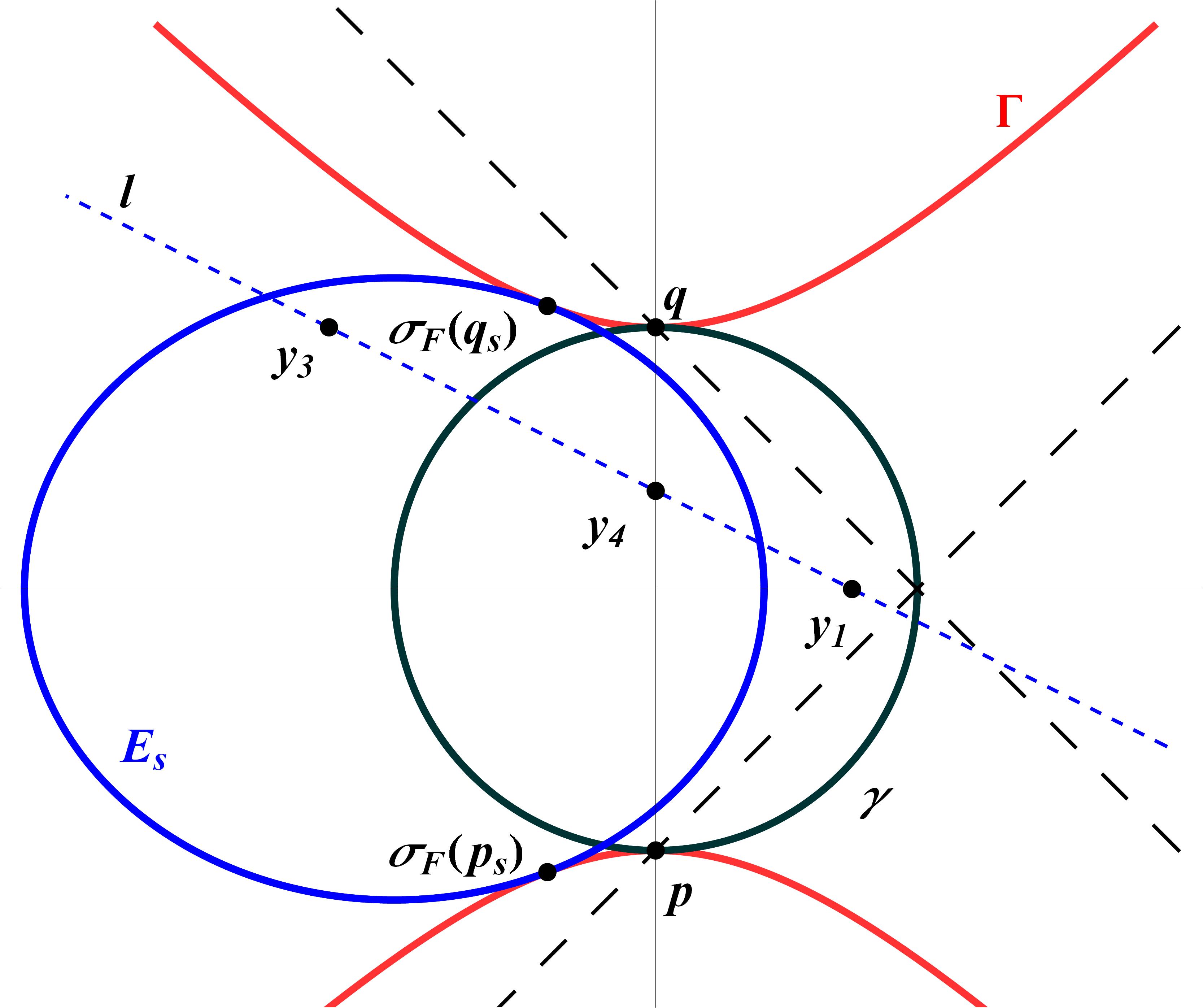}

\vspace*{3mm}

{\footnotesize

\textbf{Fig. 3.} The situation in $F$.

}

\end{wrapfigure}
Now, consider the points $p=(0,-r,h)$ and $q=(0,r,h)$ in $\gamma
\subset F$. Let $p_s=R_s(p)$ and $q_s=R_s(q)$, $s\in [0,\frac
12]$. Since the axis of rotation contains $w$, the angles $\ang
(p_s,w)$ and $\ang (q_s,-w)$ are constant for all $s$
and equal
$\frac \pi 2- \varphi_0$.
Thus, as $s$ goes from $0$ to $\frac
12$, the points $p_s,q_s$ move along arcs of vertical circles on
$\partial C(\frac \pi 2- \varphi_0,w)$. Hence, the central
projections $\sigma_F(p_s)$ and $\sigma_F(q_s)$ trace arcs of two
branches of the hyperbola
\[
\Gamma:=F\cap \partial C(\frac \pi 2- \varphi_0,w)\, .
\]
\heikodetail{

\bigskip

$$
\ang(p,w)=\ang(R_s(p),R_s(w))=\ang(R_s(p),w)=\ang(p_s,w)\quad\Foa s\in [0,1/2].
$$

\bigskip

consult also the drawing  on old p. 6

\bigskip

}%
(In fact, as $s$ goes from $0$ to $\frac 12$, the point
$\sigma_F(R_s(x))$ moves along a hyperbola in $F$ for each $x\in
D$, except the $x$'s that lie on the diameter of $D$ parallel to
$u$.)

Note also that, for each $s\in [0,\frac 12]$, the central
projection $\sigma_F(R_s(D))$ is equal to an ellipse $E_s$ which
is tangent to both arms of $\Gamma$ at $\sigma_F(p_s)$ and
$\sigma_F(q_s)$.

Suppose now that $y_2$ satisfies \eqref{pos-y2}.
Since
\[
\sigma_F \Bigl(R_s\bigl(C_\rho(\varphi_0,v)\setminus
\mathrm{int}\, B_{\rho/2} \bigr)\Bigr) = \sigma_F(R_s(D)),
\]
and the plane $P=\langle 0, y_1,y_2\rangle$ contains the line
through $0$ and $y_2$, we have $y_3:=\sigma_F(y_2)\in P$.
Therefore $P=\langle 0, y_1, y_3\rangle$.

As $ y_2\cdot u <0 <  y_2\cdot v$, the first coordinate of
$y_3=\sigma_F(y_2)$ is negative. Hence, the line $l$ which passes
through $y_3$ and $y_1$ in $F$, and satisfies $l= P\cap F$,
contains a point $y_4\in P\cap F$ on the diameter of $D$ whose
endpoints are $p$ and $q$. Thus, $\langle 0, y_1,
y_2\rangle=P=\langle 0, y_1, y_4\rangle$. If $y_4$ is not in the
center of $D$ (as on the figure above), then the desired claim
follows from the previous Lemma, applied for $P=\langle 0, y_1,
y_4\rangle$ and $\varphi_1=\pi/2$. \heikodetail{

\bigskip

With $\varphi_1=\pi/2$ and $\cos\pi/4=1/\sqrt{2}=\sin\pi/4$ one
infers
$$
c_0=\frac{1}{2}\sin2\varphi_0\cdot \left(1-\frac{1}{\sqrt{2}}\right)
> \frac{1}{16}\sin2\varphi_0.
$$

\bigskip

}
If $y_4={}$
the center of $D$, then the plane $P$ is vertical and one
can take e.g. $z=\pi_H(p)$ to conclude the proof.
In that case one has
$$
\dist(I_{h,v}(z),P)=
r=\rho\sin\varphi_0>\frac{\rho}{16}\sin2\varphi_0=\rho c_1.
$$
\hfill $\Box$

\subsection{Looking for good vertices $x_1,x_2,x_3$: the iteration}

\label{meat}

Throughout this subsection we assume that $0=x_0\in \Sigma^\ast
\subset \Sigma =
\partial U$, where $U\in \rd$ is bounded; $\Sigma$ belongs to the
class $\A$ of all admissible surfaces as defined in Definition \ref{umbrella}.
 Fix $\varphi_0=\frac \pi
4$. Proceeding iteratively, we shall construct four finite
sequences:
\begin{itemize}
\parskip -2pt

\item of compact, connected, centrally symmetric
sets $S_0\subset T_1\subset S_1 \subset T_2\subset S_2\subset \cdots \subset
S_{N-1}\subset T_N\subset S_N\subset
\rd$,

\item of unit vectors $v_0,\ldots,v_N, v_0^\ast,\ldots,v_{N-1}^\ast\in\stwo$
such that
$\ang(v_i,v_{i}^\ast) = \varphi_0/2=\pi/8$ for each $i=0,\ldots, N-1$,

\item of two-dimensional subspaces $H_i=(v_i)^\perp \subset \rd$, $i=0,\ldots,N,$

\item and of radii $\rho_0<\rho_1 < \cdots < \rho_N$, where
$\rho_N=:d_s(x_0)$, so $\rho_N$ will provide the desired stopping
distance for $x_0$ as claimed in Theorem \ref{goodtetra}.
\end{itemize}
These sequences will be shown to satisfy the following properties:
\newcounter{cond}
\newenvironment{A-conditions}{\begin{list}
{{\rm (\Alph{cond})}}{\usecounter{cond}
               \setlength{\labelwidth}{4.5em}
               \setlength{\labelsep}{1em}
               \setlength{\leftmargin}{1.8cm}
               \setlength{\rightmargin}{1em}
}}{\end{list}}

\begin{A-conditions}

\item \textbf{(Diameter of $S_i$ grows geometrically).}\quad
We have $S_i \subset B_{\rho_i}\!\equiv\! B(0,\rho_i)$ and $\diam
S_i = 2 \rho_i$ for $i=0,\ldots,N$. Moreover
\begin{equation}\label{rho-cond}
\rho_{i}> 2\rho_{i-1}\quad\Fo
i=1,\ldots,N.
\end{equation}

\item \textbf{(Large `conical caps' in $S_i$ and $T_i$).}\quad
\begin{equation}\label{conical_caps1}
S_i\setminus B_{\rho_{i-1}}=C_{\rho_i} (\varphi_0,v_i)
\setminus B_{\rho_{i-1}}\quad\Fo i=1,\ldots,N,
\end{equation}
and
\begin{equation}\label{conical_caps2}
T_{i+1}\subset B_{\rho_i}\quad\AND\quad S_i\subset T_{i+1}\quad\Fo i=0,\ldots,N-1.
\end{equation}

\item \textbf{(Relation between $S_i$ and $T_{i+1}$).}\quad
For each $i=0,\ldots,N-1$ there is a unit vector $w_i\perp v_i$ and
a continuous one-parameter family of rotations $R_s^i$ with axis
parallel to $w_i$ and rotation angle $s\varphi_0$, $s\in [0,1/2]$,
such that
\begin{equation}\label{T-Def}
T_{i+1}=S_i\cup\bigcup_{s\in [0,1/2]}R_s^i\Big(C_{\rho_i}(\varphi_0,v_i)\setminus
\INT B_{\rho_{i}/2}\Big).
\end{equation}
%

\item \textbf{($\Sigma$ does not enter the interior of
$S_i$ or $T_{i+1}$).}\quad
\begin{eqnarray}\label{interior_cond1}
\Sigma \cap \INT S_i& =& \emptyset\quad\Fo i=0,\ldots,N,\\
\label{interior_cond2} \Sigma\cap\INT T_{i+1}& =&
\emptyset\quad\Fo i=0,\ldots,N-1.
\end{eqnarray}
Moreover, we have
\begin{equation}
\label{jump} \Sigma \cap
\partial B_{r} \cap C(\varphi_0,v^\ast_{i}) = \emptyset
\qquad\mbox{for $\rho_i\le r \le 2\rho_i$,\quad $i=0,\ldots,N-1$},
\end{equation}
and
\begin{equation}\label{umbrellaposition}
\partial B_t\cap C^+(\varphi_0,v_i)\subset U\quad\AND\quad
\partial B_t\cap C^-(\varphi_0,v_i)\subset \R^3\setminus\overline{U}.
\end{equation}
for all $t\in (\rho_{i-1},\rho_i)$ and $i=1,\ldots,N.$

\item \textbf{(Points of $\Sigma\setminus\{x_0\}$ on $
\partial S_i$).}\quad
The intersection $\Sigma \cap \partial B_{\rho_i} \cap
\partial S_i$ is nonempty for each $i=1,\ldots,N$.

\item \textbf{(Big projections of $B_{\rho_i}\cap \Sigma$ onto $H_i$).}\quad
For $t\in [\rho_{i-1},\rho_i]$, $i=1,\ldots,N$  we have
\begin{equation}\label{bigproj}
\pi_{H_i} (\Sigma\cap B_t) \, \supset\, H_i \cap B_{t\sin
\varphi_0}\, .
\end{equation}
Moreover, for $r_i=\rho_i\sin\varphi_0$, $i=1,\ldots,N$,
\begin{equation}\label{segments}
I_{|z|,v_i}(z),\quad z\in A_i:= H_i \cap (B_{r_i}\setminus\INT
B_{r_i/2}) \quad\textnormal{contains at least one point of
$\Sigma$.}
\end{equation}
\heikodetail{

\bigskip

If we chose to keep the opening angle $\varphi_0$ variable,
then we would have
that
\[
I_{h(z),v_i}(z),\qquad z\in A_i:= H_i \cap (B_{r_i}\setminus\INT
B_{r_i/2})
\]
contains at least one point of $\Sigma$ for the height $h(z):=
|z|\tan((\pi/2)-\varphi_0)$.

\bigskip

}
\end{A-conditions}

\begin{center}
\includegraphics*[totalheight=11.8cm]{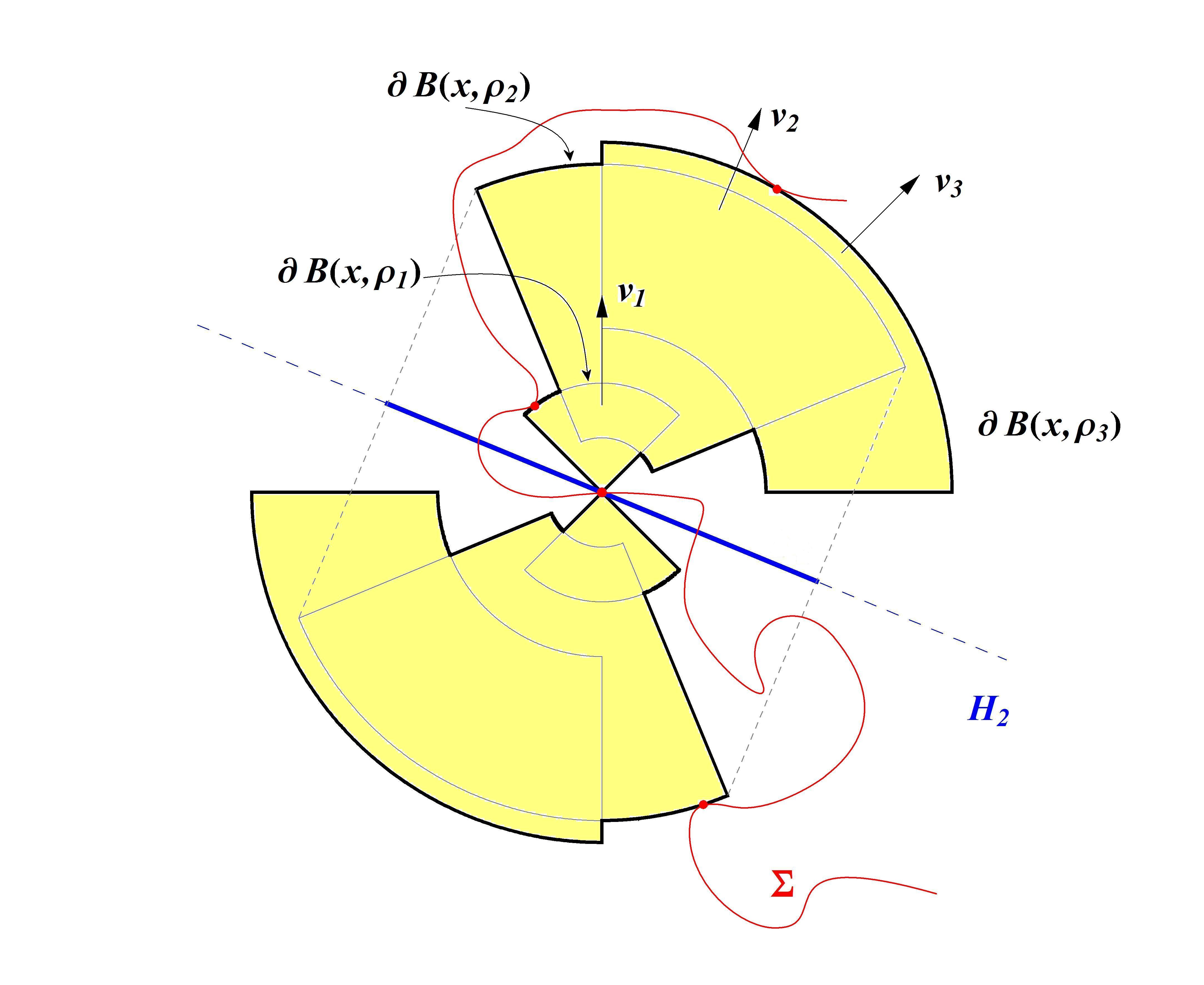}
\end{center}

\vspace*{-7mm}

{\footnotesize

\medskip\noindent\textbf{Fig. 4.} A possible outcome of the
iterative construction. Here, $x=x_0$ is at the center of the
picture and we have $N=3$. The position of the disk $B_{r_2}\cap
H_2$, containing the annulus $A_2$ mentioned in  Condition
\eqref{segments} of (F), is marked with a thick line.

}

\medskip

Once this is achieved, condition (E) implies that
\[
\H^2 (\Sigma \cap B(x_0,r)) \ge \H^2 (D^2(0,r\sin\varphi_0)) = \pi
r^2 /2 \qquad \mbox{for $0<r\le \rho_N=:d_s(x_0)$,}
\]
where $D^2(p,s)$ denotes a planar disk with center $p$ and radius
$s$. We shall also show that it is possible to select $x_j\in
B_{2\rho_N}$ ($j=1,2,3$) with the desired properties listed in
Theorem \ref{goodtetra}.

\medskip \noindent\textbf{Start of the iteration.} We set
$S_0:=\emptyset$ and $T_1:=\emptyset$, $\rho_0:=0$ and $v_0^\ast:=v_1:= v(x_0)$, where
$v(x_0)\in \stwo$ is given by Definition~\ref{umbrella}~(ii).
For $v_0$ we take any unit vector with
the angle condition $\ang(v_0,v_0^\ast)=\ang(v_0,v_1)=\varphi_0/2=\pi/8.$
Then we have $H_0:=(v_0)^\perp$ and $H_1:=(v_1)^\perp.$
Moreover, we use the convention that our closed balls
are defined as
$$
B_r=B(0,r):=\overline{\{y\in\R^3:|y|<r\}}
$$
so that the closed ball $B_0$ of radius zero is the empty set. Notice
that for a complete iteration start we need to define $\rho_1$ and
$S_1$ in order
to check Conditions \eqref{rho-cond} in (A), \eqref{conical_caps1} in (B),
\eqref{interior_cond1} and \eqref{umbrellaposition} for $i=1$,
Condition (E), and \eqref{bigproj} and \eqref{segments} constituting
Condition (F).
All the other conditions
within the list (A)--(D) are immediate for $i=0.$
\heikodetail{

\bigskip

$$
S_0=\emptyset\subset B_{\rho_0}=B_{0}=\emptyset\quad\textnormal{by our
convention about balls},
$$
and $\diam S_0=\diam\emptyset=0=2\rho_0$, which gives the first part of
(A) for $i=0.$

$T_1=\emptyset\subset \emptyset=B_{\rho_0}$ and $S_0=\emptyset
\subset T_1$, which proves \eqref{conical_caps2} in (B).

For $w_0$ we can take any unit vector perpendicular to $v_0$ ,
since no matter what rotation $R_s^i$ will result, the right-hand
side of \eqref{T-Def} will be the empty set, since
$C_{\rho_0}\subset B_{\rho_0}= \emptyset.$ The left-hand side is
by definition empty, which proves \eqref{T-Def} for $i=0.$

\eqref{interior_cond1} and \eqref{interior_cond2}  for $i=0$
are immediate since $T_1$ and $S_0$ are empty.
Also \eqref{jump} holds trivially for $i=0$, since $\partial B_r$
is empty for any $r\in [\rho_0,2\rho_0]=[0,0]$, i.e. for $r=0.$

\bigskip

}
We set
\begin{equation}
\label{K1} K^1_t := C_t(\varphi_0,v_1)\,.
\end{equation}
With  growing radii $t$ the sets $K^1_t$ describe larger
and larger double cones with constant opening angle $2\varphi_0=\pi/2$
and fixed axis $v_1$. Now we  define
\begin{equation}\label{stop_1}
\rho_1:=\inf\{t>\rho_0=0:\Sigma\cap K^1_t\cap\partial B_t\not=\emptyset\},
\end{equation}
and notice that by definition of the set $\A$ of admissible
surfaces (see Definition \ref{umbrella} (ii)) one has
$\rho_1>\delta_0(x_0)>0=2\rho_0$, which takes care of
\eqref{rho-cond} in Condition (A) for $i=1$. Set
$S_1:=K^1_{\rho_1}$, then we see that
$S_1=C_{\rho_1}(\varphi_0,v_1) \subset B_{\rho_1}$ with $\diam
S_1=2\rho_1$, so all properties of (A) hold for $i=1.$ Moreover,
$$
S_1=C_{\rho_1}(\varphi_0,v_1)=C_{\rho_1}(\varphi_0,v_1)
\setminus  B_{\rho_0},
$$
since $B_{\rho_0}=B_0=\emptyset$, thus \eqref{conical_caps1} in (B)
holds for $i=1$.
The definition of $\rho_1>0$ (see \eqref{stop_1}) implies
\eqref{interior_cond1} for $i=1$, notice that $\INT S_1$
is the union of two disjoint open cones centrally symmetric
to but not containing $x_0=0\in\Sigma$. For the proof of
\eqref{umbrellaposition} for $i=1$ we observe that for each
$t\in (0,\rho_1)$ we have by Definition of $\rho_1$ that
\begin{equation}\label{helpful}
B_t\cap C^+(\varphi_0,v_1)\subset U\cup\{0\}\quad\AND\quad
B_t\cap C^-(\varphi_0,v_1)\subset (\R^3\setminus \overline{U})\cup\{0\},
\end{equation}
which is even stronger than \eqref{umbrellaposition}. Condition
(E) holds for $i=1$, too, by definition of $\rho_1$ and the fact
that $\Sigma $ is a closed set. For $i=1$ we will prove
\eqref{segments} even for {\it all} $z\in D_1:= H_1\cap B_{r_1}$,
which would immediately imply \eqref{bigproj} of Condition
(F).\footnote{Alternatively, one could look for $t\in (0,\rho_1)$
at the (longer) vertical segments $I_{\psi(t),v_1}(z)$, $\psi(t):=
\sqrt{t^2-|z|^2}$, whose endpoints are contained in $\partial
B_t\cap C^+(\varphi_0,v_1),$ and in $\partial B_t\cap
C^-(\varphi_0,v_1),$ respectively, use \eqref{umbrellaposition}
for $i=1$ as proved just before, to conclude that
$I_{\psi(t),v_1}$ intersects $\Sigma$ for each $t\in (0,\rho_1).$
This proves \eqref{bigproj} for $t\in (0,\rho_1)$, the statement
for $t=0=\rho_0$ is trivial, and for $t=\rho_1$ use continuity,
and the fact that $\Sigma$ is a closed set. This is actually the
argument we repeat in the induction step $j\mapsto j+1$ later on,
since there we have less explicit information about $S_j$.}{} From
\eqref{helpful} we also infer that  every segment
$I_{|z|,v_1}(z)$, for $z\in  H_1\cap (B_{r_1}\setminus \{0\}) $
with $|z|<r_1$, has one endpoint in $U$, and the other in
$\R^3\setminus \overline{U}$, which implies that $I_{|z|,v_1}(z)$
intersects the closed surface $\Sigma$ in at least one point for
these $z$. For $z=0=x_0\in\Sigma$ this is trivially also true, and
for
 $z\in D_1$ with $|z|=r_1$ we approximate $z_k\to z$ as $k\to\infty$
 with points $z_k\in D_1$ and $|z_k|<r_1$ to find a sequence
 $\xi_k\in\Sigma\cap
 I_{|z_k|,v_1}(z_k)$ which converges to some surface point $\xi\in
 \Sigma\cap I_{|z|,v_1}(z). $ This  completes
the proof of \eqref{segments} even for all $z\in H_1\cap B_{r_1}$
and hence of (F) for $i=1$.

To summarize this first step,  we have defined the sets
$S_0\subset T_1\subset S_1 \subset\R^3$, the unit vectors
$v_0,v_1,v_0^\ast\in\S^2$ with $\ang(v_0,v_0^\ast)=\varphi_0/2,$
and the corresponding subspaces $H_0=(v_0)^\perp$, and
$H_1=(v_1)^\perp$, and finally radii $\rho_0=0<\rho_1$ without
having made the decision if $N=1$ or $N>1$. In addition we have
now proved the first two items in Condition (A) for $i=0,1$, and
\eqref{rho-cond} for $i=1.$ Moreover,  we have verified
\eqref{interior_cond1} for $i=0,1$, and all other statements in
the list of properties (B)--(F) are established for the respective
smallest index $i$. Note, however, that we have not defined
$v_1^\ast$ yet.

\medskip \noindent\textbf{Stopping criteria and the iteration step.}
For the decision  whether to stop the iteration or to  continue it  with
step number  $j+1$ for $j\ge 1$,  we may now assume
that the sets
$$
S_0\subset T_1\subset S_1\subset T_2\subset S_2\subset \cdots \subset
T_{j}\subset
S_{j}\subset\R^3,
$$
and unit vectors $v_0,\ldots,v_{j}$, $v_0^\ast,\ldots,v_{j-1}^\ast$
with $\ang(v_i,v_i^\ast)=\varphi_0/2$ for $i=0,\ldots,j-1$,
are defined. We also have at this point a sequence of radii
$
\rho_0=0<\rho_1 < \cdots < \rho_{j}
$
satisfying \eqref{rho-cond} for $i=1,\ldots,j.$ The
first two conditions in (A) may be assumed to hold for $i=0,\ldots,
j$. In (B) we may suppose \eqref{conical_caps1} for $i=1,\ldots,j$,
in contrast to  \eqref{conical_caps2} which holds
for $i=0,\ldots,j-1.$ Similarly,
we may now work with \eqref{T-Def} in (C),  \eqref{interior_cond2} and
\eqref{jump} in (D) for all
$i=0,\ldots, j-1,$ whereas we may use \eqref{interior_cond1} in (D) for
$i=0,\dots,j$, \eqref{umbrellaposition}, Condition (E), and
\eqref{bigproj} and \eqref{segments} in (F)
now for $i=1,\ldots,j.$

\medskip

Now we are going to study the various
geometric situations that allow us to  stop the
iteration here, in which case we set
$N:=j$, $d_s(x_0):=\rho_j=\rho_N$, so that \eqref{bigpi} and
\eqref{ahlfors-x0} stated in Theorem \ref{goodtetra} can be extracted
for $H:=H_j$ directly from Condition (F). Indeed, \eqref{bigproj}
for $t:=\rho_j=\rho_N$ yields \eqref{bigpi} since $\varphi_0=\pi/4.$
How to find the remaining vertices $x_1,x_2,x_3$ such that
Statement (i) of Theorem \ref{goodtetra}  holds for the tetrahedron
$T=(x_0,x_1,x_2,x_3)$ will be explained later in detail for each
case in which  we stop the iteration.
Moreover, we will convince ourselves that the only case
in which the iteration cannot be stopped, can happen only finitely
many times.  But each time this happens we have to define unit vectors
$v_j^\ast$,
$v_{j+1}\in\S^2$, with $\ang(v_j,v_j^\ast)=\varphi_0/2$, and $H_{j+1}:=(v_{j+1})^\perp$,
a new radius $\rho_{j+1}, $ new sets $T_{j+1}\subset S_{j+1}$
containing $S_j$,
 and then check all the properties listed in (A)--(F).

\bigskip

The different geometric situations depend on how the surface hits
the ``roof'' of the current centrally symmetric set $S_j$, that
is, where the points of the nonempty intersection in Condition (E)
lie:
\begin{description}
\item [Case 1. (Central hit.)] By this we
mean that $\Sigma \cap \partial B_{\rho_j} \cap C(\frac 34
\varphi_0, v_j)$ is nonempty.

\item [Case 2. (No central hit  but nice distribution of
intersection points.)
] By this we mean that Case~1 does
not hold but there exist two different points $x_1,x_2\in\Sigma \cap
\partial B_{\rho_j} \cap C(\varphi_0, v_j)$ such that
\begin{equation}
\label{wide} \ang \bigl(\pi_{H_j}(\sigma(x_1)),
\pi_{H_j}(\sigma(x_2))\bigr) \ge \frac \pi 3\, ,
\end{equation}
where $\pi_{H_j}$ denotes the orthogonal projection onto the
current plane $H_j=(v_j)^\perp.$
\end{description}
In Cases 1 and 2, we can find triples of points $(x_1,x_2,x_3)$
with all the desired properties and stop the iteration right away.
Below, in paragraphs~\ref{C1} and~\ref{C2}, we indicate how to
select the $x_i$'s in each of these cases, and present the
necessary estimates.

If neither Case 1 nor Case 2 occurs, then we have to deal with
\begin{description}
\item [Case 3. (Antipodal position.)] $\Sigma \cap \partial
B_{\rho_j} \cap C(\frac 34 \varphi_0, v_j)$ is empty and for \emph{any
two different points\,} $x_1,x_2\in \Sigma \cap \partial
B_{\rho_j} \cap C(\varphi_0, v_j)$ we have
\begin{equation}
\label{antipo} \ang \bigl(\pi_{H_j}(\sigma(x_1)),
\pi_{H_j}(\sigma(x_2))\bigr) < \frac \pi 3\, .
\end{equation}
\end{description}
(Intuitively, Case~3 corresponds to the situation alluded to in
the introduction to Section 4: at this stage we have to take into
account the possibility that most points of $\Sigma\cap
B_{\rho_j}$ are close to some fixed 2-plane containing the segment
with endpoints $x_0,x_1$.) Now this third case is more
complicated, we will distinguish three further subcases, of which
two will allow us to stop the iteration here. Only the third
subcase will force us to continue the iteration.

To make this precise, let us fix some point
$x_1 \in \Sigma \cap \partial B_{\rho_j} \cap
C(\varphi_0, v_j)$. Such a point does exist according to Condition
(E).
Set
$u_j:=\pi_{H_j}(x_1)/|\pi_{H_j}(x_1)| $ and
 $w_j:=u_j\times v_j$, and consider the
family of rotations
\begin{equation}\label{rotations}
R^j_s := R(s\varphi_0,w_j), \qquad s \in [0,4]\,.
\end{equation}
Consider the union of rotated conical caps
\begin{equation}\label{Gsets}
G^j_t:= \bigcup_{0\le s \le t} R_s(C_{\rho_j}(\varphi_0,v_j) \setminus
\mathrm{int}\,\, B_{\rho_j/2}
), \qquad t\in
[0,\textstyle\frac 12].
\end{equation}
Let
\begin{equation}
\label{stop-Rs} t_0:=\sup \{t\in [0,{\textstyle\frac 12}] \,
\colon \, G_t \cap \bigl(\Sigma \setminus S_j\bigr)
=\emptyset\}.
\end{equation}
(Intuitively: we rotate the conical cap  ``away from the
intersection $\Sigma\cap \partial B_{\rho_j} \cap C(\varphi_0,
v_j)$'' and look for new points of $\Sigma$ in the rotated set.)
There are now three subcases possible. To describe them, let
$v^\ast_j:= R_{1/2}(v_j)$ (this will be the new $v_{j+1}$ in the
third subcase).

\begin{description}

\item[Subcase 3~(a).] $G^j_{t_0} \cap (\Sigma \setminus
S_j) \not= \emptyset$. Then $j=N$; we stop the iteration
and select $x_2$ and $x_3$, the remaining vertices of a good
tetrahedron, using Lemma~\ref{slanted2} to obtain the desired
estimates; see subsection \ref{C3} for the computations.

\end{description}

Intuitively, Subcase 3~(a) corresponds to the situation where
we initially suspect that the surface might be similar to the one
with a little kink (see Fig. 1 at the beginning of Section~4).
Condition \eqref{antipo} alone does not exclude this -- but here,
rotating a portion of the cone slightly, we find new points of
$\Sigma$ and detect that $\Sigma$ is not flat at scale
$\rho_j$.

\begin{description}
\item[Subcase 3~(b).] We have $t_0= 1/2$ and
$G^j_{1/2} \cap (\Sigma \setminus S_j) = \emptyset$.
However,
\begin{equation}
\label{forward}
\Sigma \cap \left(C_{2\rho_j} (\varphi_0,
v_j^\ast) \setminus C_{\rho_j} (\varphi_0, v_j^\ast)\right) \not=
\emptyset.
\end{equation}
Again, $j=N$; we stop the iteration and select $x_2$ and $x_3$.
For details, see subsection \ref{C3}.
\end{description}

Informally: here we rotate a portion of the cone slightly and do
not find new points of $\Sigma$. However, there are other points
of the surface at comparable distances, again allowing us to
exclude the possibility that $\Sigma$ is close to being flat at
scale $\rho_j$.

\begin{description}
\item[Subcase 3~(c).] We have $t_0=1/2$ and
\begin{equation}\label{Gempty}
G^j_{1/2} \cap (\Sigma \setminus
S_j) = \emptyset.
\end{equation}
Moreover, \eqref{forward} is violated,
i.e.,
\begin{equation}
\label{3c} \Sigma \cap \left(C_{2\rho_j} (\varphi_0, v_j^\ast)
\setminus C_{\rho_j} (\varphi_0, v_j^\ast)\right) = \emptyset.
\end{equation}
\end{description}
If this is the case, then  we are unable to exclude the
possibility that (most of) $\Sigma$ is nearly flat at the given
scale, and the iteration goes on. We set $T_{j+1} := S_j \cup
G^j_{1/2}$, $v_{j+1}:=v_j^\ast=R_{1/2}(v_j)$,
$H_{j+1}:=(v_{j+1})^\perp$, and
\begin{equation}\label{Kj+1}
K^{j+1}_t:=C_t(\varphi_0,v_{j+1}),
\end{equation}
and define
\begin{equation}\label{stopj+1}
\rho_{j+1}:=\inf\{t>\rho_{j}:\Sigma\cap K_t^{j+1}\cap\partial B_t\not=\emptyset\}.
\end{equation}
Notice that Condition \eqref{3c} in the context of this subcase guarantees
that $\rho_{j+1} > 2\rho_j$ which verifies \eqref{rho-cond} in Condition
(A) for $i=j+1.$ Now we define
\begin{equation}\label{Sj+1}
S_{j+1}:= T_{j+1}\cup (K^{j+1}_{\rho_{j+1}}\setminus \INT
B_{\rho_j}),
\end{equation}
and
check
that Conditions (A)--(F) are satisfied.  Indeed, $S_{j+1}\subset
S_j\cup K^{j+1}_{\rho_{j+1}}\subset B_{\rho_j}\cup B_{\rho_{j+1}}$
by Condition (A) for $i=j,$ which implies that (A) holds for $i=j+1$
as well. Next,
$$
S_{j+1}
\setminus B_{\rho_j}=K^{j+1}_{\rho_{j+1}}\setminus B_{\rho_j}=C_{\rho_{j+1}}(
\varphi_0,v_{j+1})\setminus B_{\rho_j},
$$
since $S_j\subset B_{\rho_j}$ by Condition (A) for $i=j$. Hence
\eqref{conical_caps1} holds for $i=j+1$. As
$G^j_t\subset B_{\rho_j}$ for all $t\in [0,1/2]$ we have
$T_{j+1}\subset S_j\cup B_{\rho_j}\subset B_{\rho_j}$ because
of Condition (A) for $i=j.$ The second item in \eqref{conical_caps2}
is a direct consequence of the definition of $T_{j+1}$, whence
\eqref{conical_caps2} holds for $i=j$. Condition (C) holds also
for $i=j$
by definition of $T_{j+1}$. Using \eqref{interior_cond1} for $i=j$,
\eqref{Gempty},
and the definition of $\rho_{j+1}>2\rho_{j}$ in \eqref{stopj+1}
we infer that \eqref{interior_cond1} holds for $i=j+1$, and \eqref{interior_cond2} for $i=j.$
Relation \eqref{jump} for each $r\in (\rho_j,\rho_{j+1}]$
is an immediate consequence  of \eqref{3c}. For $r=\rho_j,$ however,
we have to use \eqref{Gempty} in combination with the fact that
all surface points in $\Sigma\cap\partial B_{\rho_j}\cap C(\varphi_0,v_j)$
are in antipodal position described by \eqref{antipo}, so that
$\Sigma\cap\partial B_{\rho_j}\cap C(\varphi_0,v_j^\ast)=\emptyset.$

Now we turn to the proof of \eqref{umbrellaposition} for $i=j+1.$
The definition \eqref{stopj+1} of $\rho_{j+1}$ implies that
\begin{equation}\label{alternatives}
\partial B_t\cap C^+(\varphi_0,v_{i+1})\subset U,\quad\OR\quad
\partial B_t\cap C^+(\varphi_0,v_{i+1})\subset \R^3\setminus\overline{U}
\end{equation}
for all $t\in (\rho_j,\rho_{j+1}).$ Now \eqref{Gempty} together with
\eqref{umbrellaposition} implies that
$$
\partial B_{\rho_j}\cap C^+(\varphi_0,v_{i+1})\subset U,
$$
which excludes the second alternative in \eqref{alternatives}.
Condition (E) holds for $i=j+1$ by the definition of $\rho_{j+1}$
and the fact that $\Sigma$ is a closed set. For the proof of
\eqref{bigproj} for $i=j+1$ we look
 for $t\in (\rho_j,\rho_{j+1})$
 at the  vertical segments $I_{\psi(t),v_{j+1}}(z)$, $\psi(t):=
 \sqrt{t^2-|z|^2}$, $z\in B_{t\sin\varphi_0}\cap H_{j+1}$. The  endpoints
 of these segments
lie in $\partial B_t\cap C^+(\varphi_0,v_{j+1}),$ and
 in $\partial B_t\cap C^-(\varphi_0,v_{j+1}),$ respectively.
 Now we use \eqref{umbrellaposition}
 for $i=j+1$  to conclude that $I_{\psi(t),v_{j+1}}(z)$
 intersects $\Sigma$ for each $t\in (\rho_j,\rho_{j+1}).$ This proves
 \eqref{bigproj} for $t\in (\rho_j,\rho_{j+1})$.
For  $t=\rho_j$ and $t=\rho_{j+1}$ use continuity, and the fact that $\Sigma$
 is a closed set.
 Finally, to prove \eqref{segments} for $i=j+1$ note that
 \eqref{rho-cond} together with \eqref{umbrellaposition} for $i=j+1$
 imply that the two endpoints of the
 vertical segments $I_{|z|,v_{j+1}}$ for $z\in
 A_{j+1}$ lie in the different open connected components $U$ and
 $\R^3\setminus\overline{U}.$ This suffices to conclude that
 these segments intersect $\Sigma $, which finishes the proof
 of all conditions in the list (A)--(F) in the iteration step.

Since we have established Condition (E) in the iteration step and
\eqref{rho-cond} holds, too, we can deduce that Subcase 3 (c) can
happen only finitely many times, depending on the position $x_0$
on $\Sigma$ and on the shape and size of $\Sigma$:
$$
 \diam\Sigma\ge\rho_i>2\rho_{i-1}>\cdots > 2^{i-1}\rho_1>2^{i-1}\delta_0(x_0),
$$
whence  the maximal number of iteration steps is bounded by $$
 1+\log(\diam\Sigma/\delta_0(x_0))/\log 2\, .
$$

This concludes the Subcase 3 (c). Now we have to analyze the
geometric situation in the remaining Cases 1, 2, and 3~(a)
and~(b), to extract surface points $x_1,x_2,x_3$, so that the
selected tetrahedron $T=(x_0,x_1,x_2,x_3)$ (with $x_0=0$)
satisfies Part (i) of Theorem \ref{goodtetra}. Part (ii) then
follows from an easy perturbation argument; see Corollary
\ref{cor:4.4}.

\smallskip

\subsubsection{Case 1 (Central hit): the details}

\label{C1}

We fix a point $x_1\in \Sigma \cap \partial B_{\rho_j}$ such that
\[
x_1\cdot v_j  = \pm \rho_j \cos \gamma_1 , \qquad 0\le
\gamma_1 \le \frac 34 \varphi_0,
\]
and we are going to  select suitable points  $x_2,x_3\in\Sigma\cap
B_{\rho_j}$ so that Condition (i) of Theorem \ref{goodtetra} is
satisfied. This will justify our decision to stop the iteration by
having set $N:=j$ and $d_s(x_0):=\rho_j=\rho_N.$ Without loss of
generality, rotating the coordinate system if necessary, let us
suppose that $v_j=(0,0,1)\in \rd$ and $\pi_{H_j}(x_1)\in H_j$ is
equidistant from $z_1:=(0,r_j,0)$ and $z_2:=(0,-r_j,0)$, where we
recall from Condition (F) for $i=j$ that
$r_j=\rho_j\sin\varphi_0.$ (In other words, we assume w.l.o.g.
that the second coordinate of $x_1$ is zero.)

\smallskip

Condition \eqref{segments} in (F) for $i=j$ guarantees the existence of a
point  $x_2 \in \Sigma \cap I_{h_j,v_j}(z_2)$, where $h_j=\cos\varphi_0=r_j$.
Now  let
$P:=\langle 0, x_1, x_2\rangle$. Then $\pi_{H_j}(x_2) \perp x_1$
and we have
\[
\rho_j |x_2|\, |\cos \ang(x_1,x_2)| = |x_1\cdot x_2| =
|x_1\cdot (x_2 - \pi_{H_j}(x_2))|,
\]
which yields
\[
|\cos \ang(x_1,x_2)| \le \frac{|x_2 - \pi_{H_j}(x_2)|}{|x_2|} \le
\sin (\textstyle{\frac\pi 2} - \varphi_0) = 1/\sqrt{2}\, .
\]
Thus, Definition~\ref{Def-Ved}~(iii) is satisfied for $x_0=0,x_1,x_2,$
for every  $\theta\le\pi/4$.
To select $x_3$, we consider two subcases.

\medskip\noindent\textbf{Subcase 1~(a).}
If the points $z_1,z_2$ and $\pi_{H_j}(x_1)$ are collinear, then
we simply have $P=\langle 0,x_1,z_2\rangle$. We then use (F) for $i=j$
to
select $x_3\in \Sigma \cap I_{h_j,v_j}(z_3)$, where $z_3:=
(r_j,0,0)$ belongs to the two-dimensional
disk $D_j:=D^2(0,r_j)$ in $H_j$. Thus,
\[
\rho_j\sin\varphi_0 \le |x_k-x_i|\le 2\rho_j \qquad \mbox{for
$k\not=i$, $k,i=0,1,2,3$,}
\]
which establishes Conditions (i) and (ii) of Definition \ref{Def-Ved}
for $d:=d_s(x_0)=\rho_j$ and any $\theta\le\sin\varphi_0=1/\sqrt{2}$
Finally, $\dist (x_3, P) = r_j = \rho_j\sin\varphi_0$,
and this takes care of Part (iv)  of Definition \ref{Def-Ved} so that
$T=(x_0,x_1,x_2,x_3)\in\V(\eta,d_s(x_0)$  for any $\eta<1/\sqrt{2}$,
i.e.
in this subcase Part (i) of Theorem \ref{goodtetra} is satisfied
for any $\eta<1/2.$

\medskip\noindent\textbf{Subcase 1 (b).} If the points $z_1,z_2$
and $\pi_{H_j}(x_1)$ are non-collinear, then we consider the line
segment $J:= F_j \cap B_{\rho_j} \cap P\,$ contained in the affine
plane $F_j:= H_j + h_j v_j$. Since $x_1 \in C(\frac 34 \varphi_0,
v_j)$ and $y_1:=\sigma_{F_j}(x_1)\in J$, it is easy to check that,
no matter where $x_2$ has been chosen, $J$ (and $P$) contains
points $y_2\in F_j$ such that
\[
\ang(\pi_{H_j}(y_2), \pi_{H_j} (y_1) ) \ge \arccos\left(
\cot\varphi_0 \tan \frac 34  \varphi_0\right) >  \frac \pi 5.
\]
\heikodetail{

\bigskip

Your $\sqrt{2}$ was OK, I made a correction while re-reading that
proof, trying to follow the steps of someone who wants to do all
computations precisely. --- P. I also relied on a calculator.

For the $\sqrt{2}$-factor look at the following calculation, and
for the strict  inequality above  I confess to rely on a
calculator.

Here is my computation ---P.:

$$
\cos\ang(\pi(y_2),\pi(y_1))=\frac{x}{r_j} \le\frac{h_j \tan \frac
34 \varphi_0 }{\rho_j\sin\varphi_0} \sin \frac 34
\varphi_0=\frac{\cos\varphi_0 \tan\frac 34
\varphi_0}{\sin\varphi_0}= \cot \varphi_0 \tan\frac 34 \varphi_0.
$$
(for $x$ see image on my old p. 16)
\bigskip

}
Therefore, we may apply Lemma~\ref{slanted1} with
$\varphi_0=\frac \pi 4$ and $\varphi_1:=\pi/5$  to select a point $x_3\in
\Sigma$ on one of the vertical segments $I_{h_j,v_j}(z)$, $z\in
\gamma_j:={}$ the boundary of $D_j$ in $H_j$, so that
\[
\eta \rho_j < \dist(x_3,P)\quad\AND\quad \eta\rho_j  <
|x_k-x_i|\le 2\rho_j \qquad
\mbox{for $k\not=i$, $k,i=0,1,2,3$.}
\]
where $\eta:=1/100< \pi/200\le \frac 12
(1-\cos \frac \pi {10}) = c_0(\pi/4,\pi/5)
$ (and we used $1-\cos x\ge x^2/\pi$, $x\in
[0,\frac\pi 2]$, for the first inequality).
\heikodetail{

\bigskip

We have
\begin{eqnarray*}
|x_1-x_0| & = & \rho_j>\eta\rho_j,\\
|x_3-x_i| & \ge & \dist(P,x_3)>\eta\rho_j,\\
|x_2-x_0| & \ge & r_j=\rho_j\sin\varphi_0=\frac{1}{\sqrt{2}}\rho_j>
\eta\rho_j,\\
|x_2-x_1| & \ge & |z_2-z_1|\ge
|z_2|=r_j=\rho_j\sin\varphi_0=\frac{1}{\sqrt{2}}\rho_j>\eta\rho_j,
\end{eqnarray*}
by choice of $z_2$.

\bigskip

}
This verifies Conditions (i), (ii), and (iv) of Definition
\ref{Def-Ved}  for each $\theta\le\eta=1/100$, and we have seen
before that Part (iii) of that Definition holds for
all $\theta\le \pi/4.$ Hence Part (i) of Theorem \ref{goodtetra}
is also satisfied for $\eta:=1/100$ in this subcase, which completes
our considerations for Case 1.

\subsubsection{Case 2 (No central hit but nice distribution
of intersection points): the details}

\label{C2}

\textbf{The setting.} As in Case 1, we have stopped the iteration,
set $N:=j$, $d_s(x_0):=\rho_j=\rho_N$.
Let $H_j =(v_j)^\perp$ and $F_j = H_j + h_j v_j$,
and let $\sigma\equiv \sigma_{F_j}$ denote the central projection
from $0$ to $F_j$.

Recall that we now have
\begin{equation} \Sigma \cap \partial B_{\rho_j} \cap
C(\frac 34 \varphi_0, v_j) = \emptyset \label{rim}
\end{equation}
but we assume that there are \emph{two different\/} points
$x_1,x_2\in \Sigma \cap \partial B_{\rho_j} \cap C(\varphi_0,
v_j)$ such that
\begin{equation}
\label{wide2} \ang \bigl(\pi_{H_j}(\sigma(x_1)),
\pi_{H_j}(\sigma(x_2))\bigr) \ge \frac \pi 3\, .
\end{equation}
Let $y_k=\sigma(x_k)$, $k=1,2$. Since the plane $P=\langle
0,x_1,x_2\rangle = \langle 0, y_1, y_2\rangle$, we can apply
Lemma~\ref{slanted1}  with $\varphi_0=\pi/4$, $ \varphi_1=\frac
\pi 3$  to select a third point $x_3\in \Sigma$ on a vertical
segment $I_{h_j,v_j}(z_3)$ (using \eqref{segments} for $i=j$),
where $z_3\in \gamma_j:=\partial B_{r_j}\cap H_j$, the outer
boundary of $A_j$ in $H_j$. This gives
\[
\eta_1 \rho_j \le \dist(x_3,P)\quad\AND\quad\eta_1\rho_j
\le |x_k-x_i|\le 2\rho_j \qquad
\mbox{for $k\not=i$, \quad $k,i=0,1,2,3$.}
\]
where now we have $\eta_1 = c_0(\pi/4,\pi/3) = \frac 12 (1-\cos
\frac \pi 6) = \frac 14$.
\heikodetail{

\bigskip

We have
\begin{eqnarray*}
|x_1-x_0| & = & \rho_j>\eta_1\rho_j,\\
|x_3-x_i| & \ge & \dist(P,x_3)\ge\eta_1\rho_j,\\
|x_2-x_0| & = & \rho_j>\eta_1\rho_j,\\
|x_2-x_1| & \ge & 2x\sin\pi/6=2h_j\tan(\frac 34 \varphi_0)\sin\pi/6=
2\rho_j\cos\varphi_0\tan(\frac 34 \varphi_0)\sin\pi/6\\
& \ge & 2\rho_j\frac 1 {\sqrt{2}} \frac 34 \varphi_0 \frac 2\pi \pi/6=
\rho_j\sqrt{2}\pi/16=\rho_j
\sqrt{2\pi^2}\pi/16>\rho_j\sqrt{18}\pi/16>\rho_j\sqrt{16}\pi/16=\rho_j
\pi/4
\end{eqnarray*}
(for $x$ see my drawing on my old p. 16)
\bigskip

}

\smallskip

It remains to verify that the angle $\ang (x_1-x_0,x_2-x_0) = \ang
(x_1,x_2)$ is in $[\eta_2, \pi - \eta_2]$ for some absolute
constant $\eta_2>0$ (possibly smaller than $\eta_1$), to verify
Condition (iii) in Definition \ref{Def-Ved}. This is intuitively
obvious but we give the details (without aiming at the best
possible bounds).

\smallskip

Let us suppose first that the two scalar products $ x_k\cdot
v_j$ ($k=1,2$) have the same sign. Write
\[
x_k=u_k+w_k, \qquad u_k:=\pi_{H_j}(x_k) \quad\mbox{for $k=1,2$,}
\]
and let $a_k:= {|w_k|}/{\rho_j} \equiv {|w_k|}/{|x_k|}$ for
$k=1,2$. Note that since $|x_1|=|x_2|=\rho_j$ and \eqref{rim} is
satisfied, we have in fact
\begin{equation}
a_k \le \sin\left(\frac\pi 2 - \frac 34 \varphi_0\right) =
\frac{5\pi}{16}, \qquad k=1,2. \label{ak-bd}
\end{equation}
\heikodetail{

\bigskip

see little drawing on my old p. 17:
$$
a_k=\frac{|w_k|}{\rho_j}=\sin\ang(x_k,u_k)\overset{\eqref{rim}}{\le}
\sin \left(\frac \pi 2-\frac 34 \varphi_0\right)=\sin\frac 5{16}\pi.
$$

\bigskip

}
Moreover, we have
\begin{equation}
\ang(u_1,u_2) =
\ang\bigl(\pi_{H_j}(\sigma(x_1)),\pi_{H_j}(\sigma(x_2))\bigr) \ge
\frac \pi 3 , \label{u1u2}
\end{equation}
(the first equality in \eqref{u1u2} holds since the scalar
products of $x_k,\ k=1,2,$ with $v_j$ are of the same sign). Set
$\psi:= \ang (x_1,x_2)$. Then, since the scalar products $
x_k\cdot v_j$ ($k=1,2$) have the same sign, we have $
w_1\cdot w_2=|w_1|\cdot |w_2| >0$, and therefore
\begin{eqnarray*}
0\, \le \cos\psi =\frac{x_1\cdot x_2}{|x_1|\cdot |x_2|} & = &
\frac{(u_1\cdot u_2) + (
w_1\cdot  w_2) }{\rho_j^2} =\frac{|u_1|\cdot |u_2|\cos\ang(u_1,u_2)}{\rho_j^2}+a_1a_2\\
& {\le} & \frac 12 (1-a_1^2)^{1/2}
(1-a_2^2)^{1/2} + a_1a_2 \qquad \mbox{by \eqref{u1u2}}\\
& \stackrel{(\ast)}{\le} & (1-\lambda)\Bigl( (1-a_1^2)^{1/2}
(1-a_2^2)^{1/2} + a_1a_2\Bigr) \\
& \le & 1-\lambda \qquad\mbox{by Young's inequality,}
\end{eqnarray*}
\heikodetail{

\bigskip

where we have also used $u_k^2+w_k^2=\rho_j^2$, or
$|u_k|=\rho_j\sqrt{1-a_k^2}$ by Pythagoras...

and
$$
(1-a_1^2)^{1/2}
(1-a_2^2)^{1/2} + a_1a_2\le \frac{1-a_1^2}{2}+
\frac{1-a_2^2}{2}+a_1a_2=1-\frac{(a_1-a_2)^2}{2}.
$$

\bigskip

}
provided that we can choose $\lambda \in (0,\frac 12)$ so that
$(\ast)$ holds, i.e., equivalently,
\begin{equation}
\label{aim-C2} \lambda a_1a_2 \le \left(\frac 12 - \lambda\right)
(1-a_1^2)^{1/2} (1-a_2^2)^{1/2}\, .
\end{equation}
Now, \eqref{ak-bd} implies that the left-hand of \eqref{aim-C2}
does not exceed $\lambda \sin^2 \frac{5\pi}{16}$ whereas the right-hand
side is certainly greater than $(\frac 12 - \lambda) \cos^2
\frac{5\pi}{16}$. Thus, \eqref{aim-C2} holds for every $\lambda
\le \frac 12 \cos^2 \frac{5\pi}{16}$, e.g. for $\lambda = \frac 12
\cos^2 \frac{\pi}{3}=\frac 18$ and then with strict inequality.
This gives $\cos \psi\in [0, \frac
78)$, i.e.,
\[
\eta_2\le \psi=\ang(x_1,x_2) \le \frac \pi 2\, .
\]
for $\eta_2 := \arccos \frac 78 \simeq 0.505 > 1/4$.

\smallskip
If the two scalar products $x_k\cdot v_j$ ($k=1,2$)
have different signs, we consider $\tilde x_2 = - x_2$. Since the
central projections $\sigma(x_2)$ and $\sigma(\tilde x_2)$
coincide, we can apply the previous reasoning to $x_1$ and $\tilde
x_2$, to obtain $\ang(x_1,\tilde x_2)\in [\eta_2,\frac \pi 2]$,
i.e. $\ang(x_1,x_2)\in [\frac \pi 2, \pi - \eta_2]$.

\medskip

With the choice $\eta:=\min\{\eta_1,\eta_2\}=\eta_1=1/4$ we have
verified that the tetrahedron $T=(x_0,x_1,x_2,x_3)$ satisfies
all conditions of Definition \ref{Def-Ved}, hence is of class
$\V(\eta,d_s(x_0)$ for $\eta=1/4$, which proves Part (i) of
Theorem \ref{goodtetra} also in Case 2.
This concludes the proof in Case 2.

\subsubsection{Case 3 (Antipodal position): the details}

\label{C3}

We deal with Subcases 3 (a) and 3(b), where we have stopped
the iteration, have set $N:=j$, with stopping distance $d_s(x_0):=
\rho_j=\rho_N.$
Recall that $H_j =(v_j)^\perp$, $F_j = H_j + h_j
v_j$, and $\sigma\equiv \sigma_{F_j}$ is the central projection
from $0$ to $F_j$.

As in Case 2, $\Sigma \cap \partial B_{\rho_j} \cap C(\varphi_0,
v_j)$ is nonempty but we have
\[
\Sigma \cap \partial B_{\rho_j} \cap C(\frac 34 \varphi_0, v_j) =
\emptyset.
\]
However, in this Case condition \eqref{wide} is violated, i.e.
\emph{for every two\/} points $x_1,x_2\in \Sigma \cap \partial
B_{\rho_j} \cap C(\varphi_0, v_j)$ we have
\begin{equation}
\label{antipo2} \ang \bigl(\pi_{H_j}(\sigma(x_1)),
\pi_{H_j}(\sigma(x_2))\bigr) < \frac \pi 3\, .
\end{equation}

We have already fixed $x_1\in\Sigma\cap\partial B_{\rho_j}\cap
C(\varphi_0,v_j)$  and assume now without loss of generality that
$v_j=(0,0,1)$, \,$x_1\cdot v_j>0$, and $u:=u_j=
\pi_{H_h}(x_1)/|\pi_{H_j}(x_1)| =(1,0,0).$  Hence the unit vector
$w:=w_j=(0,-1,0)$ determines the axis of the rotations $R^j_s$
defined in \eqref{rotations} which in turn were used to rotate
conical caps to obtain the sets $G_t^j$ and the stopping
rotational angle $t_0$ (see \eqref{Gsets} and \eqref{stop-Rs}). On
this basis the three subcases in Case 3 were distinguished. Let us
describe in some detail how we choose $x_2$ and $x_3$ in Subcase 3
(a) and (b).

\subsubsection*{Stopping the iteration in Subcase~3~(a)}

Let us first note that $t_0>0$. To see this, set
\[
X^j  :=  \{ y\in \rd\, \colon \,  (y\cdot v_j)(
y\cdot u ) \le 0\}, \qquad Y^j  :=  X^j \cap (C_{\rho_j}(\varphi_0,v_j)
\setminus \INT B_{\rho_j/2}),
\]
and note that if $R_s(C_{\rho_j}(\varphi_0,v_j)
\setminus \INT B_{\rho_j/2})$ contains a new point $y$ of $\Sigma$,
i.e. a point $y\in \Sigma \setminus S_j$, then we have in
fact $y\in R_s(Y^j)$. However, this cannot happen for $s$
arbitrarily close to $0$, as in Case~3 we have
\[
\dist \bigl(Y^j, \Sigma \cap X^j\bigr) >0
\]
due to \eqref{antipo2}, \eqref{interior_cond1}
and \eqref{umbrellaposition} for $i=j$ in (D),
and
\eqref{rho-cond} for $i=j$.
\heikodetail{

\bigskip

We have
$
\Sigma\cap\partial B_{\rho_j}\cap Y_j=\emptyset $ by virtue of \eqref{antipo2},and \eqref{umbrellaposition} for $i=j$ prohibits any surface point
on the outer
half of the the cone walls $\partial C_{t}(\varphi_0,v_j)\setminus\INT
B_{\rho_{j}/2}$ (notice that by \eqref{rho-cond} $\rho_j/2>\rho_{j-1}$,
and even no surface point in the interior of $S_j$ by \eqref{interior_cond1}
for $i=j$. That gives us positive distance of the closed sets $Y^j$
and $\Sigma\cap X_j$.

\bigskip

}
\smallskip

We choose $x_2 \in G_{t_0} \cap (\Sigma\setminus K^j_{\rho_j})$.
It is easy to see that if $x_2\cdot  v_j$ and $
x_1\cdot  v_j$ have the same sign, then
\begin{equation}
\begin{split}
  \frac{3}{16}\pi=\frac 34 \varphi_0  \le
\ang(x_1,x_2) & \le\ang(x_1,v_j)+\ang(v_j,R_{t_0}(v_j))+
\ang(R_{t_0}(v_j),x_2)\\
& \le  \varphi_0 + t_0\varphi_0 + \varphi_0 \le \frac 52 \varphi_0
= \frac{5}{8}\pi\,.
\end{split}
\label{angle-3a}
\end{equation}
If the scalar products $x_2\cdot  v_j$ and $x_1\cdot
v_j$ have different signs, then \eqref{angle-3a} holds with
$\tilde{x}_2=(-x_2)$ instead of $x_2$, so that in either case we
have
\begin{equation}
\label{ii-3a}  \frac{3}{16}\pi\le \ang(x_1,x_2) \le \pi-
\frac{3}{16}\pi = \frac {13}{16}\pi,
\end{equation}
and Condition (iii) of Definition~\ref{Def-Ved} holds with
$\theta:=3\pi/16$.

Now, take $P=\langle 0, x_1, x_2\rangle = \langle 0,
\sigma_{F_j}(x_1), \sigma_{F_j}(x_2)\rangle$ and apply
Lemma~\ref{slanted2} in connection with \eqref{segments} for $i=j$
in (F) to find the last good vertex $x_3$ on one of
the segments $I_{h_j,v_j}(z)$, where $z$ runs along the circle
$\gamma_j$ bounding the disk $H_j \cap B_{r_j}$,
$r_j=\rho_j\sin\varphi_0$. Then $\dist(x_3,P)\ge c_1(\varphi_0)\rho_j$
where $c_1(\varphi_0)=\frac 1{16}\sin2\varphi_0=\frac 1{16},$ which
verifies Condition (iv) of Definition \ref{Def-Ved}
with $\theta:=1/16$. Conditions (i) and (ii) of that definition
are easily checked, so that $T=(x_0,x_1,x_2,x_3)\in\V(\eta,d_s(x_0)$
(and therefore Part (i) of Theorem \ref{goodtetra}
is shown) for $\eta=1/16$ in Subcase 3 (a).
\heikodetail{

\bigskip

(i) is immediate since all points (as in the previous cases
are contained in $B_{\rho_j}$ and $d_s(x_0)=\rho_j.$
For (ii)
\begin{eqnarray*}
|x_1-x_0| & = & \rho_j > \eta\rho_j,\\
|x_2-x_0| & \ge & \rho_j/2 >\eta\rho_j\quad\textnormal{by \eqref{Gsets}},\\
|x_3-x_i| & \ge & \dist(x_3,P)\ge\eta\rho_j\quad\Fo i=0,1,2,\\
|x_2-x_1| & \ge & \dist(G_{t_0}\setminus S_j,x_1)\ge
 \dist(G_{t_0}\setminus S_j,I_{\rho_j,v_j}(0))
 =\rho_j\sin\varphi_0=\frac 1{\sqrt{2}}\rho_j>\eta\rho_j.
 \end{eqnarray*}

 \bigskip

 }

\subsubsection*{Stopping the iteration in Subcase~3~(b)}

Use \eqref{forward} to select a point $x_2 \in \Sigma \cap
\left(C_{2\rho_j} (\varphi_0, v_j^\ast) \setminus C_{\rho_j}
(\varphi_0, v_j^\ast)\right)$.

Assume first that $x_2\cdot  v^\ast_j >0$. Since, by the
definition of $R_s$ and $v^\ast_j=R_{1/2}(v_j)$, we have
\[
\ang(x_1,v^\ast_j) = \ang(x_1,v_j) + \ang(v_j,v^\ast_j) \in
[{\textstyle\frac{5}{4}}\varphi_0,{\textstyle\frac{3}{2}}\varphi_0],
\]
and $\ang (x_2, v_j^\ast) \le \varphi_0$,
\heikodetail{

\bigskip

$$
\ang(x_1,v_j)\ge \frac 34 \varphi_0\quad\AND\quad \ang(v_j,v^\ast_j)
=\frac 12\varphi_0,
$$
on the other hand, $\ang(x_1,v_j)\le \varphi_0.$ So that above we
have  $\frac 54 = \frac 34 + \frac 12$. (Note that $\frac 74 >
\frac 32$.)

\bigskip

}
two applications of the
triangle inequality for the spherical metric give
\[
\ang(x_1,x_2) \in
[{\textstyle\frac{1}{4}}\varphi_0,{\textstyle\frac{5}{2}}\varphi_0]
= [\pi/16, 5\pi /8]
\]
in that case.
\heikodetail{

\bigskip

$$
\ang(x_1,x_2)\ge \ang(x_1,v_j^\ast)-\ang(v_j^\ast,x_2)\ge\frac 5
4\varphi_0 - \varphi_0,
$$
and
$$
\ang(x_1,x_2)\le \ang(x_1,v_j^\ast)+\ang(v_j^\ast,x_2)\le\frac 52\varphi_0.
$$

\bigskip

}
If $x_2\cdot  v^\ast_j <0$, then we
estimate the angle $\ang(x_1,-x_2)$ in the same way. This yields
\[
\ang(x_1,x_2) \in [\pi/16, 15\pi/16],
\]
no matter what is the sign of $x_2\cdot  v_j^\ast$, which yields
Condition (iii) of Definition \ref{Def-Ved} for $\theta=\pi/16.$
Note that this estimate for the angle implies an estimate for the
distance, $\rho_j \sin(\pi/16) \le |x_2-x_1|$ being part of
Condition (ii) in Definition \ref{Def-Ved} for
$\theta=\sin\pi/16$.

\smallskip

To select $x_3$, we argue similarly to the proof of
Lemma~\ref{slanted2}.

Consider the affine plane $F\equiv F_j = H_j + h_j v_j$,
$h_j=\rho_j \cos \varphi_0$. Let $\sigma\equiv\sigma_{F}$ be the
central projection from $0$ to $F$. Set
\[
E := \sigma \left(C_{2\rho_j} (\varphi_0, v_j^\ast) \right) \,
\subset \, F\, ;
\]
this is a filled ellipse in $F$. We have $y_2=\sigma(x_2)\in E$.
Consider now the point $y_1=\sigma(x_1)\in F$. The plane
$P=\langle 0, x_1, x_2\rangle $ is equal to $\langle 0, y_1,
y_2\rangle$. The straight line $l=P\cap F$ passes through
$y_1,y_2$, and has to intersect $\partial E$ and $l_2$, where the
straight line
\[
l_2 :=P_2 \cap F \qquad\mbox{for}\quad P_2:= (R(7\pi/8,w)(v_j))^\perp=
(R_{7/2}(v_j))^\perp,
\]
is tangent to $\partial E$ in $F$, and the direction of $l_2$ is
perpendicular to $v_j$ and to $u=(1,0,0)$. Let $y_3$ be that point in $\partial E
\cap l$ --- which in general contains two points --- which is
closer to $y_1$, and let $\{y_4\}:=l_2\cap l$. Then it is easy to
see that $y_4$ lies on $l$ between $y_3$ and $y_1$. Therefore, $l$
contains a point $y_5$ such that (see the figure below)
\begin{equation}
\label{y5} \ang\bigl(\pi_{H_j}(y_5), \pi_{H_j}(y_1)\bigr) = \phi
:= \arccos \left[\cot\varphi_0\left(\tan \frac \pi 8\right)\right]
=\arccos \left(\tan \frac \pi 8\right)= 1.1437\ldots \, ,
\end{equation}
and we have $P=\langle 0, x_1, x_2 \rangle = \langle
0, y_1,y_5\rangle$.
\heikodetail{

\bigskip

look at my drawing from the side (old p. 20) to find
that the distance  $\xi$ of $l_2$ from the center in in Figure 2
is given by
$\xi=h_j\tan\pi/8$. Then consult again Figure 2 to estimate
$$
\cos\phi=\frac{\xi}{r_j}=\frac{h_j}{r_j}\tan\frac{\pi}{8}=
\frac{\rho_j\cos\varphi_0}{\rho_j\sin\varphi_0}\tan\pi/8.
$$

\bigskip

}%

\begin{center}
\includegraphics*[totalheight=8cm]{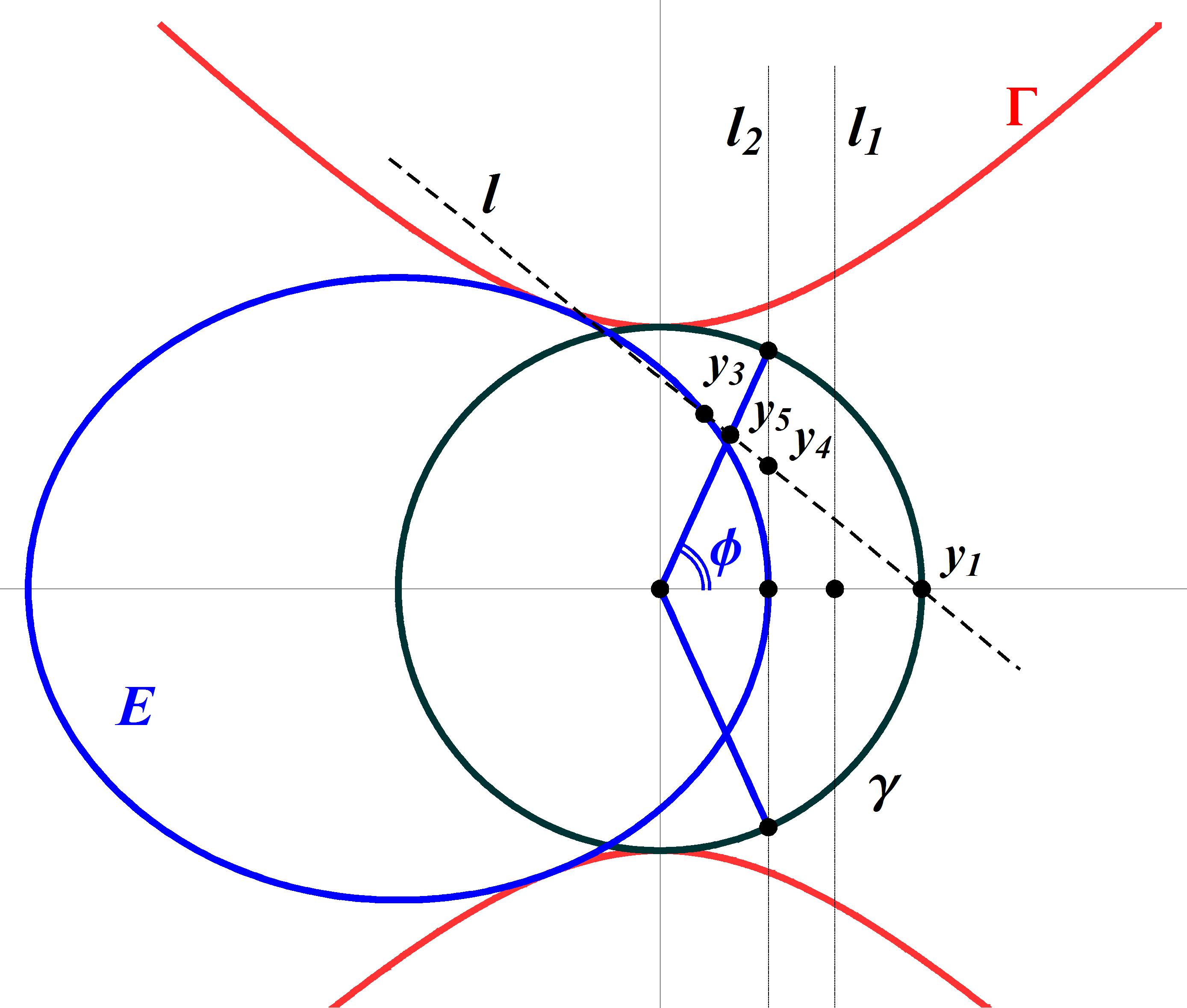}
\end{center}

{\footnotesize

\medskip\noindent\textbf{Fig. 5.} The configuration in $F$ discussed
above. The (slanted, dashed) line $l$ passes through
$y_1=\sigma(x_1)$ and some other point (not shown) belonging to
the ellipse $E$. The four points depicted on $l$ are, from right
to left, $y_1$, $y_4$, $y_5$ and $y_3$. Note that $y_4$ is always
situated between $y_3$ (which is on the boundary of the ellipse)
and $y_1$. The position of $y_5$, which is chosen on $l$ so that
the angle $\ang(\pi_{H_j}(y_1),\pi_{H_j} (y_5))=\phi$, may change,
depending on the slope of $l$ and position of $y_1=\sigma(x_1)$ (a
special case $\sigma(x_1)=x_1\in F$ is shown here). For some
positions of $x_2$ considered in Subcase~3~(b), when $l$ is not so
close to a tangent to $E$, we might obtain the order: $y_1$, then
$y_4\in l_2$, then $y_3\in
\partial E$, and finally $y_5$ satisfying \eqref{y5}.

}

\smallskip

Applying Lemma~\ref{slanted1} with $\varphi_1:=\phi$, we find a
point $z_3\in  H_j\cap\partial  B_{r_j}$, $r_j=\rho_j
\sin\varphi_0$, and because of \eqref{segments} for $i=j$ in (F)
the last vertex $x_3\in I_{h_j,v_j}(z_3)\cap\Sigma\subset
B_{\rho_j}$ of a good tetrahedron. The estimate from
Lemma~\ref{slanted1} gives now
$$\dist(x_3,P) \ge
c_0(\varphi_0,\phi) \rho_j = 0.0795\ldots \cdot \rho_j\, . $$
Since $c_0(\varphi_0,\phi) < \cos(\pi/16)$, it is easy to see that
all the distances $d_{ik}:=|x_i-x_k|$, $i\not=k$, satisfy
\[
0.0795\ldots \cdot \rho_j\,\,\le\,\, d_{ik}\,\, \le 3 \rho_j.
\]
\heikodetail{

\bigskip

Let's check:
\begin{eqnarray*}
|x_2-x_1| & \ge & \rho_j\sin(\pi/16)\quad\textnormal{estimated before,}\\
|x_3-x_i| & \ge & \dist(x_3,P)\ge c_0\rho_j\quad\Fo i=0,1,2,\\\
|x_1-x_0| & = & \rho_j > c_0\rho_j,\\
|x_2-x_0| & \ge & \rho_j > c_0\rho_j,
\end{eqnarray*}
and all points (including $x_2$ lie in $B_{2d_s(x_0)}=B_{2\rho_j}$
establishing (i), (ii) of Definition \ref{Def-Ved}.

\bigskip

}%

\smallskip

All the conditions of Definition \ref{Def-Ved} are verified now,
and we conclude that $T=(x_0,x_1,x_2,x_3)\in\V(\eta,d_s(x_0))$
for $\eta:=c_0(\varphi_0,\phi)=0.0795$ and $d_s(x_0)=\rho_j=\rho_N,$
which implies the validity
of Part (i) of Theorem \ref{goodtetra} for this last Case where the
iteration was stopped.  Part (ii) follows from
 Corollary \ref{cor:4.4}
 below.
 \hfill $ \Box $

\subsection{Estimates for perturbed tetrahedra}
\label{Sec:3.4}

\begin{lemma}
\label{TTprim}
Assume that $x_0=0,x_1,x_2,x_3\in\rd$ satisfy

\begin{enumerate}
\renewcommand{\labelenumi}{{\rm (\roman{enumi})}}

\item  $\eta d \le |x_i-x_j| \le \eta^{-1} d$
for all $i\not =j$, $i,j=0,1,2,3$,
\item $\dist(x_3,\langle x_0,x_1,x_2\rangle) \ge \eta d$;

\item $\eta\le \ang(x_1-x_0,x_2-x_0) \le \pi-\eta$,
\end{enumerate}
where $\eta\in (0,\frac 12)$ and $d>0$. Then, there exists a
number $\eps=\eps(\eta)\in (0,1/4)$ such that
\begin{equation}
\label{pert} \dist \bigl(y_3, \langle y_0,y_1,y_2\rangle\bigr)\ge
\frac 12 \eta d
\end{equation}
whenever $y_i \in B_{\eps d}(x_i)$ for $i=0,1,2,3$.
\end{lemma}

\medskip\noindent\textbf{Proof.} W.l.o.g. we may assume that
$x_0=0.$ Let $y_i = x_i+ v_i$ with
$|v_i|\le \eps d$ for $i=0,1,2,3$; we shall fix $\eps \in (0,1/4)$ later on.
Since the left-hand side of \eqref{pert} is invariant under
translations, it is enough to prove \eqref{pert} for the quadruple
$(y_0,y_1,y_2,y_3)$ shifted by $- v_0$. Thus, from now on we
suppose that
\[
y_0=x_0=0, \qquad y_j = x_j + w_j, \quad\mbox{where } |w_j| \le 2
\eps d \quad\mbox{for $j=1,2,3$.}
\]
By (iii), (i), and the fact that $\eta<1/2$,  we have
\[
d^2\eta^4 \le d^2 \eta^2 \sin \eta \le |x_1\times x_2| \le |x_1|\,
|x_2| \le d^2 \eta^{-2}\, .
\]
\heikodetail{

\bigskip

Use on the left-hand side that
$$
\sin\eta > (2/\pi)\eta > (1/2)\eta >\eta^2
$$

\bigskip

}
Moreover, $ y_1 \times y_2 = (x_1 \times x_2) + v$, where the
remainder vector $v$ satisfies by Assumption (i)
\[
|v| = |w_1\times x_2 + x_1 \times w_2 + w_1\times w_2| \overset{\textnormal{\rm (i)}}{\le} 2 \cdot
2\eps d \cdot d\eta^{-1} + (2\eps d)^2\le d^2\eta^{-1}(4\eps+4\eps^2) <
5\eps d^2\eta^{-1}
\]
(the last inequality is satisfied for all $\eps \in (0,1/4)$ and
$0<\eta \le 1$). Thus,
\[
|y_1\times y_2| \le \frac 32 |x_1\times x_2|
\]
if $|v| \le \frac 12 d^2\eta^4 \le \frac 12 |x_1\times x_2|$, and
the last condition is satisfied whenever
\begin{equation}\label{first_eps_condition}
10\eps \le \eta^5.
\end{equation}
Since $y_0=0=x_0$, for all such choices of $\eps$ we have
according to Assumption (ii)
\begin{eqnarray*}
\dist \bigl(y_3, \langle y_0,y_1,y_2\rangle\bigr) & = &
\frac{|\langle y_3, y_1\times y_2\rangle|}{|y_1\times y_2|} \\
& \ge & \frac{2|\langle y_3, y_1\times y_2\rangle|}{3|x_1\times
x_2|} \,\, \ge\,\, \frac{2|\langle x_3, x_1\times
x_2\rangle|}{3|x_1\times x_2|} - R\\
& \overset{\textnormal{\rm (ii)}}{\ge} & \frac 23 d\eta - R ,
\end{eqnarray*}
where, by the triangle inequality,
\begin{eqnarray*}
0\, \le \, R & \le & \frac{2}{3|x_1\times x_2|} (|w_3|\, |x_1|\,
|x_2| +
|w_3|\, |v| + |x_3|\, |v|) \\
& \le & \frac 23 (d^2\eta^4)^{-1}( 2\eps d \cdot d^2 \eta^{-2} +
\eps
d \cdot d^2\eta^4 + d\eta^{-1} \cdot 5\eps d^2 \eta^{-1}) \\
& < & 6\eps d\eta^{-6}
\end{eqnarray*}
as $0<\eta < 1$ hence $\eta^4<\eta^{-2}$ for the last inequality.
Choosing $\eps=\eps(\eta)\in (0,1/4)$ so small
that $R \le 6\eps
d\eta^{-6} \le \frac 16 d\eta$ in addition to the
requirement in \eqref{first_eps_condition}, we conclude the proof. \hfill
$\Box$

\begin{corollary}\label{cor:4.4}
Given $d>0$
one finds for any $\eta\in (0,1/2)$
a constant $\alpha=\alpha(\eta)\in (0,\eta/20)$
such that for all tetrahedra $T\in\V(\eta,d)$ one has
$$
T'\in\V(\frac \eta 2,\frac 32 d)\quad\Foa \|T-T'\|\le\alpha d.
$$
\end{corollary}
We omit the proof since it relies on simple distance estimates
using the triangle inequality and on Lemma \ref{TTprim}.

\subsection{Large projections and forbidden conical sectors}

It is clear that conditions (A)--(F) stated at the beginning of
Section~\ref{meat} combined with the lower bound for stopping
distances obtained in Proposition~\ref{prop:3.5} imply the
statement of Proposition~\ref{sectors} for all points $x\in
\Sigma^\ast$.

Using density of $\Sigma^\ast$ and closedness of $\Sigma$ it is
easy to see that Proposition~\ref{sectors} does hold also for all
$x\in \Sigma\setminus \Sigma^\ast$.

Indeed, fix $x\in \Sigma$ and $r< R_0=R_0(E,p)$. Choose a
sequence of $x_i\to x$, $x_i\in \Sigma^\ast$. For each $x_i$, let
$H_i$ and $v_i$ be the plane and unit vector whose existence is
given by Proposition~\ref{sectors} for points of $\Sigma^\ast$.
Set $D_i:=H_i\cap B(x_i, r/\sqrt{2})$.

Passing to subsequences if necessary, we can assume that $H_i$ and
$v_i$ converge as $i\to \infty$  to a plane $H$ and a unit
vector $v$. We shall show that $H$ and $v$ satisfy the
requirements of Proposition~\ref{sectors} for~$x$ and~$r$.

For each $w\in D:=H\cap B(x,r/\sqrt{2})$ we select $w_i\in D_i$
with $|w_i-x_i|=|w-x|$ such that $w_i\to w$ as $i\to\infty$. By
\eqref{bigpi-1} applied for $x_i$, $\Sigma$ contains points
$y_i=w_i+t_i v_i$ where the coefficients $t_i$ satisfy
\[
|t_i|^2 \le r^2 - |w_i-x_i|^2 = r^2 - |w-x|^2.
\]
Again, without
loss of generality we can assume that $t_i\to t$ as $i\to\infty$,
so that
\[
y_i=w_i + t_i v_i \to y= w+tv, \qquad |t|^2 \le r^2 - |w-x|^2.
\]
It is clear that $y\in \Sigma \cap B(x,r)$ and $\pi_H(y)=w$ so
that \eqref{bigpi-1} holds at $x$.

Finally, if one of \eqref{C+in}--\eqref{C-out} were violated with
our choice of $H$ and $v$, then the respective condition would be
violated for $x_i$, $r$, $H_i$ and $v_i$ for all $i$ sufficiently
large, a contradiction.


\section{Uniform flatness and oscillation of the tangent planes}

\label{sec:5} \setnumbers

Throughout this section we assume that $\Sigma = \partial U$ is a
closed, compact admissible surface in $\rd$, with
\[
\M_p(\Sigma) < E <  \infty
\]
for some $p>8$. As was shown before in Theorem \ref{Thm:3.1}, all
such $\Sigma$ are Ahlfors regular with bounds depending only on
the energy, i.e. there exists an $R_0=R_0(E,p)>0$ whose precise
value was given in \eqref{R0Ep} such that
\begin{equation}\label{5.1}
\H^2(\Sigma \cap B(x,R))\ge \frac \pi 2 R^2 \qquad \mbox{for all
$x\in \Sigma$ and $R\in (0, R_0]$.}
\end{equation}
We shall show that each such $\Sigma$ is in fact a manifold of
class $C^1$. To this end, we shall show that the tangent plane to
$\Sigma$ exists and satisfies an a priori H\"{o}l\-der estimate.
This a priori estimate allows to cover $\Sigma$ by a finite number
of balls, with radii depending only on $p$ and the bound for
energy, such that in each of these balls $\Sigma$ is a graph of a
$C^1$ function with H\"{o}lder continuous derivatives, see
Corollary~\ref{cor:5.4new}. This fact will be used also later in
Section \ref{sec:7} when dealing with sequences of admissible
surfaces with equi\-bounded energy.

Our aim in this section  will be to estimate the so-called {\it
beta numbers}; see e.g. the introductory chapter of
\cite{davidsemmes},
\begin{equation}
\label{Jones-b} \beta_\Sigma(x,r) := \inf\left\{ \sup_{y\in
\Sigma\cap B(x,r)} \frac{\dist (y,F)}{r} \quad \colon \quad
\mbox{$F$ is an affine plane through $x$} \right\}
\end{equation}
for small radii $r$ and points $x\in \Sigma\, $, and to show that
\begin{equation}
\beta_\Sigma(x,r) \le C(E,p) r^\kappa \label{beta-small}
\end{equation}
where $\kappa= \kappa(p) = (p-8)/(p+16)>0$. One of the issues is
that we want to have such estimates for all $r < R_1(E,p)$ where
$R_1(E,p)$ is a constant that does not depend on $\Sigma$.

It is known that for the class of Reifenberg flat sets with
vanishing constant uniform estimates like \eqref{beta-small} imply
$C^{1,\kappa}$ regularity, cf. for example David, Kenig and Toro
\cite[Section 9]{davidkenigtoro}, or Preiss, Tolsa and Toro
\cite[Def. 1.2 and Prop. 2.4]{ptt}. In our case, we a priori know
that $\Sigma\in \A$ and this information by itself does not imply
Reifenberg flat\-ness. However, we establish \eqref{beta-small}
inductively; while doing that, we can simultaneously ensure that
$\Sigma$ is Reifenberg flat with a vanishing constant in a scale
depending only on the energy.

In order to show precisely what is the role of energy bounds, we
give all details of that reasoning. Everything is based on
iterative applications of Proposition~\ref{sectors} and of the
following simple lemma.

\begin{lemma}[Flat boxes]
\label{flat} Suppose that $\M_p(\Sigma) < E$ for some $p> 8$.
Then, for any given number $1> \eta>0$ there exist two positive
constants $\eps_0 = \eps_0(\eta) > 0$ and $c_1 = c_1(\eta,p) >0$
such that whenever a triple of points $\Delta = (x_0, x_1, x_2)
\in \Sigma^3$ satisfies
\[
\Delta \in \Spac (\eta,d), \qquad d\le  R_0(E,p)
\]
where $R_0(E,p)$ is given by \eqref{R0Ep}, then we have
\begin{equation}
\label{coin} \Sigma \cap B(x_0,3d)\,\, \subset \,\, U_{\eps
d}(\langle x_0,x_1,x_2 \rangle )
\end{equation}
for each $\eps \in (0,\eps_0(\eta))$ which satisfies the balance
condition
\begin{equation}
\label{eps-d} \eps^{16+p} d^{8-p} \ge c_1(\eta,p) E.
\end{equation}
\end{lemma}
In other words, we have
\[
\beta_\Sigma(x_0,3d) \le \frac{\eps}3
\]
(and also a slightly weaker inequality $\beta_\Sigma(x_0,d) \le
\eps$) whenever we can find an appropriate triple of points of
$\Sigma$ and \eqref{eps-d} is satisfied. Note that the balance
condition \eqref{eps-d} is satisfied for $\eps\approx E^{1/(p+16)}
d^\kappa$, so that the `boxes' $B(x_0,3d)\cap U_{\eps d} (\langle
x_0,x_1,x_2\rangle)$ become indeed flatter and flatter as the
scale $d\to 0$. \heikodetail{

\bigskip

Look at
\begin{eqnarray*}
\beta_\Sigma(x_0,3d) &:= &\inf\left\{
\sup_{y\in\Sigma\cap B(x_0,3d)}\frac{\dist(y,F)}{3d}:F\,\,\textnormal{affine}\right\}\\
& \le & \sup_{y\in\Sigma\cap B(x_0,3d)}\frac{\dist(y,\langle
x_0,x_1,x_2\rangle)}{3d}\overset{\eqref{coin}}{\le}\frac{\eps
d}{3d}=\frac \eps 3,
\end{eqnarray*}
or a little weaker (since more local)
\begin{eqnarray*}
\beta_\Sigma(x_0,d) &:= &\inf\left\{
\sup_{y\in\Sigma\cap B(x_0,d)}\frac{\dist(y,F)}{d}:F\,\,\textnormal{affine}\right\}\\
& \le & \sup_{y\in\Sigma\cap B(x_0,d)}\frac{\dist(y,\langle
x_0,x_1,x_2\rangle)}{d}\overset{\eqref{coin}}{\le}\frac{\eps d}{d}
=\eps.
\end{eqnarray*}

\bigskip

}

\begin{REMARK}\label{badscaling}\rm This lemma and its iterative
applications in the proof of Theorem~\ref{thm:5.3} are one of the
main reasons behind our choice of definition of $\M_p$. The proof
presented below shows that for any inte\-grand $\K_s(T)$
satisfying
\[
\K_s(T)\approx \frac{h_{\text{min}}(T)}{(\diam T)^{2+s}}, \qquad
s>0,
\]
for which the scaling invariant exponent equals $8/(1+s)$, the
appropriate balance condition replacing \eqref{eps-d} would be
\[
\eps^{16+p} d^{8-(1+s)p} \gtrsim \text
{Energy}:=\int_{\Sigma^4}\K_s(T)^p\, d\mu\, .
\]
For $p>8/(1+s)$ this would yield, instead of \eqref{beta-small}
above, an inequality of the form $\beta_\Sigma(x,r) \lesssim
r^{\kappa(s,p)}$ with $\kappa(s,p) = (p+sp-8)/(p+16)$. However,
for $p>24/s$ we have $\kappa(s,p)> 1$, and reasoning as in the
proof of Theorem~\ref{thm:5.3} below one could show that the
normal to $\Sigma$ is H\"{o}lder continuous with exponent
$\kappa(s,p)>1$, i.e. constant! Because of that we do not work
with the $c_{\text{MT}}$ curvature introduced by Lerman and
Whitehouse in \cite{LW08b}: for sufficiently large $p$, the only
surface with finite energy would be a plane.
\end{REMARK}

\medskip
\noindent\textbf{Proof.} We argue by contradiction. Suppose that
some point $x_3 \in \Sigma \cap B(x_0,3d)$ does not belong to
$U_{\eps d}(P)$, $P:=\langle x_0,x_1,x_2\rangle$. Fix
$\eps_0=\eps_0(\eta)>0$ so small that if $\eps<\eps_0$, then for
all tetrahedra $T'$ with vertices $x_i'\in B(x_i,\eps^2 d)$,
$i=0,1,2,3$ one has
\begin{equation}\label{T'_est}
\dist (x_3',\langle x_0',x_1',x_2'\rangle )\ge \frac{\eps d}{4}
=\frac{\eps}6 \cdot \frac{3d}2 \qquad\mbox{and}\qquad
\Delta(T')=(x_0',x_1',x_2')\in \Spac (\eta/2, 3d/2).
\end{equation}
(An exercise, similar to the proof of Lemma~\ref{TTprim}, shows
that one can take e.g. $\eps_0(\eta)=\eta^2/200$.)
Now, since $\eps^2 d < d\le
R_0(E,p)$, we have by \eqref{5.1}
\[
\H^2 (\Sigma \cap B(x_i,\eps^2 d)) \ge \frac \pi 2 (\eps^2 d)^2
> \eps^4 d^2
\]
for $i=0,1,2,3$. Invoking Lemma~\ref{R-flat} with $\kappa=\eps/6$
as suggested by \eqref{T'_est},
 we obtain an estimate of the integrand,
\[
\K(T') \ge \frac{1}{50^2} \left(\frac \eta 2\right)^3
\frac{\eps}6\cdot \frac 2{3d} \,=\, \frac{1}{18\cdot
10^4}\frac{\eta^3\eps}{d}\, , \qquad T'=(x_0',x_1',x_2',x_3')\, .
\]
Integrating this inequality w.r.t. $T'\in \Sigma^4 \cap \B_{\eps^2
d}(T)$, we immediately obtain
\begin{eqnarray*}
E & > & \M_p(\Sigma) \ge \int_{\Sigma^4 \cap \B_{\eps^2 d}(T)}
\K^p(T') \, d\mu(T') > \left(\eps^4 d^2\right)^4 \,
\left(\frac{\eta^3\eps}{18\cdot 10^4d}\right)^p \\
& = &\eta^{3p}(18\cdot 10^4)^{-p} \eps^{16+p}d^{8-p},
\end{eqnarray*}
which is a contradiction to \eqref{eps-d} if we choose
$c_1(\eta,p)=\eta^{-3p} (18\cdot 10^4)^{p}$. \hfill $\Box$

\medskip\noindent\textbf{Remark.} From now on, we fix $\eta>0$ to
be the constant whose existence is asserted in
Theorem~\ref{goodtetra}, and we write
\begin{equation}\label{c1-eta-p}
c_1(p):=c_1(\eta,p)
\end{equation}
for that fixed value of $\eta$.

\begin{lemma}[Good triples of points of $\Sigma$]
\label{triples} Let $\Sigma\in \A$, $p>8$ and $\M_p(\Sigma) <
\infty$. Suppose that $x\in \Sigma$, $y\in \Sigma$ and $0<d=|x-y|<
d_s(x)$, where $d_s(x)$ is the stopping distance from
Theorem~\ref{goodtetra}. Then there exists a point $z\in \Sigma
\cap B(x,d)$ and an affine plane $H$ passing through $x$ such that

\begin{enumerate}
\renewcommand{\labelenumi}{{\rm (\roman{enumi})}}

\item $\Delta=(x,y,z) \in \Spac (\eta,d)$, where $\eta$ is the
constant from Theorem~\ref{goodtetra};

\item $\pi_H (\Sigma \cap B(x,d)) \supset H \cap B(x,d\sin\varphi_0
)$, where $\varphi_0=\frac \pi 4$;

\item $\ang(H, P) \le \alpha_0^\ast$, where $P=\langle x,y,z\rangle$
and
\begin{equation}
\alpha_0^\ast := \frac\pi2 - \arctan \frac 1{\sqrt{2}} =
0.955\ldots < \frac \pi 3\, . \label{alphazero}
\end{equation}
\end{enumerate}
\end{lemma}

\medskip\noindent\textbf{Proof.} W.l.o.g. we suppose that $x=0\in \rd$.
Applying Proposition~\ref{sectors}, we find $v\in \stwo$ and
$H=(v)^\perp$ such that \eqref{bigpi-1}, \eqref{C+in} and
\eqref{C-out} do hold for $r=d=|x-y|$, $H$ and $v$. In particular,
\begin{equation}
D:= H\cap B(x,d/\sqrt{2}) \, \, \subset \,\, \pi_{H}(B(x,d) \cap
\Sigma)\, ,\label{bigpi-2}
\end{equation}
and by \eqref{C+in}--\eqref{C-out}
\begin{equation}
\frac \pi 4 \le \ang(y-x,v) \le \frac{3\pi}{4}\, . \label{yx-vi}
\end{equation}
By \eqref{bigpi-2}, for each $w$ in the boundary circle of the
disk $D$ the segment $I(w):=I_{d/\sqrt{2},v}(w)$ (cf.
Section~\ref{basic} for the definition) contains at least one
point of $\Sigma$. Choose $w_0\in D$ such that $w_0-x \perp
\pi_{H}(y-x)$ and $|w_0-x|=d/\sqrt{2}$ and then choose any point
$z\in \Sigma \cap I(w_0)$. We claim that the conditions of the
lemma are satisfied by that point $z$ and $H$.

Indeed, we have $z\in B(x,d)$ and $\min (|z-x|, |z-y|)\ge
d/\sqrt{2}\ge \eta d$. \heikodetail{

\bigskip

just used $\eta < 1/2 $ ---  not more

\bigskip

}%
By choice of $z$ and $w_0$, we also have
\[
(z-x)\cdot (y-x) =  z\cdot y = (z-\pi_H(z))\cdot (y-\pi_H(y)) =
\pm |z-\pi_H(z)|\, |y-\pi_H(y)|\, .
\]
\heikodetail{

\bigskip

For the second equality :
$$
z\cdot \pi_H(y) =0 \quad\textnormal{by construction $w_0-x\perp
\pi_H(y-x)$}
$$
and therefore
$$
\pi_H(z)\cdot y=0 \AND \pi_H(z)\cdot\pi_H(y)=0.
$$

\bigskip

}%
Thus, $|\cos\ang(z,y)| = (|z-\pi_H(z)|/|z|)\, (|y-\pi_H(y)|
/|y|) \le (\cos \varphi_0)^2 = \frac 12$, so that $\ang(z,y) \in
[\frac \pi3, \frac{2\pi}3]$. This implies that $\Delta=(x,y,z)$ is
$(\eta,d)$-wide, i.e. $\Delta\in\mathscr{S}(\eta,d).$
\heikodetail{

\bigskip

use
$$
\cos\ang(z,y)=\frac{z\cdot y}{|z||y|}=\frac{(z-\pi_H(z))\cdot
(y-\pi_H(y)}{|z||y|}
$$
and since we do not know the exact height of $z$ on that vertical
segment $I(w_0)$ we can only say
$$
\frac{|z-\pi_H(z)|}{|z|}=\cos\alpha\le\cos\varphi_0
$$
for some $\alpha\in [\varphi_0,\pi/2).$ (Look at cone drawing from
the side, see back of old page 23)

\bigskip

}%

To check (iii), one solves an exercise in elementary geometry. For
that let $P:=\langle x,y,z\rangle$. It is enough to check that
$\frac \pi 2\ge \ang(P,v_i) \ge \arctan (1/\sqrt{2})$ and then use
$\ang(P,H)=\frac{\pi}2 - \ang(P,v)$. To compute $\ang(P,v)$, let
$F= H+ hv$, $h=d\cos\varphi_0=d/\sqrt{2}$ and note that the
distance $\delta:=\dist(l_1,l_2)$ between the two straight lines
$l_1:= P\cap F$ and $l_2:= \{x+ sv \colon s\in \R\}\perp F$
satisfies $\delta \ge h/\sqrt{2} = d/2$. This gives the desired
estimate of the angle. \hfill $\Box$ \heikodetail{

\bigskip

An image here might help (see back of old page 23). One uses $h\le
|\sigma_F(z)-\sigma_F(x)|$ where $\sigma_F$ is the central
projection from $x$ onto $F$. On the image one draws the worst
case position of $l_1$ closest possible to the vertical $l_2$, and
then uses Pythagoras to estimate this minimal distance between
$l_1$ and $l_2$. From there one estimates
$$
\tan\ang(P,v_i)=\frac{\delta}{h}\ge
\frac{h}{h\sqrt{2}}=\frac{1}{\sqrt{2}}.
$$

\bigskip

}%

\begin{theorem}[Existence and oscillation of the tangent plane]\label{thm:5.3}
Assume that $\Sigma\in \A$ and $\M_p(\Sigma) < E$ for some $p>8$.
Then, for each $x\in\Sigma$ there exists a unique plane
$T_x\Sigma$ (which we refer to as \emph{tangent plane of $\Sigma$
at $x$}) such that
\begin{equation}\label{distance_from_tangent-plane}
\dist(x',x+T_x\Sigma)\le C(p,E)|x'-x|^{1+\kappa} \quad\Foa
x'\in\Sigma\cap B_{\delta_1}(x),
\end{equation}
where $\kappa:={(p-8)/(p+16)}>0$ and $\delta_1=\delta_1(E,p) >0$.
Moreover, there is a constant $A=A(p)$ such that whenever $x,y\in
\Sigma$ with $0<d=|x-y|\le \delta_1(E,p)$, then
\begin{equation}
\label{osc-normal} \ang(T_x\Sigma,T_y\Sigma) \le A(p)\,
E^{1/(p+16)}\,  d^\kappa\, .
\end{equation}
\end{theorem}
\begin{REMARK}
In fact a possible choice for $\delta_1(E,p)$ is
\begin{equation}\label{delta_1}
\delta_1(E,p):=\min\left\{1, R_0(E,p),
\left(\frac{\mu_0\kappa}{400} \right)^{1/\kappa} \bigl(c_1(p)
E\bigr)^{-1/(p-8)} \right\},
\end{equation}
where $R_0(E,p) $ is the absolute constant given in \eqref{R0Ep}
of Theorem \ref{Thm:3.1}, $c_1(p)$ is defined in \eqref{c1-eta-p},
and $\mu_0:=\frac 14 \left(\frac \pi3 - \alpha_0^\ast\right)\ $.
\end{REMARK}

\medskip
\noindent\textbf{Proof of Theorem \ref{thm:5.3}.} Let us describe
first a rough idea of the proof.

To begin, we use Lemma~\ref{triples} and select $z\in \Sigma\cap
B(x,d)$ such that the triple $\Delta =  (x,y,z) \in \Spac
(\eta,d)$. Then, fixing $\delta_1(E,p)$ small and setting
\[
d_N := d/{10}^{N-1}, \qquad \eps_N \quad\mbox{such that }
\eps_N^{16+p}d_N^{8-p} \equiv c_1(p) E \quad\mbox{for
$N=1,2,\ldots$},
\]
we shall find triples of points, $\Delta_N = (x,y_N,z_N)\in
\Sigma^3$, such that $y_N,z_N\in B(x,2d_N)\,$ and the angle
$\gamma_N = \ang(y_N-x, z_N-x) \approx \frac \pi 2$ with a small
error bounded by $C\sum \eps_N$ where $C$ depends only on $p$ and
$E$. The crucial tool needed to select $y_N,z_N$ is the knowledge
that $\Sigma\cap B(x,d_1)$ has large projections onto some fixed
plane.

Thus, an application of Lemma~\ref{flat} shall give
\begin{equation}
\label{N-boxes} \Sigma \cap B(x,3d_N) \,\, \subset\,\, U_{\eps_N
d_N}(P_N), \qquad P_N=\langle x,y_N,z_N\rangle.
\end{equation}
Moreover, we shall check that the planes $P_N$ satisfy $\ang
(P_{N+1},P_N) \le C \eps_N$, and $P_1$ is close to $P_0=\langle
x,y,z\rangle$. Thus, the sequence $(v_N)$ of normal vectors to
$P_N$ is a Cauchy sequence in $\stwo$. This allows us to set the
(affine) tangent plane $P\equiv T_x\Sigma +x$ to be the limit
plane of the $P_N$, and to prove that $P$ does not depend on the
choice of $y_N,z_N$ and $P_N$ (which is by no means unique). (It
is intuitively clear that $P=\lim P_N$ should be equal to the
affine tangent plane to $\Sigma$ at all points where $\Sigma$ a
priori happens to have a well defined tangent plane.) The whole
reasoning gives
\[
\ang (T_x\Sigma,P_0)\le C \eps_1 = C' d^\kappa.
\]
Reversing the roles of $y$ and $x$, we run a similar iterative
reasoning to obtain the above inequality with $x$ replaced by $y$.
An application of the triangle inequality, combined with a routine
examination of the constants, ends the proof.

\smallskip

Let us now pass to the details.

\medskip

Again, we assume for the sake of convenience that $x=0$. Set
\begin{equation}
d_N: = \frac{d}{10^{N-1}}, \quad d=|x-y|, \qquad N=1,2,\ldots,
\label{d-N}
\end{equation}
and let $\eps_N$ be defined by
\begin{equation}
\eps_N^{16+p}d_N^{8-p} \equiv c_1(p) E, \qquad N=1,2,\ldots
\label{eps-N}
\end{equation}
Note that
$$
\eps_N=\left(\frac{c_1(p)E}{d^{8-p}}\right)^{\frac{1}{16+p}} \cdot
\Big(10^{N-1}\Big)^\frac{8-p}{16+p}\rightarrow 0\quad\As
N\to\infty.
$$
Moreover, by our choice of $\delta_1$ in \eqref{delta_1},
 \begin{eqnarray}
  200\sum_{N=1}^\infty \eps_N
 & = & 200 \bigl(c_1(p) E\bigr)^{1/(p+16)} \sum_{N=1}^\infty
 d_N^\kappa, \qquad \kappa :=\frac{p-8}{p+16}>0, \nonumber  \\
 & = & 200 \bigl(c_1(p) E\bigr)^{1/(p+16)}
 \biggl(\sum_{N=0}^\infty 10^{-N\kappa}\biggr) d^\kappa\nonumber  \\
 & \le & \frac{400}{\kappa} \bigl(c_1(p) E\bigr)^{1/(p+16)}
 d^\kappa \nonumber  \\
 & \le &  \mu_0\ = \ \frac 14 \left(\frac \pi3 - \alpha_0^\ast\right)
\, , \label{small-angle}
 \end{eqnarray}
 where $\alpha_0^\ast \in (0,\frac \pi 3)$ is given by
 \eqref{alphazero}. (We have used $\sum 10^{-j\kappa} =
 10^\kappa/(10^\kappa-1) \le 2/\kappa$ in the second inequality
 above.) In particular $\eps_N\ll 1$ for all $N\in\N.$
\heikodetail{

 \bigskip

 Look at
 $$
 \frac{10^x}{10^x-1}\le \frac{2}{x}, \quad\textnormal{\rm or}\quad
 x10^x\le 2(10^x-1)\Fo x\in (0,1),
 $$
 compute derivatives to find that the left-hand side's derivative
 is dominated by the right-hand side's, and notice that
 LHS(0) equals RHS(0).

 \bigskip

}%

Proceeding inductively, we shall define two sequences of points
$y_N,z_N\in \Sigma$ which converge to $x=0$ and satisfy the
following conditions for each $N=1,2,\ldots$.

\begin{gather}
\frac{d_N}2 \le |y_N|, \, \, |z_N|\le \frac{3d_N}2\, . \label{yNzN} \\[4pt]
\mbox{An initial plane $P_0$ and planes $P_N= \langle
0,y_N,z_N\rangle$ satisfy} \quad
\alpha_{N}:= \ang(P_{N},P_{N-1}) \le 200 \eps_N. \label{alpha-N} \\
\mbox{The angle $\gamma_N:= \ang (y_N, z_N)\in [0,\pi]$ satisfies
} \Bigl|\gamma_N-\frac \pi 2\Bigr|\le
6\eps_1+40(\eps_1+\cdots+\eps_{N-1})\, .\label{50sum}
\end{gather}
We shall also show that there exists a fixed plane $H$ (given by
an application of Lemma~\ref{triples} at the first step of the
whole construction) through $x$ such that, for each
$N=1,2,\ldots,$
\begin{equation}
\pi_H(B(x,d_N)\cap \Sigma) \supset D_H(x,d_N/2):=B(x,d_N/2)\cap
H\, . \label{bigpi-N}
\end{equation}

Here is a short description of the order of arguments: we first
apply Lemma~\ref{triples} to select $P_0$ and then correct it
slightly to have two points $y_1,z_1$ satisfying \eqref{50sum}.
This is done in Steps~1 and~2 below. Next, proceeding inductively,
we first select $y_{N+1},z_{N+1}$ very close to the intersection
of segments $[0,y_N]$ and $[0,z_N]$ with the boundary of $\partial
B_{d_{N+1}}$ (Step 3). Finally, we estimate the angle $\alpha_N$
(Step 4) and prove that $P=\lim P_N$ does not depend on the choice
of $P_0$ (Step 5).

\medskip\noindent \textbf{Step 1.} For given $x$ and $y$ use
Lemma~\ref{triples} to select $z\in B_d(\Sigma)$ and the plane $H$
satisfying conditions (i)--(iii) of that lemma. (Notice that
$|x-y|=d\le\delta_1(E,p)\le R_0(E,p)<d_s(x_0)$ by our choice
\eqref{delta_1} and \eqref{bound-Mp} in
Proposition~\ref{prop:3.5}, so that Lemma~\ref{triples} is indeed
applicable.)

Let $P_0 =\langle x,y,z \rangle = \langle 0,y,z\rangle$; by (iii),
we have
\begin{equation}
\alpha_0':=\ang (P_0, H) \le \alpha_0^\ast = \frac \pi 2 - \arctan
\frac 1{\sqrt{2}} < \frac\pi{3}. \label{alpha0'}
\end{equation}
Lemma~\ref{flat} gives $ \beta_\Sigma(x,d_1)\le \eps_1\, . $ Set
\begin{equation}
F_0:= \{ z'\in B(0,d_1)\, \colon\, \dist(z', P_0) \le
\eps_1d_1\}=U_{\eps_1d_1}(P_0)\cap B_{d_1}.
\end{equation}
We know that $\Sigma \cap B(x,d_1)\subset F_0$. The goal will be
to prove that one can choose $y_N,z_N$ so that for $P_N:=\langle
x,y_N,z_N \rangle $
\begin{equation}
\label{goal} \Sigma \cap B(x,d_N) \subset F_N :=\{ z'\in
B(0,d_N)\, \colon\, \dist(z', P_N) \le \eps_Nd_N\}=
U_{\eps_Nd_N}(P_N)\cap B_{d_N}
\end{equation}
also for $N=1,2\ldots$, and to provide an estimate for
$\alpha_N=\ang(P_{N},P_{N-1})$ showing that for large $N$ the
center planes of the sets $F_N$ stabilize around a fixed affine
plane.

Note that \eqref{bigpi-N} for $N=1$ follows from
Lemma~\ref{triples}~(ii) since $\sin\varphi_0=1/\sqrt{2}>1/2.$

\medskip\noindent\textbf{Step 2 (choice of $P_1$).} We shall choose
$y_1,z_1$ with $\gamma_1=\ang(y_1,z_1) \approx \frac \pi 2$, and
we shall show that the plane $P_1=\langle 0,y_1,z_1\rangle$
satisfies $\alpha_1=\ang(P_1,P_0)\le 12\eps_1$. To this end,
select a point $x_0\in F_0$ such that
\[
h_0:=\dist(x_0,H) = \max_{\xi\in F_0} \, \dist (\xi, H) >0.
\]
It is clear that $x_0$ exists since $F_0 $ is compact, and that
$x_0\in\partial B_{d_1}$; see Figure~6.

Let $\alpha_0'':=\ang(x_0,P_0)$ denote the angle between $x_0$ and
its orthogonal projection $\pi_{P_0}(x_0)$ onto the plane $P_0$.
We have $\sin \alpha_0'' = \eps_1d_1/d_1 = \eps_1$. Hence,
$\alpha_0''\le (\pi/2)\sin\alpha_0''< 2\eps_1$.

\begin{center}
\includegraphics*[totalheight=6.7cm]{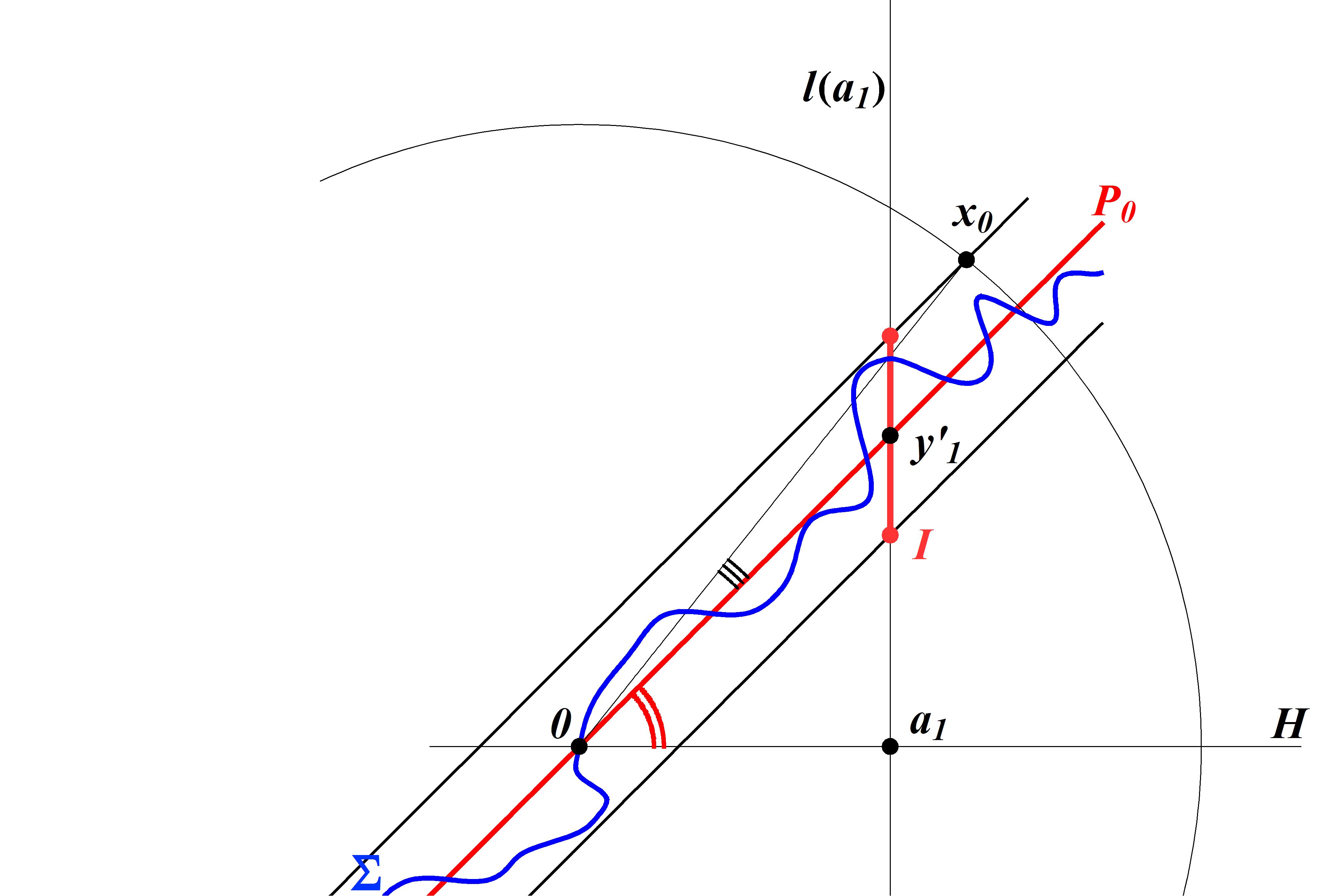}
\end{center}

{

\footnotesize

\medskip\noindent\textbf{Fig. 6.} The initial configuration in
$B(x,d_1)$; cross-section by a plane which is perpendicular to $H$
and $P_0$. A priori, at this stage we do not control the topology
of $\Sigma$ and we cannot even be sure that $\Sigma$ is a graph
over $H$ (or $P_0$). The angle $\alpha_0''$ is marked with a
triple line. \label{fig3}

}

\medskip

Now, since $2\eps_1< 200\sum\eps_N\le \frac 14(\frac \pi 3-
\alpha_0^\ast)$ by \eqref{small-angle}, we can use \eqref{alpha0'}
twice to obtain
\begin{eqnarray}
h_0 = d_1 \sin(\alpha_0'+ \alpha_0'') &
\overset{\eqref{alpha0'}}{<} &
d_1\sin(\alpha_0^*+\frac{1}{4}(\frac{\pi}{3}-\alpha_0^*))\nonumber \\
& = & d_1\sin(\frac{3}{4}\alpha_0^*+\frac{1}{4}\frac{\pi}{3})
\overset{\eqref{alpha0'}}{<} d_1\sin\frac \pi 3 =
d_1\frac{\sqrt{3}}2.\label{h0-est}
\end{eqnarray}
This implies that each straight line $l=l(w)$ which is
perpendicular to $H$ and passes through a point $w$ in the disk
\[
D_0 = D_H(0,r_0)\equiv H\cap B_{r_0}, \qquad \mbox{where $r_0^2+
h_0^2=d_1^2$,}
\]
intersects the finite slab $F_0$ along a segment $I$ of length
$2l_0$, where $\eps_1 d_1/l_0 = \cos\alpha_0'$, which gives $l_0=
(\eps_1d_1)/\cos\alpha_0' < 2\eps_1d_1$ by virtue of
\eqref{alpha0'}. Since $r_0^2=d_1^2-h_0^2 > d_1^2/4$ according to
\eqref{h0-est}, we have $D:=D_H(0, d_1/2)\subset D_0$ in $H$.
Choose two points $a_1,b_1$ in the circle which bounds $D$ in $H$
so that $a_1\perp b_1$ and $b_1\in P_0 \cap H$. Take the lines
$l(a_1)$, $l(b_1)$ passing through these points and perpendicular
to $H$, and select
\begin{equation}
\label{choice-y1z1} y_1 \in \Sigma \cap l(a_1) \cap F_0, \qquad
z_1 \in \Sigma \cap l(b_1) \cap F_0
\end{equation}
(such points do exist since $\Sigma \cap B(x,d_1)\subset F_0$ and
the projection of $\Sigma\cap B(x,d_1)$ onto $H$ contains $D$ by
\eqref{bigpi-N} already verified for $N=1$).

Note that $y_1'$, $z_1'=b_1$ given by
\begin{equation}
\label{y1z1prime} \{ y_1'\} = l(a_1)\cap P_0, \qquad \{ z_1'\} =
l(b_1)\cap P_0
\end{equation}
satisfy $y_1'\perp z_1'$. Let $\psi_0 :=\ang(y_1',y_1)$,
$\theta_0:=\ang(z_1',z_1)$. We have \heikodetail{

\bigskip

For the following tan-estimate look at worst case scenario for the
positions of $y_1\in l(a_1)\cap F_0\cap\Sigma$ and $y_1'\in
l(a_1)\cap P_0$: The angle $\psi$ becomes largest when $y_1$ is
below $y_1'$ on the vertical segment $l(a_1)$ and on the (lower)
boundary of $F_0$; see picture on old p. 27 To estimate the
Ankathete of $\psi$ from below use a rough estimate by $d_1/2-l_0$
which is only equality if $P_0\perp H$ which is excluded by
\eqref{alpha0'}.

\bigskip

}
\begin{eqnarray}
\psi_0\le\tan \psi_0 & \le & \frac{\eps_1 d_1}{({d_1}/{2})-l_0} \nonumber \\
& \le & \frac {2\eps_1 d_1}{d_1-4\eps_1 d_1} \qquad \mbox{as
$l_0\le 2\eps_1d_1$}\label{tan-psi} \\ [4pt] %
& \le & 3 \eps_1 \nonumber
\end{eqnarray}
since $\eps_1 \le \sum\eps_N \le (200)^{-1}\mu_0 \ll 1/12$ by
\eqref{small-angle}. Similarly, we have $\tan \theta_0\le 3
\eps_1$, so that both angles $\psi_0$ and $\theta_0$ do not exceed
$3\eps_1$. Therefore,
$0\le\gamma_1=\ang(y_1,z_1)\le\ang(y_1,y_1')+\ang(y_1',z_1')+\ang(z_1',z_1)$
satisfies
\begin{equation}\label{g0}
\Bigl| \gamma_1-\frac \pi 2\Bigr| \le \psi_0 + \theta_0 \le 6
\eps_1\, ,
\end{equation}
which gives \eqref{50sum} for $N=1.$ By choice of $a_1,b_1$,
\eqref{yNzN} is satisfied for $N=1$. Thus, the triangle
$\Delta=(x,y_1,z_1)$ is $(\eta,d)$-wide, i.e. $\Delta
\in\Spac(\eta,d)$ for $\eta:=\min\{1/2,\pi-(\pi/2+6\eps_1)\}=1/2$
(by \eqref{small-angle}), and $d:=d_1\le R_0(E,p).$ Consequently,
by virtue of \eqref{eps-N} we can derive \eqref{goal} for $N=1$
with the help of Lemma~\ref{flat}.

Finally, normalizing $y_1',z_1'\in P_0$, we easily check that
\begin{equation}
\alpha_1:=\label{angP1P0} \ang(P_1,P_0) < 12 \eps_1\quad\Fo
P_1:=\langle x,y_1,z_1\rangle\, ,
\end{equation}
which gives \eqref{alpha-N} for $N=1$. \heikodetail{

\bigskip

As on old p. 29 (Step 4) we proceed as follows: Set
$u_0:=y_1'/|y_1'|,$
 $w_0:=z_1'/|z_1'|,$
$u_1:=y_1/|y_1|,$
 $w_1:=z_1/|z_1|,$
 and set
\begin{eqnarray*}
 1 & \ge & M_1:=|u_1\times w_1|=|\sin\gamma_1|\ge \frac{1}{\sqrt{2}}, \\
1 & =  & M_0:=|u_0\times w_0|,
\end{eqnarray*}
 since $z_1'\perp y_1'$, and for the first inequality we know by
 \eqref{g0} that $\gamma_1=\ang(u_1,w_1)\in (\frac{5}{12}\pi,\frac{7}{12}\pi)$
 so that $\sin\gamma_1\ge\sin\pi/4=1/\sqrt{2}$.
 Estimate the difference of the unit normal vectors of $P_0=
 \langle x,y_1',z_1'\rangle $ and
 $P_1=\langle x,y_1,z_1\rangle $:
 $$
 \frac{u_1\times w_1}{|u_1\times w_1|}-
 \frac{u_0\times w_0}{|u_0\times w_0|}
 = \frac{(u_1\times w_1)-(u_0\times w_0)}{M_1}+
 \frac{M_0-M_1}{M_1}(u_0\times w_0)=:T_1+T_2
 $$
 as on old p. 29:
 \begin{eqnarray*}
 |T_1| & \le &  \sqrt{2}\Big[|u_1-u_0||w_1|+
 |w_1-w_0||u_1|\Big] =2\sqrt{2}
 \Big[\sqrt{\frac{1-\cos\psi_0}{2}}+\sqrt{\frac{1-\cos\theta_0}{2}}\Big]\\
 & = & 2\sqrt{2}\Big[\sin\psi_0/2+\sin\theta_0/2\Big]
 \le 2\sqrt{2}\Big[\psi/2+\theta_0/2\Big]\\
 & \overset{\eqref{tan-psi}}{\le} & 2\sqrt{2}\cdot 3\eps_1,
\end{eqnarray*}
and
\begin{eqnarray*}
|T_2| & \le & \frac{1-M_1}{M_1}\le\sqrt{2}|\sin\pi/2-\sin\gamma_1|
\\
 & \le & \sqrt{2}|\pi/2-\gamma_1|\overset{\eqref{g0}}{\le}6\sqrt{2}\eps_1.
 \end{eqnarray*}
So, $|T_1|+|T_2|\le 8\sqrt{2}\eps_1 <12\eps_1.$

\bigskip

} Moreover, by \eqref{angP1P0}, \eqref{alpha0'}, and
\eqref{small-angle} we have
\[
\ang(P_1,H) \le \ang(P_1,P_0) + \alpha_0'\quad
\overset{\eqref{angP1P0},\eqref{alpha0'}}{<}\quad
12\eps_1+\alpha_0^* \overset{\eqref{small-angle}}{<}
\frac{1}{4}\left(\frac{\pi}{3}-\alpha_0^*\right)+\alpha_0^*
\overset{\eqref{alpha0'}}{<}\pi/3.
\]
To summarize, we have now proven \eqref{yNzN}, \eqref{alpha-N},
\eqref{50sum}, \eqref{bigpi-N}, and \eqref{goal} for $N=1.$

\medskip\noindent\textbf{Step 3 (induction).} Suppose now that
$y_1, \ldots, y_N$, $z_1, \ldots, z_N$ have already been selected
so that conditions  \eqref{yNzN}, \eqref{alpha-N}, \eqref{50sum},
\eqref{bigpi-N}, and \eqref{goal}
 are satisfied for $j=1,\ldots,N$.
Note that since \eqref{goal} is satisfied for all indices up to
$N$, we have
\begin{equation}
\beta_\Sigma(x,d_j) \le \eps_j = O(d_j^\kappa), \qquad j=1, \ldots
, N.
\end{equation}
We shall select two new points $y_{N+1},z_{N+1}$ such that
\eqref{yNzN}, \eqref{alpha-N}, \eqref{50sum}, \eqref{bigpi-N} and
\eqref{goal} are satisfied with $N$ replaced by $N+1$.

Choose first two auxiliary points,
\begin{equation}\label{aux-points}
\{y_{N+1}'\} := [0,y_N] \cap \partial B(0,d_{N+1})\, ,\qquad
\{z_{N+1}'\} := [0,z_N] \cap \partial B(0,d_{N+1})\, .
\end{equation}
Since $P_N=\langle 0,y_N,z_N\rangle$, we have
$y_{N+1}',z_{N+1}'\in P_N\cap B_{d_{N+1}}\subset F_N$. Fix $x_N\in
F_N$ such that
\[
h_N:=\dist(x_N,H) = \max_{\xi\in F_N} \, \dist (\xi, H).
\]
Set $\alpha_{N}':=\ang (P_N, H)$, $\alpha_{N}'':= \ang(x_N,P_N)$.
We note that $x_N\in \partial B_{d_N}$ and by \eqref{small-angle}
\[
\alpha_{N}''= \arcsin \eps_N < 2 \eps_N \le 2 \eps_1
\overset{\eqref{small-angle}}{<} \frac 14 \left(\frac \pi 3 -
\alpha_0^\ast\right).
\]
Applying the triangle inequality and using the induction
hypothesis \eqref{alpha-N} up to $N$, and \eqref{alpha0'}, we
estimate
\begin{eqnarray}\label{aN-1'}
\alpha_{N}'& = & \ang (P_N, H) \\ [4pt] %
& \le &  \ang (P_0,H) + \ang(P_1,P_0) + \ang(P_2,P_1)
+ \cdots + \ang(P_N,P_{N-1})\nonumber \\[3pt]
& = & \alpha_0' +  \alpha_1 + \cdots + \alpha_{N} \qquad
\nonumber \\
& \le & \alpha_0^\ast +  \frac 14 \left(\frac \pi 3 -
\alpha_0^\ast\right) \qquad \mbox{by \eqref{alpha0'},
\eqref{alpha-N}, and \eqref{small-angle}.}\nonumber
\end{eqnarray}
Thus, $\alpha_{N}' + \alpha_{N}''\le \alpha_0^\ast + \frac 1 2
\left(\frac \pi 3 - \alpha_0^\ast\right)< \frac \pi 3$ and, as in
the second step, we have
\[
h_N = d_N \sin(\alpha_{N}'+ \alpha_{N}'') < d_N\sin\frac \pi 3 =
d_N\frac{\sqrt{3}}2, \qquad N\ge 1\, .
\]
Hence, $d_N^2 = h_N^2 + r_N^2$ for some $r_N > d_N/2$; as
previously, we conclude that each straight line $l=l(w)$ which is
perpendicular to $H$ and passes through a point $w$ in the disk
\[
D_N = D_H(0,r_N)\equiv H\cap B_{r_N},
\]
intersects the finite slab $F_N$ along a segment $I$ of length
$2l_N$, where $\eps_N d_N/l_N = \cos\ang(P_N,H)$, which gives $l_N
< 2\eps_N d_N$ by virtue of \eqref{aN-1'}. Moreover, by
\eqref{bigpi-N} (which, by the inductive assumption, holds for
$N$), each segment $I(w)$ for $w\in D_H(0,d_N/2)$ vertical to $H$
must contain at least one point of $\Sigma$.

We now choose $y_{N+1},z_{N+1} \in F_N \cap \Sigma$ such that
\begin{equation}
\pi_H(y_{N+1}) = \pi_H(y_{N+1}')\, , \qquad \pi_H(z_{N+1}) =
\pi_H(z_{N+1}')\, .
\end{equation}
To establish the desired estimate of $\ang(P_{N+1},P_N)$, we show
first that
\begin{gather}
\psi_N :=\ang (y_{N+1},y_N) \le 20 \eps_N,\label{psi-N} \\[4pt]
\theta_N :=\ang (z_{N+1},z_N) \le 20 \eps_N. \label{theta-N} \\[4pt]
\end{gather}
Indeed, \heikodetail{

\bigskip

For the following tan-estimate look at worst case scenario for the
positions of $y_{N+1}\in l(\pi_H(y_{N+1}'))\cap F_N\cap\Sigma$ and
$y_{N+1}'\in l(\pi_H(y_{N+1}'))\cap P_N$: The angle $\psi_N$
becomes largest when $y_{N+1}$ is below $y_{N+1}'$ on the vertical
segment $l(\pi_H(y_{N+1}'))$ and on the (lower) boundary of $F_N$;
see picture on old p. 28 To estimate the Ankathete of $\psi_N$
from below use a rough estimate by $d_{N+1}-l_N$ which is only
equality if $P_N\perp H$ which is excluded by \eqref{aN-1'}.

\bigskip

}
\begin{eqnarray*}
\tan\psi_N & = & \tan \ang (y_{N+1},y_N)\ \overset{\eqref{aux-points}}{=}\ %
\tan \ang (y_{N+1},y_{N+1}') \\
& < & \frac{\eps_N d_N}{d_{N+1}-l_N} \\
& < & \frac{\eps_N d_N}{d_{N+1}/2} \ =\   20 \eps_N,
\end{eqnarray*}
where the last inequality holds since $l_N < 2\eps_Nd_N \le 2
\eps_1 d_N < d_N/300 < d_{N+1}/2$; remember that $2\eps_1 \le 2
\sum_N \eps_N\le (100)^{-1} \mu_0 <(100)^{-1}\pi/12<1/300$ by
\eqref{small-angle}.

Thus, $\psi_N\le \tan\psi_N \le 20\eps_N$. Similarly, $\theta_N\le
\tan\theta_N \le 20\eps_N$. This proves \eqref{psi-N} and
\eqref{theta-N}. Moreover, the triangle inequality gives an
estimate of the angle $\gamma_{N+1}=\ang(y_{N+1},z_{N+1})$,
\begin{equation}
|\gamma_{N+1}-\gamma_N| \le \theta_N + \psi_N \le 40 \eps_N,
\label{g-NN+1}
\end{equation}
\heikodetail{

\bigskip

Triangle inequality for angles gives
$$
\ang(y_{N+1},z_{N+1})\le \ang(y_{N+1},y_{N}) + \ang(y_{N},z_{N})+
\ang(z_{N},z_{N+1})
$$
and also
$$
\ang(y_{N},z_{N})\le \ang(y_{N},y_{N+1}) + \ang(y_{N+1},z_{N+1})+
\ang(z_{N+1},z_{N})
$$

\bigskip

} and consequently
\[
\Bigl|\gamma_{N+1}-\frac \pi 2\Bigr| \, \le \,
\Bigl|\gamma_{N}-\frac \pi 2\Bigr| + 40\eps_N.
\]
By induction, this inequality implies \eqref{50sum} with $N$
replaced by $N+1$.
We also have
\[
\frac{d_{N+1}}{2}\le d_{N+1} - l_N \le |y_{N+1}| \le d_{N+1} + l_N
\le \frac{3d_{N+1}}{2},
\]
\heikodetail{

\bigskip

see worst case picture on old p. 28 bottom as for the tan-estimate

\bigskip

}%
and a similar estimate for $|z_{N+1}|$, which gives \eqref{yNzN}
with $N$ replaced by $N+1$. Therefore  the triangle
$\Delta=(x,y_{N+1}, z_{N+1})$ is $(\eta_1, d_{N+1})$-wide, i.e.
$\Delta\in\Spac(\eta_1,d_{N+1})$ for
$\eta_1:=\min\{1/2,(\pi/2)-50\sum_N\eps_N\}=1/2$ according to
\eqref{small-angle}. Since $\eta_1 = 1/2 >\eta={}$ the constant
from Theorem~\ref{goodtetra},  Lemma \ref{flat} is again
applicable to
 obtain
\[
\Sigma \cap B_{3d_{N+1}} \subset U_{\eps_{N+1}d_{N+1}}(P_{N+1})\,
,
\]
which implies \eqref{goal} with $N$ replaced by $N+1$.

To check \eqref{bigpi-N} with $N+1$ instead of $N$, we fix $z\in
\Sigma \cap (B_{d_N}\setminus B_{d_{N+1}})$ and estimate
$|\pi_H(z)|$. Since then $z\in F_N = U_{\eps_Nd_N}(P_N) \cap B_N$
and the angle $\alpha_N'= \ang(P_N,H)$ satisfies \eqref{aN-1'}, we
check that $\sin(\ang(z,P_N))\le\eps_Nd_N/|z|\le\eps_Nd_N/d_{N+1}$
and consequently
\begin{eqnarray*}
|\pi_H(z)| & = & |z|\cos \ang(z,H) \\
& \ge & |z| \cos \left(\alpha_N' + \arcsin
\frac{\eps_Nd_N}{d_{N+1}}\right) \\
& > & |z| \cos (\alpha_N' + 20\eps_N) \qquad \mbox{as
$d_N=10d_{N+1}$} \\
& > & |z| \cos \frac \pi 3 \qquad \mbox{by \eqref{small-angle}
and \eqref{aN-1'}} \\
& \ge & d_{N+1}/ 2\, .
\end{eqnarray*}
Since by the inductive assumption (\eqref{bigpi-N} up to index
$N$) the projection $\pi_H(\Sigma \cap B_{d_N})$ contains the
whole disk $D_H(0,d_N/2)$,  we do obtain $\pi_H(\Sigma \cap
B_{d_{N+1}}) \supset D_H(0,d_{N+1}/2)$. \heikodetail{

\bigskip

The argument is:

\begin{eqnarray*}
D_H(0,d_N/2)\subset \pi_H(\Sigma\cap B_{d_N}) & = &
\pi_H(\Sigma\cap (B_{d_N}\setminus B_{d_{N+1}}))\cup
\pi_H(\Sigma\cap B_{d_{N+1}}) .
\end{eqnarray*}
The just established estimate implies that for all
$\xi\in\pi_H(\Sigma\cap (B_{d_N}\setminus B_{d_{N+1}}))$ one has
$\xi\not\in D_H(0,d_{N+1}/2)$, which implies
$D_H(0,d_{N+1}/2)\subset\pi_H(\Sigma\cap B_{d_{N+1}}).$

\bigskip

}

It remains to verify \eqref{alpha-N} with $N$ replaced by $N+1$,
i.e., the desired inequality for the angle $\alpha_{N+1}=\ang
(P_{N+1}, P_{N})$.

\medskip\noindent\textbf{Step 4. Estimates of $\alpha_N$.}
We normalize the vectors spanning $P_j$ and set $u_j :=
{y_j}/{|y_j|}$, $w_j := {z_j}/{|z_j|}\,$. We also set $M_j = |u_j
\times w_j|$, noting that by \eqref{50sum} which we already have
shown to hold up to $N+1$, and by \eqref{small-angle} that
$$
\gamma_j\in (\frac{5}{12}\pi,\frac{7}{12}\pi)
$$
so that
\begin{equation}
1 \ge M_j = \sin \gamma_j \ge \frac{\sqrt{2}}2  \quad\Foa
j=1,\ldots,N+1.
 \label{MN}
\end{equation}
Now, we compute the difference of unit normals to $P_{N+1}$ and
$P_N$,
\[
\frac{u_{N+1} \times w_{N+1}}{|u_{N+1} \times w_{N+1}|} -
\frac{u_N \times w_N}{|u_N \times w_N|} =: T_1 + T_2,
\]
where
\begin{eqnarray*}
T_1 & := &
\frac{M_N(u_{N+1} \times w_{N+1} - u_N \times w_N)}{M_N M_{N+1}}\, ,\\
T_2 & := & \frac{M_N-M_{N+1}}{M_N M_{N+1}}\, u_N \times w_N\, .
\end{eqnarray*}
Since $u_N,w_N\in \stwo$, we can use \eqref{psi-N},
\eqref{theta-N} (which yield the estimates of $u_{N+1}-u_N$ and
$w_{N+1}-w_N$), and in addition \eqref{MN} and \eqref{g-NN+1}, to
obtain
\begin{eqnarray*}
|T_1| & \overset{\eqref{MN}}{\le} & \sqrt{2} |u_{N+1} \times w_{N+1} - u_N \times w_N)|\\
& \le & \sqrt{2} \bigl(|u_{N+1} - u_N | + |w_{N+1} - w_N)|\bigr)
\\
& \le & 40 \sqrt 2 \eps_N \ < \ 60 \eps_N\, ; \\
|T_2| & \overset{\eqref{MN}}{\le} & 2 |\sin \gamma_N - \sin
\gamma_{N+1}|\qquad\textnormal{since
$\sin$ is $1$-Lipschitz}\\
& \le & 2 |\gamma_N -\gamma_{N+1}| \overset{\eqref{g-NN+1}}{\le}
80\eps_N\, .
\end{eqnarray*}
This implies
\begin{equation}
\label{140} \alpha_{N+1}=\ang(P_{N+1},P_N) \le 140\eps_N\, ,
\end{equation}
i.e., \eqref{alpha-N} holds also with $N$ replaced by $N+1$.

Finally, a computation similar to \eqref{aN-1'} shows that
\begin{equation}\label{angle-H}
\ang(H,P_{N+1})< \pi/3\, .
\end{equation}
To summarize,  under the inductive hypothesis that \eqref{yNzN},
\eqref{alpha-N}, \eqref{50sum}, \eqref{bigpi-N}, and \eqref{goal}
hold up to $N$, we have shown  that \eqref{yNzN}, \eqref{alpha-N},
\eqref{50sum}, \eqref{bigpi-N}, and \eqref{goal} do hold up to
$N+1$,  which yields \eqref{yNzN}--\eqref{bigpi-N} and
\eqref{goal} for all $N\in\N$ by the induction principle.

\medskip\noindent\textbf{Step 5 (existence and uniqueness of $\lim P_N$).}
The unit vectors $u_N=y_N/|y_N|$ and $w_N=z_N/|z_N|$ spanning the
affine planes $P_N$ with unit normals $\nu_N:=u_N\times w_N$
satisfy $\ang(u_N,w_N)=\gamma_N\in (\frac{5}{12}\pi,
\frac{7}{12}\pi)$ for all $N$, so that  subsequences again denoted
by $u_N$ and $w_N$ converge to unit vectors $u,w\in\S^1$ with
$\ang(u,v)\in [\frac{5}{12}\pi, \frac{7}{12}\pi]$ spanning a
limiting affine  plane $P$ with unit normal vector $\nu:=u\times
w$, so that we can say $P_N\to P$ as $N\to\infty.$ Since all $P_N$
contain $x=0$ so does $P$. As in \eqref{small-angle}, summing the
tail of a geometric series, we obtain by \eqref{alpha-N}:
\begin{eqnarray}
\ang(P,P_N)=\lim_{k\to\infty}\ang(P_k,P_N) &   \le &
\sum_{j=N}^\infty \alpha_{j+1}
 \le  200\sum_{j={N+1}}^\infty \eps_{j} \nonumber \\
& \overset{\eqref{alpha-N}}{\le} & \frac{400}{\kappa} \bigl(c_1(p)
E\bigr)^{1/(p+16)}
d_{N+1}^\kappa \nonumber  \\
& =: & C_2(p) E^{1/(p+16)} d_{N+1}^\kappa\, \quad\Foa
N=0,1,2,\ldots \label{speed}
\end{eqnarray}
In particular,
\begin{equation}
\ang(P,P_0) \le C_2(p) E^{1/(p+16)} d_1^\kappa \equiv C_2(p)
E^{1/(p+16)} d^\kappa\, .\label{speed2}
\end{equation}

However, as we cannot a priori claim that $\Sigma$ is a graph over
$H$, the choice of $y_N$ and $z_N$ for small values of $N$ does
not have to be unique. Suppose that for two different choices of
sequences $y_N,z_N\in\Sigma$ and $y_N',z_N'\in\Sigma$ (satisfying
\eqref{yNzN}--\eqref{50sum}, and \eqref{goal} for all $N\in\N$),
we obtain
\[
P_N = (x,y_N,z_N)\to P, \qquad P'_N = (x,y_N',z_N')\to P' \quad\As
N\to\infty,
\]
but $P\not=P'$ and $\pi/2\ge \ang(P,P') = \vartheta >0$. Fix $N$
so large that $\eps_N< \vartheta/10$ and
\[
\max\bigl(\ang( P,P_N),\ang(P',P_N')\bigr) <\vartheta/10\, .
\]
Since $y_N'\in P_N'$ and $d_N/2 \le |y_N'|\le 3d_N/2$
by\eqref{yNzN}, we obtain $\ang (y_N',P') < \vartheta/10$. Hence,
the angle between $y_N'$ and $P_N$ cannot be too small: $\ang
(y_N', P) \ge \ang(P',P)-\ang(y_N',P')>9\vartheta/10$ and $\ang
(y_N',P_N)\ge \ang(y_N',P)-\ang(P,P_N)
>4\vartheta /5$. Therefore,
\[
\dist (y_N',P_N) = |y_N'|\, \sin\ang(y_N',P_N) > \frac{d_N}2 \cdot
\frac 2\pi\cdot \frac{4\vartheta}{5} >\frac{d_N\vartheta}{5}
 > 2 \eps_N d_N,
\]
which is a contradiction to
\[
\Sigma \cap B(x,3d_N) \subset U_{\eps_Nd_N}(P_N),
\]
as $|y_N'|\le 3d_N /2 < 3d_N$. Thus, $P=\lim P_N$ is unique and
does not depend on the choices of $y_N,z_N$.

We set  $P=:x + T_x\Sigma$ to define the tangent plane $T_x\Sigma$
of $\Sigma $ at the point $x$, and we set $n_\ast(x):=\nu$ to
obtain a well-defined unit  normal to $\Sigma $ at $x$; the
estimate \eqref{speed} gives in fact
\eqref{distance_from_tangent-plane} (justifying the term ``tangent
plane'')
\begin{eqnarray}
\dist (x', P) & \le & 2\eps_Nd_N+d_N\sin\ang(P,P_N)\nonumber \\
& = & E^{1/(p+16)}\, O(d_N^{1+\kappa}), \quad N\to \infty, \quad
\mbox{for all $x'\in B(x,d_N)\cap \Sigma$,} \label{tan-box}
\end{eqnarray}
where the constant in `big O' above depends only on $p$.
\heikodetail{

\bigskip

see back of old p. 29 for a picture

\bigskip

}%

\medskip\noindent\textbf{Step 6 (conclusion of the proof).}
Reversing the roles of $x$ and $y$, running the whole procedure
one more time, and using \eqref{speed2} twice, we obtain
\begin{equation}
\label{angle-TxTy} \ang(T_x \Sigma, T_y\Sigma)\ \le\ \ang(T_x
\Sigma, P_0) + \ang (P_0, T_y \Sigma)\ \le\ 2C_2(p)E^{1/(p+16)}
d^\kappa\, .
\end{equation}

\smallskip

\noindent $\Box$

\medskip

We state one corollary which easily follows from the last result
and its proof. It tells us that it is not really important how we
choose $P_0$; there are many choices which give a similar
approximation of $T_x \Sigma$.

\begin{corollary}
\label{tan-sec} Assume that $\Sigma\in \A$ and $\M_p(\Sigma) < E$
for some $p>8$. Let $T_x\Sigma$ and $\delta_1=\delta_1(E,p)
>0$ be given by Theorem~\ref{thm:5.3}.

Whenever $x,y,\zeta\in \Sigma$ with $0<d=|x-y|\le \delta_1(E,p)$,
$d/2\le |x-\zeta|\le d$ and $\ang(\zeta-x,y-x)\in [\pi/3,
2\pi/3]$, then $T_x\Sigma$ and the plane $P=\langle
x,y,\zeta\rangle$ satisfy
\begin{equation}
\ang (T_x\Sigma, P) \le C_3(p) E^{1/(p+16)} d^\kappa, \qquad
\kappa= \frac{p-8}{p+16}\, , \label{any-P}
\end{equation}
where the constant $C_3(p)$ depends only on $p$.
\end{corollary}

\noindent\textbf{Proof.} We use the notation introduced in the
proof of Theorem~\ref{thm:5.3}. Since $\ang(T_x\Sigma,P_0)\lesssim
E^{1/(p+16)} d^\kappa$ by \eqref{speed2}, it is enough to show
that the angle $\ang(P_0,P)$ does not exceed a constant multiple
of $E^{1/(p+16)} d^\kappa$. Noting that $d/2\le |\zeta-x|\le
d=d_1$ and $\zeta$ belongs to the slab $U_{\eps_1d_1}(P_0)$, we
easily compute this angle and finish the proof. The computational
details, very similar to the proof of \eqref{140}, are left to the
reader. \qed \heikodetail{

\bigskip

for example one argues as follows for the angle-estimate:
$$
\sin\ang(P.P_0)=\frac{\dist(\zeta,P_0)}{|\zeta|}\le\frac{\eps_1d_1}{
d_1/2}=2\eps_1\cong E^{\frac{1}{p+16}}d^\kappa
$$
by choice of $\eps_1$ in \eqref{eps-N}.

\medskip

}%

In order to deal with sequences of surfaces with equibounded
energy in Section \ref{sec:7} we establish a local graph
representation of one such surface  $\Sigma$ of finite
$\M_p$-energy on a scale completely determined by the energy value
$\M_p(\Sigma)$ and with a priori estimates on the
$C^{1,\kappa}$-norm of the graph function.
\begin{corollary}\label{cor:5.4new}
Assume that $p>8$, $\M_p(\Sigma)< E<\infty.$ Then there exist two
constants, $0<a(p)<1<A(p)<\infty$, such that for each $x\in\Sigma$
there is a function
$$
f:T_x\Sigma\to \Big( T_x\Sigma\Big)^\perp \simeq\R
$$
with the following properties:
\begin{enumerate}
\item[\rm (i)]
$f(0)=0,$ $\nabla f(0)=(0,0)$,
\item[\rm (ii)]
$|\nabla f(y_1)-\nabla f(y_2)|\le
A(p)E^{\frac{1}{p+16}}|y_1-y_2|^\frac{p-8}{p+16}$,
\item[\rm (iii)]
If $R_1\equiv R_1(E,p):=a(p)E^{-1/(p-8)}\le R_0(E,p)$ (where
$R_0(E,p)$ has been defined in \eqref{R0Ep} of Theorem
\ref{Thm:3.1}) and if
$$
\Phi(y):=x+(y,f(y)),\quad y\in T_x\Sigma\simeq\R^2,
$$
then
\begin{equation}\label{graph-patch-size}
\Phi(D_{\frac{3}{4}R_1})\subset
\Big[B(x,R_1)\cap\Sigma\Big]\subset\Phi(D_{R_1}),
\end{equation}
where $D_{R_1}=B(0,R_1)\cap T_x\Sigma$ is a disk in $T_x\Sigma$
around $0\in T_x\Sigma,$ and
\begin{equation}\label{graph-estimate}
|D\Phi(y_1)-D\Phi(y_2)|\le
A(p)E^{\frac{1}{p+16}}|y_1-y_2|^\frac{p-8}{p+16}.
\end{equation}
\end{enumerate}
In particular, $\Sigma$ is an orientable $C^{1,\kappa}$-manifold
for $\kappa=(p-8)/(p+16).$
\end{corollary}
\noindent\textbf{Proof.} Basically, we mimic here the proof of
Theorem 5.1 from \cite{svdm-global}. (In \cite{svdm-global}, we
knew that the surface cannot penetrate two balls of \emph{fixed}
radius, touching $\Sigma$ at every point; this is replaced here by
angle estimates \eqref{speed} and \eqref{speed2}, and the
existence of forbidden conical sectors, cf.~
Proposition~\ref{sectors}.)

Fix $x\in \Sigma$. Without loss of generality suppose that $x=0$.

\smallskip\noindent\textbf{Step 1 (the definition of $f$).} We use
the notation from the proof of Theorem~\ref{thm:5.3}. Recall the
the plane $P=x+T_x\Sigma$ (used to \emph{define\/} $T_x\Sigma$)
has been obtained as a limit of planes $P_N$ satisfying
\eqref{speed}; for all $x,y\in \Sigma$ with $|x-y|=d\le
\delta_1(E,p)$ given by \eqref{delta_1} we had the angle estimate
\eqref{angle-TxTy}. Using \eqref{tan-box}, one can easily show
that \heikodetail{

\bigskip

choosing $N$ such that $d_{N+1}<d\le d_N$ in \eqref{tan-box}

\bigskip

}
\begin{equation}
\label{tan-box2} \dist (x',P) \le A_1(p) E^{1/(p+16)}
d^{1+\kappa}\,
\end{equation}
whenever $x'\in B(x,d)\cap \Sigma$ for some $d\le \delta_1(E,p)$.
We shall use this estimate and Proposition~\ref{sectors} to
show that if $r\le a(p) \delta_1(E,p)$ for a sufficiently small
constant $a(p)\in (0,1)$, then
\begin{equation}
\left( \pi_P ( B(x,{4r/3}) \cap \Sigma\right)
)\qquad\mbox{contains the disk $D_r:=B(x,r)\cap P$.}\label{onto}
\end{equation}
Indeed, otherwise there would be a point $z\in D_r$ and a segment
$I=I_{h,w}(z)\perp P$ (we fix a unit vector $\Sphere^2\ni w\perp
P$) of length
\begin{eqnarray*}
2h &:=& 2A_1(p) E^{1/(p+16)} (4r/3)^{1+\kappa} \\
&\le & \frac{r}{100} \qquad\mbox{if $a(p)$ is small enough}
\end{eqnarray*}
such that $I\cap \Sigma=\emptyset$. By \eqref{tan-box2}
all points of $\Sigma$ in $B(x,d)$, $d=4r/3$, are in fact located
in the thin slab $U_h(P)$. Thus, it is easy to use
Proposition~\ref{sectors}, \eqref{C+in}--\eqref{C-out}, and check
that --- no matter what is the angle between $P$ and the vector
$v$ given by that Proposition --- the sets
$C^\pm_{2r}(\varphi_0,v) \setminus B_r$ contain two open balls
$B^\pm$ which are in two \emph{different\/} components of
$B(x,d)\setminus U_h(P)$. Hence,
\[
B^+ \subset C^+_{2r}(\varphi_0,v) \cap U, \qquad B^-\subset
C^-_{2r}(\varphi_0,v) \cap (\bbbr^3\setminus \overline{U})
\]
Now, one could use the segment $I$  to construct a curve which
contains no point of $\Sigma$ but nevertheless joins a point in
$B^-$ to a point in $B^+$. This is a contradiction proving
\eqref{onto}.

Next, using \eqref{angle-TxTy}, one proves that $\pi_P$ is
injective on $B(x,{4r/3})\cap P$. Otherwise, there would be a
point $z'\in P$, $4r/3> |z'-x|=\rho>0$, and a segment
$I':=I_{h',w}(z')$ with
\[
h'= A_1(p) E^{1/(p+16)} \rho^{1+\kappa} \le \rho/100
\]
such that $I'\cap\Sigma$ would contain two different points
$y_1\not=y_2$. Then, letting $P_1=T_{y_1}\Sigma$,
$v_1=(y_1-x)/|y_1-x|$ and $v_2=(y_2-y_1)/|y_2-y_1|$, we would use
\eqref{angle-TxTy} to obtain
\begin{eqnarray*}
\ang(v_1,v_2) & \le & \ang (v_1, P) + \ang (P,P_1) + \ang
(P_1,v_2) \\
& \le & A_2(p) E^{1/(p+16)} \rho^\kappa \\
& < & \frac\pi 4  \qquad\mbox{if $a(p)$ is small enough.}
\end{eqnarray*}
\heikodetail{

\bigskip

See my drawing on old p. 48 and for the previous estimate use
$$
\sin\ang(v_1,P)\simeq \frac{\dist(y,P)}{\rho}\simeq \rho^{\kappa}
$$
and for the next estimate similarly
$$
\sin\ang(v_1,P)\simeq\rho^{\kappa}.
$$

\bigskip

}
Since $v_2\perp P$ we have on the other hand for sufficiently
small $a(p)$
\[
\ang(v_1,v_2)\ge \frac \pi 2 - \ang (P,v_1) \ge \frac \pi 2 -
A_3(p) E^{1/(p+16)} \rho^\kappa > \frac \pi 4,
\]
a contradiction.

For $y\in U$, where $U$ denotes the interior in $P$ of
$\pi_P(\Sigma \cap B(x,{4r/3})$ we now define
\[
f(y) = w\cdot \left(\Bigl.\pi_P\Bigr|_{\Sigma\cap
B(x,{4r/3}}\right)^{-1}(y)\, ,
\]
and let $\Phi(y)$ be defined by the formula given in Part (iii) of
the Corollary. Note that $U\supset D_r$ by \eqref{onto}. It is
clear that $f(0)=0$ and $\nabla f(0)=(0,0)$. \heikodetail{

\bigskip

the latter and the following stem from looking at how the graph
is squeezed in between graphs of $C^{1,\kappa}$ functions, which
we know already by estimates like \eqref{tan-box2}. Do not
use the explicit formula defining $f$ to calculate this! See
my sketches on my old pp. 49,50.

\bigskip

}
The differentiability of
$f$ at other points follows from \eqref{tan-box2} which implies
that for $\varrho\to 0$ $\mbox{Graph}\, f \cap B(x,\varrho)$ is
trapped in a flat slab of height ${}\lesssim \varrho^{1+\kappa}$
around a \emph{fixed\/} plane (depending on $x$ but independent
from $\varrho$).

\smallskip\noindent\textbf{Step 2 (bounds for $|\nabla f|$).}
The vector $(\nabla f(y),-1)$ is parallel to the normal direction
to $\Sigma$ at $x$ when $y=\pi_P(x)$.  Taking $y\in U$, we have by
\eqref{osc-normal} of Theorem~\ref{thm:5.3}
\[
\alpha(y) \equiv\ang (T_{\Phi(y)}\Sigma, T_{0}\Sigma) \le \pi/4
\]
Since $\tan\alpha(y)=|\nabla f(y)|$, we have $|\nabla f(y)|\le 1$
everywhere in $D_r$. Thus, $f$ is Lipschitz with Lipschitz
constant 1.

\smallskip\noindent\textbf{Step 3 (the oscillation of  $\nabla f$).}
 Fix two points $y_1,y_2\in
U$ and set $a=D_1f(y_1)$, $b=D_2f(y_1)$, $c=D_1f(y_2)$,
$d=D_2f(y_2)$ where $D_i$ stands for the $i$-th partial
derivative. The angle $\alpha$ between the tangent planes to
$\Sigma$ at $x_i=\Phi(y_i)$, $i=1,2$, satisfies
\begin{eqnarray}
\sin^2\alpha &=& \frac{(a-c)^2+(b-d)^2+(ad-bc)^2}
{(1+a^2+b^2)(1+c^2+d^2)}
\nonumber  \\
&\stackrel{\mbox{\scriptsize (Step 2)}}{\ge}& \frac
{(a-c)^2+(b-d)^2} {4}=\frac{|Df(y_1)-Df(y_2)|^2}{4}\, .
\label{sinestim}
\end{eqnarray}
\heikodetail{

\bigskip

simply use the cross product to express $\sin\alpha$

\bigskip

} An upper bound for $\alpha$ is also given by \eqref{osc-normal}.
Combining the two, and noting that $|x_1-x_2|\le 2|y_1-y_2|$, we
obtain the desired estimate for $y_1,y_2\in U$ and conclude the
proof, extending $f$ to the whole tangent plane by well-known
extension theorems; see e.g. \cite[Chapter 6.9]{GT98}. \hfill
$\Box$

\begin{REMARK}\label{rem:5.4newA}
 \rm Assume that some absolute small constant $\eps_0$ is given a
priori, say $\eps_0=\frac{1}{100}. $ Then, shrinking $a(p)$ in the
previous corollary if necessary, we have above  for $y_1,y_2 \in
D_{R_1}$
\begin{eqnarray*}
|\nabla f(y_1)-\nabla f(y_2)| & \le &
A(p)E^{\frac{1}{p+16}}|y_1-y_2|^\frac{p-8}{p+16}
 \le   A(p)E^{\frac{1}{p+16}}(2R_1)^\frac{p-8}{p+16}\\
& \le &  2 A(p)E^{\frac{1}{p+16}}a(p)^\frac{p-8}{p+16} \Big(
E^{-1/(p-8)}\Big)^\frac{p-8}{p+16}  =  2A(p)a(p)^\frac{p-8}{p+16}
<\eps_0.
\end{eqnarray*}
\end{REMARK}

\begin{REMARK}\label{orient} \rm It is now clear that if
$\Sigma \in \A$ with $\M_p(\Sigma)<\infty$ for some $p>8$, then
$\Sigma=\partial U$ is a closed, compact surface of class
$C^{1,\kappa}$. Thus, $\Sigma$ is orientable and has a well
defined global normal, $n_\Sigma$.

For a discussion of issues related to orienta\-bility, we refer
the reader to \cite{lima} and to Dubrovin, Fomenko and Novikov's
monograph,  \cite[Chapter 1]{dnf}.
\end{REMARK}


\section{Improved H\"{o}lder regularity of the Gau{\ss} map}
\label{sec:6} \setnumbers

In this section we prove
\begin{theorem}
Let $\Sigma\in \A$; assume that $p>8$ and $\M_p(\Sigma)\le E <
\infty$. Then $\Sigma$ is  an orientable manifold of class
$C^{1,\lambda}(\Sigma)$ for $\lambda = 1-\frac 8p$.
\label{improved} Moreover, the unit normal $n_\Sigma$ satisfies
the local estimate
\begin{equation}\label{normal-hoelder-estimate}
|n_\Sigma (x_1)-n_\Sigma(x_2)|\le C(p) \biggl( \int_{[\Sigma\cap
B(x_1,10|x_1-x_2|)]^4} \K^p \, d\mu \biggr)^{1/p}|x_1-x_2|^\lambda
\end{equation}
for all $x_1,x_2\in\Sigma$ such that $|x_1-x_2|\le \delta_2(E,p):=
a_2(p) E^{-1/(p-8)}$.

\end{theorem}

\medskip\noindent\textbf{Remark.} Once
\eqref{normal-hoelder-estimate} is established, the global
estimate $|n_\Sigma(x_1)-n_\Sigma(x_2)|\le \mathrm{const}
|x_1-x_2|^\lambda$ follows.

Before passing to the proof of the theorem, let us explain
informally what is the main qualitative difference between the
estimates in Sections 5 and 6. In Section 5, to prove that the
surface is in fact $C^{1,\kappa}$, we were iteratively estimating
the contribution to the energy of tetrahedra with vertices on
patches that were very small when compared with the edges of those
tetrahedra. A priori, this might be a tiny fraction of
$\M_p(\Sigma)$. Now, knowing already that locally the surface is a
(flat) $C^{1,\kappa}$ graph, we can use a slicing argument to
gather more information from energy estimates
--- this time, considering not just an insignificant portion of the
local energy but the whole local energy to improve the estimates
of the oscillation of the normal vector.

The whole idea is, roughly speaking, similar to the proof of
Theorem~1.2 in our joint paper with Marta Szuma\'{n}ska, see
\cite[Section 6]{ssvdm-triple}. Since the result is local, we
first use Theorem~\ref{thm:5.3} to consider only a small piece of
$\Sigma$ which is a (very) flat graph over some plane, and then we
use the energy to improve the H\"{o}lder exponent from
$\kappa=(p-8)/(p+16)$ to $\lambda=1-\frac 8p>\kappa$.

\medskip\noindent\textbf{Proof of Theorem \ref{improved}. \, Step 1.
The setting.} W.l.o.g. we consider a portion of $\Sigma$ which is
a graph of $f\colon \R^2 \supset 5Q_0 \to \R$, where $Q_0$ is some
fixed (small) cube centered at $0$ in $\R^2$ and $f\in
C^{1,\kappa}$ satisfies $\nabla f(0)=(0,0)$ and has a very small
Lipschitz constant, say
\begin{equation}
\label{lip-f} |f(x)-f(y)| \le \eps_0 |x-y|, \qquad \qquad x,y \in
5 Q_0\, .
\end{equation}
By an abuse of notation, we write $n_\Sigma(x)$ to denote the
normal to $\Sigma$ at the point $F(x)\in \Sigma$, where
\begin{equation}
\label{def-F} F\colon \R^2 \supset 5Q_0\ni x \longmapsto (x,f(x))
\in \bbbr^3\,
\end{equation}
is the local parametrization of $\Sigma$ given by the graph of
$f$, compare with Corollary \ref{cor:5.4new}. To ensure
\eqref{lip-f}, just use Remark~\ref{rem:5.4newA}.

We shall write $\K(x_0,x_1,x_2,x_3)$ to denote the integrand of
$\M_p$ (without the power $p$) evaluated at the tetrahedron with
four vertices $F(x_i)\in \Sigma$ for $x_i$ in the domain of the
para\-me\-trization $F$.

Since \eqref{lip-f} implies that $|\nabla f|\le \eps_0$, we  also
have $|F(x)-F(y)|\le (1+\eps_0)|x-y|^\kappa$; hence
\begin{equation}
(1+\eps_0)^2 \H^2(U) \ge \H^2 \bigl(\Sigma \cap F(U)\bigr) \ge
\H^2(\pi_{\R^2}(\Sigma\cap F(U)))= \H^2(U) \quad
\label{flat-meas}
\end{equation}
for every open set $U\subset 5 Q_0$. For the sake of convenience,
we assume in the whole proof
\begin{equation}
\eps_0 < \frac 1{100}\, . \label{small-lip}
\end{equation}
It is an easy computation to check that for every two points
$x,y\in 5Q_0$ we have
\begin{equation}
\label{n-grad} (1-2\eps_0) |\nabla f(x)-\nabla f(y)|\, \le \,
|n_\Sigma(x)-n_\Sigma(y)|\, \le\, (1+2\eps_0) |\nabla f(x)-\nabla
f(y)|\, .
\end{equation}
\heikodetail{

\bigskip

see back of my old p. 31 to see that

\bigskip

} We fix an orthonormal basis $(e_1,e_2,e_3)$ of $\rd$ so that
$e_1,e_2$ are parallel to the sides of $Q_0$.

\medskip\noindent\textbf{Step 2.}
Set, for $r\le \diam Q_0<1$, and any subset $S\subset Q_0$
\begin{eqnarray*}
\Phi_1^\ast(r, S) & :=  &
\max_{{\scriptstyle \|y-z\, \|\le r}\atop {\scriptstyle y,z \in Q_0\cap S}} |n_\Sigma(y)-n_\Sigma(z)|\, , \\
\Phi_2^\ast(r,S) & :=  & \max_{{\scriptstyle \|y-z\, \|\le r}\atop
{\scriptstyle y,z \in Q_0\cap S}} |\nabla f (y)-\nabla f (z)|\, , \\
\Phi^\ast(r,S) & := & \Phi_1^\ast(r,S) + \Phi_2^\ast(r,S)\, ,
\end{eqnarray*}
where $\|\cdot\|$ denotes the $\ell^\infty$ norm in $\bbbr^2$,
i.e. $\|x\|:= \max (|x_1|,|x_2|)$ for $x=(x_1,x_2)$. Shrinking
$Q_0$ if necessary, we may assume that
\begin{equation}
\Phi^\ast(\diam Q_0,Q_0) \le \frac 1{100}\,  \label{1-100}
\end{equation}
(by continuity of $n_\Sigma $ and of $\nabla f$.)

As in \cite[Section 6]{ssvdm-triple}, we want to prove the
following

\medskip\noindent\textbf{Key estimate.} Assume that $u,v\in Q_0$ and let $Q(u,v):={}$
the cube centered at $(u+v)/2$ and having edge length $2|u-v|$.
There exist positive numbers
$\delta_2=\delta_2(E,p)=a_2(p)E^{-1/(p-8)}$ and $C(p)>0$ such that
whenever $0<|u-v|\le \delta_2$, then
\begin{equation}
|n_\Sigma (u)-n_\Sigma (v)| \le 40\, \Phi^\ast
\Bigl(\frac{2|u-v|}{N}, Q(u,v)\Bigr) + C(p) E(u,v)^{1/p}
|u-v|^\lambda, \label{key}
\end{equation}
where $N$ is a (fixed) large natural number such that
$(N/2)^\kappa > 240 $ and
\[
E(u,v):= \int_{[F(Q(u,v))\cap \Sigma]^4} \K^p \, d\mu\, .
\]
One should view the second term on the right-hand side of
\eqref{key} as the main one. The first one is just an error term
that can be iterated away by scaling the distances down to zero.
\smallskip

We now postpone the proof of \eqref{key} for a second and show
that it yields the desired result upon iteration.

Note that \eqref{n-grad} and \eqref{small-lip} imply $\Phi^\ast
\le 3 \Phi_1^\ast$. \heikodetail{

\bigskip

$$
\Phi^\ast= \Phi_1^\ast+ \Phi_2^\ast\le
\Phi_1^\ast+\frac{1}{1-2\eps_0}\Phi_1^\ast\le
\left(1+\frac{1}{1-\frac{1}{50}}\right)\Phi_1^\ast=
\left(1+\frac{50}{50-1}\right)\Phi_1^\ast=\frac{99}{49}\Phi_1^\ast<3\Phi_1^\ast.
$$

\bigskip

}%
Moreover, if $u,v\in B(a,R)$ and $\|u-v\|=r\le R$, then
$Q(u,v)\subset B(\frac{u+v}2,\sqrt 2 \, |u-v|)\subset
B(\frac{u+v}2,2\|u-v\|)\subset B(a,R+2r)$.
\heikodetail{

\bigskip

$$\|u-v\|\ge\frac{1}{\sqrt{2}}|u-v|$$

is in the worst case (of a square) equality by Pythagoras

\bigskip

}
Thus, denoting
\[
M_p (a,\rho) := \biggl(\int_{[F(B(a,\rho))\cap \Sigma]^4} \K^p \,
d\mu\, \biggr)^{1/p}, \qquad a\in Q_0, \rho>0,
\]
and taking the supremum over $u,v\in B(a,R)$ with $|u-v|\le r\le R$,
one checks that \eqref{key} implies
\begin{eqnarray}
\label{pre-morrey} \Phi^\ast(r,B(a,R)) & \le & 120\,
\Phi^\ast(r/n,B(a,R+2r)) \\
& &{}+ 3C(p) M_p(a,R+2r) \, r^\lambda\, , \qquad n\equiv N/2\, .
\nonumber
\end{eqnarray}
\heikodetail{

\bigskip

In fact, \eqref{key} implies that
\begin{eqnarray*}
\Phi^\ast(r,B(a,R)) & \le & 3\Phi^*_1(r,B(a,R)) =
3 \max_{\|y-z\|\le r\atop y,z\in Q_0\cap B(a,R)}
|n_\Sigma(y)-n_\Sigma(z)|\\
& \overset{\eqref{key}}{\le} &
\max_{\|y-z\|\le r\atop y,z\in Q_0\cap B(a,R)}
\left[120\Phi^*\Big(\frac{2|y-z|}{N},Q(y,z)\Big)+
3C(p)E(y,z)^(1/p)|y-z|^\lambda\right],
\end{eqnarray*}
and now use monotonicity of $\Phi^*$ in both of its entries
together with the previously proven inclusions

\bigskip

}
A technique which is standard in PDE allows to get rid of the
first term on the right-hand side of this inequality. Indeed, upon
iteration \eqref{pre-morrey} implies
\begin{eqnarray*}
\Phi^\ast(r,B(a,R)) & \le & 120^j\,
\Phi^\ast(r/n^j,B(a,R+2\sigma_j)) \\
& & {}+ 3C(p)M_p(a,R+2\sigma_j)\,  r^\lambda\,  \sum_{i=0}^{j-1}
\left(\frac{120}{n^\lambda}\right)^i\, ,\qquad j=1,2, \ldots
\end{eqnarray*}
where
\[
\sigma_j:= r\sum_{i=0}^{j-1} n^{-i} \le 2r.
\]
\heikodetail{

\bigskip

\begin{eqnarray*}
\Phi^\ast(r,B(a,R)) & < & 120\Big[120\Phi^\ast(r/n^2,B(a,R+2r+\frac{2r}{n}))+3C(p)
M_p(a,R+2r+\frac{2r}{n})\Big]+3C(p)M_p(a,R+2r)r^\lambda\\
& \le &
120^2\Phi^\ast(\frac{r}{n^2},B(a,R+2r\sum_{i=0}^{2-1}n^{-i}))+
3C(p)
M_p(a,R+2r+2r\sum_{i=0}^{2-1}n^{-i})r^\lambda \Big[\frac{120}{n^\lambda}+1\Big]
\end{eqnarray*}

\bigskip

} As $n^\lambda = (\frac N2)^\lambda > ( \frac N2 )^\kappa > 240$,
we obtain $120/n^\lambda<1/2$ which implies
$\sum_i(120/n^\lambda)^i<2$ and hence
\[
\Phi^\ast(r,B(a,R)) < 120^j \Phi^\ast\bigl(r/n^j,B(a,R+4r)\bigr) +
6 C(p)M_p(a,R+4r)\, r^\lambda\, ,\qquad j=1,2, \ldots
\]
Now by Corollary \ref{cor:5.4new} we have a priori
$\Phi^\ast(r,S)\le \Phi^\ast(r,Q_0) \le C r^\kappa$ for every set
$S\subset Q_0$ and $r\le\diam Q_0$. Thus,
\[
120^j \Phi^\ast(r/n^j,B(a,R+4r)) \le C r^\kappa
\bigl(120/n^\kappa\bigr)^j < C r^\kappa 2^{-j}
\]
by choice of $N$. \heikodetail{

\bigskip

$$
\Phi^\ast(r/n^j,...))\le C(r/n^j)^\kappa <C/n^{\kappa j}=C2^{\kappa j}/N^{
\kappa j}
$$
since $r<1$.

\bigskip

} Passing to the limit $j\to \infty$ and setting $R=r$, we obtain
\begin{equation}
\label{6K} \Phi^\ast(r,B(a,r)) \le 6\, C(p)M_p(a,5r)  r^\lambda\,
,
\end{equation}
and this oscillation estimate immediately implies the desired
H\"older estimate \eqref{normal-hoelder-estimate} for the unit
normal  vector. In the remaining part of the proof, we just verify
\eqref{key}.

\medskip\noindent\textbf{Step 3: bad and good parameters.} From
now on, we assume that $u\not= v\in Q_0$ are fixed. We pick the
sub\-cube $Q=Q(u,v)$ of $5Q_0$ with edges parallel to those of
$Q_0$, so that the center of $Q(u,v)$ is at $(u+v)/2$ and the edge
of $Q(u,v)$ equals $2|u-v|$. Set
\begin{equation}
\label{m-CN} m = (20N)^{-2},  \qquad C_m=m^{-4},
\end{equation}
and
consider the sets of \emph{bad parameters\/} defined as follows:
\begin{eqnarray}
\Sigma_0 & = & \{ x_0 \in Q \, \colon \, \H^2(\Sigma_1(x_0)) \ge
m |u-v|^2\}, \label{sig0}\\
\Sigma_1(x_0) & = & \{ x_1 \in Q \, \colon \,
\H^2(\Sigma_2(x_0,x_1)) \ge
m |u-v|^2\}, \label{sig1}\\
\Sigma_2(x_0,x_1) & = & \{ x_2 \in Q \, \colon \,
\H^2(\Sigma_3(x_0,x_1,x_2)) \ge
m |u-v|^2\}, \label{sig2}\\
\Sigma_3(x_0,x_1,x_2) & = & \{ z \in Q \, \colon \,
\K(x_0,x_1,x_2,z) > \bigl(C_m
E(u,v)\bigr)^{1/p}|u-v|^{-8/p}\}.\label{sig3}
\end{eqnarray}
A word of informal explanation to motivate the above choices:
\emph{if we already knew\/} that $\Sigma$ is of class
$C^{1,\lambda}$, $\lambda=1-8/p$, then close to $u$ we would have
lots of tetrahedra with two perpendicular edges of the base having
length $\approx |u-v|$, and the height $\lesssim
|u-v|^{1+\lambda}$. For such tetrahedra our curvature integrand
does not exceed, roughly, a multiple of
$|u-v|^{\lambda-1}=|u-v|^{-8/p}$.
\heikodetail{

\bigskip

$V(T)\simeq d^{2+1+\lambda}=d^{3+\lambda}$ and $A(T)\simeq d^2$ and
$\diam(T)\simeq d^2$

\bigskip

}
Of course, there is no reason to
believe a priori that it is indeed the case. But it helps, as we
shall check, to look at tetrahedra that violate this naive
estimate, and to try and estimate how many of them there are.

We first estimate the measure of $\Sigma_0$. Using
\eqref{flat-meas} which gives a comparison of $d\H^2$ on
$\Sigma\cap F(5 Q_0)$ with the Lebesgue measure in $5Q_0$, we
obtain
\begin{eqnarray*}
E(u,v) & \ge & \int_{\Sigma_0} \int_{\Sigma_1(x_0)}
\int_{\Sigma_2(x_0,x_1)} \int_{\Sigma_3(x_0,x_1,x_2)}
\K^{p}(x_0,x_1,x_2,z)\, d\H^2_z\, d\H^2_{x_2}\, d\H^2_{x_1}\,
d\H^2_{x_0} \\
& > & C_m E(u,v)\,  m^3 {|u-v|^{-2}} \H^2(\Sigma_0) \\[3pt]
& = &  E(u,v)\,  m^{-1} {|u-v|^{-2}} \H^2(\Sigma_0),
\end{eqnarray*}
\heikodetail{

\bigskip

\begin{eqnarray*}
E(u,v) & \ge & \int_{\Sigma_0} \int_{\Sigma_1(x_0)}
\int_{\Sigma_2(x_0,x_1)} \int_{\Sigma_3(x_0,x_1,x_2)}
\K^{p}(x_0,x_1,x_2,z)\, d\H^2_z\, d\H^2_{x_2}\, d\H^2_{x_1}\,
d\H^2_{x_0} \\
& \overset{\eqref{sig3}}{>} & \int_{\Sigma_0} \int_{\Sigma_1(x_0)}
\int_{\Sigma_2(x_0,x_1)} C_mE|u-v|^{-8}\H^2(\Sigma_3(x_0,x_1,x_2))\,d\H_{x_2}
d\H^2_{x_1}d\H^2_{x_0}\qquad \\
& \overset{\eqref{sig2}}{\ge} & \int_{\Sigma_0}
\int_{\Sigma_1(x_0)}\int_{\Sigma_2(x_0,x_1)}
mC_mE|u-v|^{-6}d\H^2_{x_2}d\H^2_{x_1}d\H^2_{x_0}\\
&= &
\int_{\Sigma_0}\int_{\Sigma_1(x_0)}
\H^2(\Sigma_2(x_0,x_1))mC_mE|u-v|^{-6}d\H^2_{x_1}d\H^2_{x_0}\\
& \overset{\eqref{sig1}}{\ge} & \int_{\Sigma_0} \int_{\Sigma_1(x_0,x_1)}
m^2C_mE|u-v|^{-4}d\H^2_{x_1}d\H^2_{x_0}
=
\int_{\Sigma_0} \H^2(\Sigma_1(x_0)m^2C_mE|u-v|^{-4}\,d\H^2_{x_0}\\
& \overset{\eqref{sig0}}{\ge} & m^2C_mE|u-v|^{-2}\H^2(\Sigma_0)
\end{eqnarray*}

\bigskip

}%
which yields
\begin{equation}
\label{small-bad} \H^2(\Sigma_0) < m |u-v|^2 =
\frac{|u-v|^2}{400N^2} \ll |Q(u,v)|=4|u-v|^2\, .
\end{equation}

\medskip\noindent\textbf{Step 4: auxiliary good points.} In a small
neighbour\-hood of  $u$ we select $x_0 \in Q(u,v) \setminus
\Sigma_0$ so that $\|x_0-u\|\le (20N)^{-1}|u-v|$.
\heikodetail{

\bigskip

If there were no such ``good'' point $x_0$ then the whole square
$Q^\ast$ centered at $u$ and edges of length $2(20N)^{-1} \|u-v\|$
parallel to the ones of $Q$ would only contain points of the
``bad'' set $\Sigma_0$, but then $\Sigma_0$ would have too large
measure:
$$
\H^2(\Sigma_0)\ge
\H^2(Q^\ast)=(2/20N)^2|u-v|^2\overset{\eqref{small-bad}}{>}\H^2(\Sigma_0).
$$
(see also the small images on old p.33)

\bigskip

} Once $x_0$ is chosen, we select $x_1\in Q(u,v)\setminus
\Sigma_1(x_0)$ and then $x_2\in Q(u,v)\setminus \Sigma_2(x_0,x_1)$
so that
\[
\|x_1-x_0\|\approx \|x_2-x_0\|\approx \frac{|u-v|}{N}
\qquad\mbox{and}\qquad \ang (x_2-x_0,x_1-x_0) \approx \frac \pi
2\, .
\]
More precisely, let $Q(x_0)$ be the cube with one vertex at $x_0$
and two other vertices at
\[
a_1:= x_0 + \frac{|u-v|}N e_1, \qquad a_2 :=  x_0 + \frac{|u-v|}N
e_2\, .
\]
We select $x_1,x_2\in Q(x_0)$ such that
\begin{gather}
\label{x1} x_1 \in Q(x_0) \setminus \Sigma_1(x_0), \qquad
\|x_1-a_1\|\le \frac{|u-v|}{20N}, \\
\label{x2} x_2 \in Q(x_0) \setminus \Sigma_2(x_0,x_1), \qquad
\|x_2-a_2\|\le \frac{|u-v|}{20N}\, .
\end{gather}
(See also the figure below.) Since $x_0\not \in \Sigma_0$, we can
use \eqref{sig0}--\eqref{sig1} to check that $x_1,x_2$ satisfying
\eqref{x1}--\eqref{x2} do exist. \heikodetail{

\bigskip

$$
\H^2(\Sigma_1(x_0))<m|u-v|^2
$$
and the cube $\{\xi\in\R^2:\|\xi-a_1\|\le |u-v|/20N\}$ has
$\H^2$-measure $ 2^2|u-v|^2/(400N)^2>m|u-v|^2$, which implies the
existence of $x_1$ with \eqref{x1}. Similarly,
$$
\H^2(\Sigma_2(x_0,x_1))<m|u-v|^2
$$
implies the existence of a point $x_2\in
Q(x_0)\setminus\Sigma_2(x_0,x_1)$ with
$\|a_2-x_2\|\le\frac{|u-v|}{20N},$ since the cube
$\{\xi\in\R^2:\|\xi-a_2\|\le |u-v|/20N\}$ has $\H^2$-measure $
2^2|u-v|^2/(400N)^2$. ($x_2$ could be contained in $\Sigma_1(x_0)$
though!)

\bigskip

}

\begin{center}
\hfill
\includegraphics*[totalheight=7.5cm]{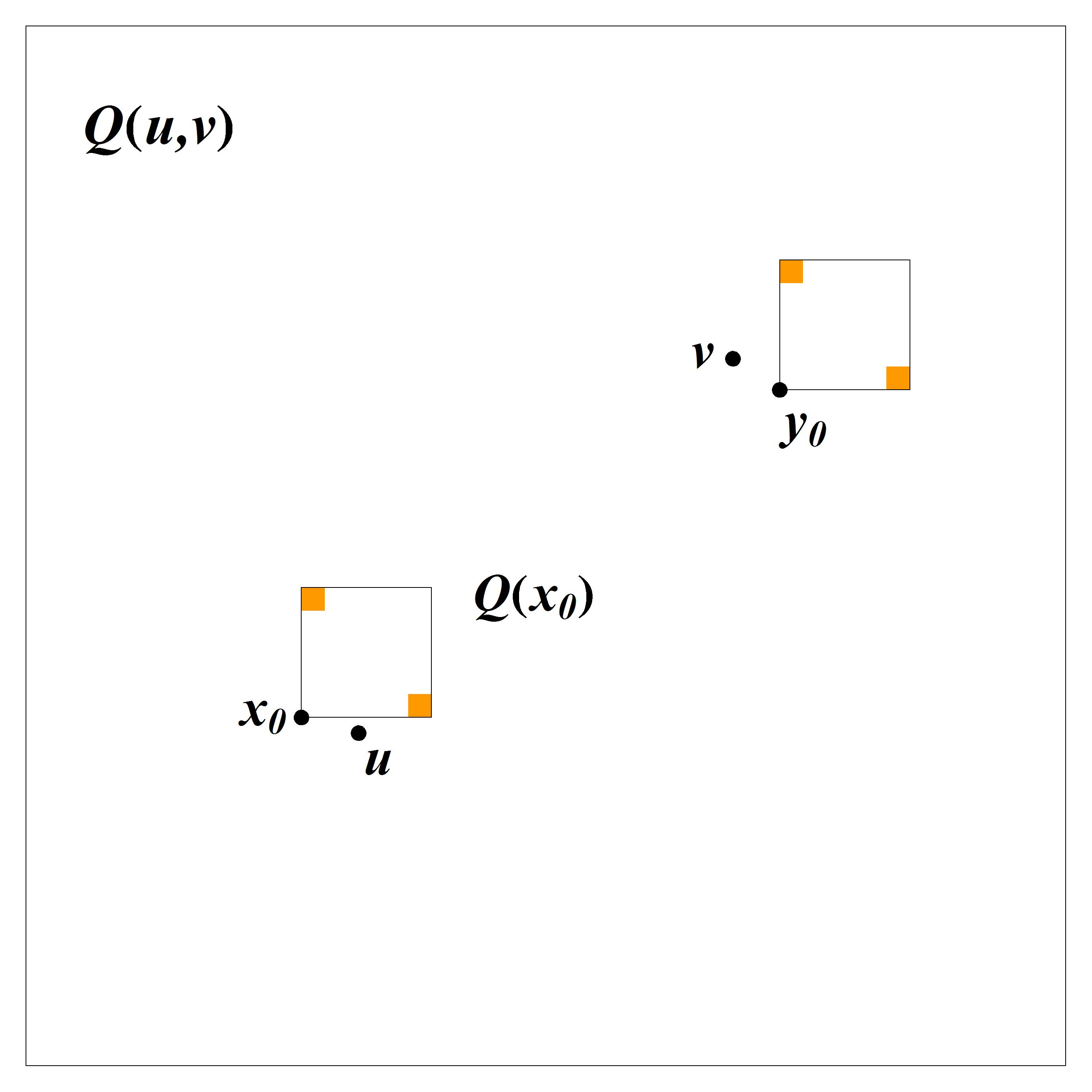}
\hfill
\includegraphics*[totalheight=7.5cm]{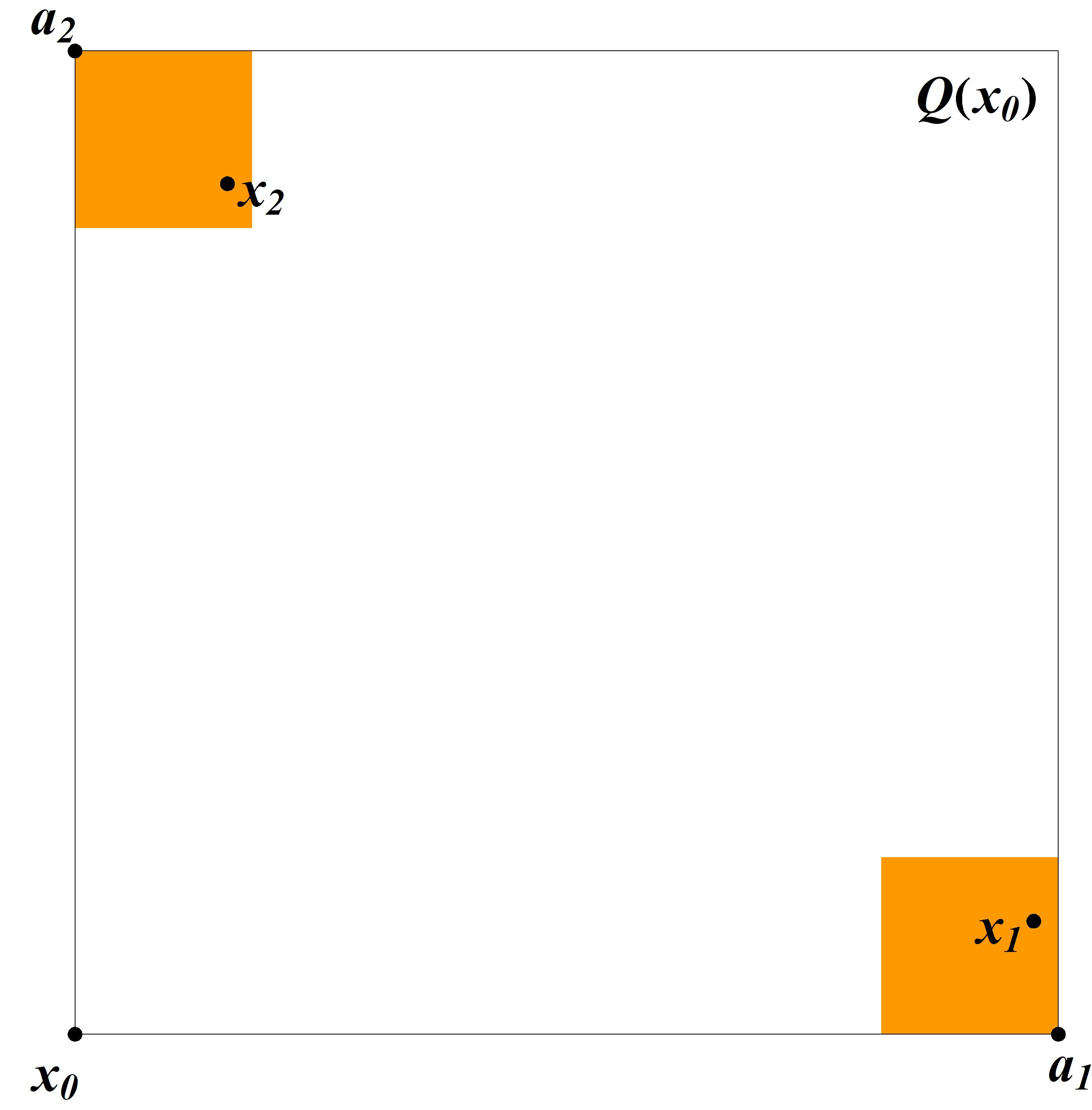}
\hfill
\end{center}

{\footnotesize

\medskip\noindent\textbf{Fig. 7.} The position of auxiliary good
parameters in the domain of $f$. Left: $Q(u,v)$ and two sub\-cubes
$Q(x_0),Q(y_0)$, with lower left-hand corners at $x_0,y_0$. Right:
$Q(x_0)$ magnified. We fix $x_0\not\in \Sigma_0$, close to $u$,
and $x_1,x_2$ are selected in the little shaded sub\-cubes of
$Q(x_0)$. Since the Lip\-schitz constant of $f$ is small, $\Sigma$
is a flat graph over $Q(u,v)$. Thus, the vectors
$v_j:=F(x_j)-F(x_0)$ ($j=1,2$) are nearly orthogonal and have
lengths very close to $|u-v|/N={}$ the edge of $Q(x_0)$, see
Step~5 below for the details.

\medskip

}

In a fully analogous way we select $y_0,y_1,y_2$ close to $v$ ---
using \eqref{small-bad} initially again but then by defining
 sets $\Sigma_1(y_0)$,
$\Sigma_2(y_0,y_1),$ and $\Sigma_3(y_0,y_1,y_2)$ as in
\eqref{sig1}, \eqref{sig2}, and \eqref{sig3}. Thus, $y_0\in
Q(u,v)\setminus \Sigma_0$, $y_1\in Q(y_0)\setminus \Sigma_1(y_0)$
and $y_2\in Q(y_0) \setminus \Sigma_2(y_0,y_1)$, where $Q(y_0)$ is
a copy of $Q(x_0)$ translated by $y_0-x_0$, satisfy
\[
\|y_0-v\|\le  \frac{|u-v|}{20N}\, , \qquad \|y_1-y_0\|\approx
\|y_2-y_0\|\approx \frac{|u-v|}{N}\, , \qquad \ang
(y_2-y_0,y_1-y_0) \approx \frac \pi 2\, .
\]
Then we set $P_x:=\langle F(x_0),F(x_1),F(x_2)\rangle$,
$P_y:=\langle F(y_0), F(y_1), F(y_2) \rangle$, and we let $n_x,
n_y$ denote the unit normal vectors of these two planes. By the
triangle inequality,
\begin{eqnarray*}
|n_\Sigma(u)-n_\Sigma(v)| & \le & |n_\Sigma(u)-n_\Sigma(x_0)| + |n_\Sigma(x_0)-n_x| \\
& & {}+ |n_x-n_y| \\
& & {}+ |n_y-n_\Sigma(y_0)|+|n_\Sigma(y_0)-n_\Sigma(v)|.
\end{eqnarray*}
The non-obvious term is the middle one, $|n_x-n_y|\le \ang
(P_x,P_y)$; \heikodetail{

\bigskip

Indeed, since $n_x$ and $n_y$ are unit vectors one has for the
difference
$$
|n_x-n_y|=\sqrt{2-2\cos\ang(P_x,P_y)}=2\sqrt{\frac 12
(1-\cos\ang(P_x,P_y))} =2\sin(\ang(P_x,P_y)/2)\le \ang(P_x,P_y).
$$

\bigskip

}%
the remaining four terms give a small contribution which does not
exceed a constant multiple of $ \Phi^\ast (20\,|u-v|/N,Q(u,v))$.
But due to the choice of $\Sigma_3$ the planes $P_x$ and $P_y$
turn out to be almost parallel: their angle is $\lesssim
|u-v|^\lambda$.

Since $u,v$ are now fixed and will not change till the end of the
proof, from now we use the abbreviations
\[
\Phi_i^\ast(r)\equiv \Phi_i^\ast(r,Q(u,v)), \qquad
\Phi^\ast(r)\equiv \Phi^\ast(r,Q(u,v))\, .
\]
We shall check that
\begin{eqnarray}
|n_\Sigma(x_0)-n_x| & \le & 16 \Phi^\ast(2|u-v|/N)\, , \label{x-plane}\\
|n_\Sigma(y_0)-n_y| & \le & 16 \Phi^\ast(2|u-v|/N)\, ,
\label{y-plane} \\
|n_x-n_y| & \le & K |u-v|^\lambda\, . \label{planes}
\end{eqnarray}
Combining these estimates with the obvious ones,
\[
|n_\Sigma(u)-n_\Sigma(x_0)|\le \Phi^\ast (|u-v|/N)\, , \qquad
|n_\Sigma(v)-n_\Sigma(y_0)|\le \Phi^\ast (|u-v|/N)\, ,
\]
\heikodetail{

\bigskip

$$
|n_\Sigma(u)-n_\Sigma(x_0)|\le \Phi^\ast (\|u-x_0\|)\le \Phi^\ast
((20N)^{-1}|u-v|) \le \Phi^\ast ((N)^{-1}|u-v|)
$$
by monotonicity of $\Phi^\ast$, so it is a  rather rough estimate

\bigskip

} and using monotonicity of $\Phi^\ast$, one immediately obtains
\eqref{key}. \heikodetail{

\bigskip

\begin{eqnarray*}
|n_\Sigma(u)-n_\Sigma(v)|&\le & |n_\Sigma(u)-n_\Sigma(x_0)|+
|n_\Sigma(x_0)-n_x|+ |n_x-n_y|+ |n_y-n_\Sigma(y_0)|+
|n_\Sigma(y_0)-n_\Sigma(v)|\\
& \le  & \Phi^\ast(|u-v|/N)+16
\Phi^\ast(2|u-v|/N)+K|u-v|^\lambda+16 \Phi^\ast(2|u-v|/N)+
\Phi^\ast(|u-v|/N)\\
& \le &  34 \Phi^\ast(2|u-v|/N) \qquad\textnormal{by monotonicity
of $\Phi^\ast$}.
\end{eqnarray*}

\bigskip

}

\bigskip\noindent\textbf{Step 5: proofs of
\eqref{x-plane} and \eqref{y-plane}.} We only prove
\eqref{x-plane}; the other proof is identical. Let
\[
v_j:= F(x_j)-F(x_0), \qquad j=1,2.
\]
By the fundamental theorem of calculus,
\begin{eqnarray}
v_j & = & \int_0^1 \nabla F(x_0+t(x_j-x_0)) (x_j-x_0)\, dt \nonumber \\
& = & \nabla F(x_0) (x_j-x_0) + \int_0^1(\nabla F(x_0+t(x_j-x_0))
-\nabla F(x_0)) (x_j-x_0)\,dt\nonumber\\
& =: & w_j+\sigma_j,\quad\Fo j=1,2, \label{arms}
\end{eqnarray}
where the error terms $\sigma_j$ satisfy
\begin{eqnarray}
|\sigma_j| & \le & \left|\int_0^1 (\nabla F(x_0+t(x_j-x_0))-\nabla
F(x_0))\, dt (x_j-x_0|)\right| \nonumber \\
& \le & \Phi^\ast (2|u-v|/N) \, \diam Q(x_0)\nonumber \\
& \le & 2 \Phi^\ast (2|u-v|/N) \, \frac{|u-v|}N, \qquad j=1,2.
\label{errors}
\end{eqnarray}
\heikodetail{

\bigskip

Here, even $\Phi_1^\ast$ would have been sufficient...

\bigskip

More details for the above calculation: First observe that
\begin{eqnarray*}
|(\nabla F(x_0+t(x_j-x_0))- \nabla F(x_0) (x_j-x_0)| & = &
\left|\left[\left(\begin{array}{c}
\Id_{\R^2}\\
\nabla f(x_0+t(x_j-x_0))\end{array}\right)- \left(\begin{array}{c}
\Id_{\R^2}\\
\nabla f(x_0)\end{array}\right)
\right](x_j-x_0)\right|\\
& = & \left|\left( \begin{array}{c}
x_j-x_0-(x_j-x_0) \\
\big[ \nabla f(x_0+t(x_j-x_0))-\nabla f(x_0)\big] \cdot (x_j-x_0)
\end{array}\right)   \right| \\
& = & \left| \left( \begin{array}{c}
0 \\
\big[\nabla f(x_0+t(x_j-x_0))-\nabla f(x_0)\big]\cdot (x_j-x_0)
\end{array} \right) \right| \\
& \le & \Phi_2^\ast (\|x_j-x_0\|)|x_j-x_0|
\end{eqnarray*}
and to estimate the argument in the $\Phi^\ast$-function
\begin{eqnarray*}
\|x_1-x_0\| & = & \left\|
x_1-(x_0+(|u-v|/N)e_1)+(|u-v|/N)e_1\right\|
\\
& \le &
\|x_1-a_1\|+\left\|\frac{|u-v|}{N}e_1\right\|\\
& \overset{\eqref{x1},\eqref{x2}}{\le} & \frac{|u-v|}{20N}+
 \frac{|u-v|}{N}<\frac{2}{N}|u-v|.
 \end{eqnarray*}

 \bigskip

 }
With $w_j=\nabla F(x_0)\cdot (x_j-x_0)$, $j=1,2\,$ we have
\begin{equation}
\label{nornor} n_x=\frac{v_1\times v_2}{|v_1\times v_2|}\, ,
\qquad n_\Sigma(x_0)= \frac{w_1\times w_2}{|w_1\times w_2|}\, .
\end{equation}
To estimate the difference of these two vectors, we first estimate
$|v_j|,|w_j|$ and the angles $\ang(v_1,v_2)$, $\ang(w_1,w_2)$.
This is an elementary computation; we give some details
below.\footnote{If you do not want to check the details of our
arithmetic, please note the following: we use $N$ \emph{only\/} to
fix the scale and to control the ratio of $\diam Q(x_0)$ and
$\diam Q(u,v)$. Thus, $N$ \emph{does not\/} influence the ratio of
lengths of $v_1,v_2,w_1,w_2$ (which are all ${}\approx |u-v|/N$)
and the angles between these vectors (which are absolute since we
assume \eqref{lip-f} and \eqref{1-100}).

Therefore, the constant `$16$' in \eqref{x-plane}--\eqref{y-plane}
is not really important. Any absolute constant would be fine; one
would just have to adjust $N$ to derive \eqref{6K} from
\eqref{key}.}

Using the fact that $|\nabla f|$ is bounded by $\eps_0<1/100$ by
Remark \ref{rem:5.4newA} and \eqref{small-lip},
$x_1-x_0$ and $x_2-x_0$ are close to two
perpendicular sides of $Q(x_0)$, and both error terms $\sigma_j$
are smaller than $|u-v|/50N$ by \eqref{1-100}, \heikodetail{

\bigskip

Since $2|u-v|/N<2\diam Q_0$ we have by monotonicity of $\Phi^\ast$
$$
2\Phi^\ast(2|u-v|/N)\le 2\Phi^\ast(2\diam
Q_0)\overset{\eqref{1-100}}{\le } 2/100.
$$

\bigskip

} one can check that
\[
\frac{9}{10} \frac{|u-v|}N  \le  \min(|v_j|, |w_j|) \\
\le  \max(|v_j|, |w_j|) \le \frac{11}{10} \frac{|u-v|}N
\]
for $j=1,2$. \heikodetail{

\bigskip

For example
\begin{eqnarray*}
|v_j| & \ge & |\nabla F(x_0) (x_j-x_0)|-|\sigma_j|\\
& \ge & \left| \left(\begin{array}{c}
x_j-x_0\\
\nabla f(x_0)\cdot (x_j-x_0)
\end{array}\right)\right|-\frac{|u-v|}{50N}\\
& = &
\sqrt{|x_j-x_0|^2+|\nabla f(x)\cdot (x_j-x_0)|^2}-\frac{|u-v|}{50N}\\
& \ge & |x_j-x_0|-\frac{|u-v|}{50N}\\
& \ge & \|x_j-x_0\|-\frac{|u-v|}{50N}\\
& \ge & \|a_j-x_0\|-\|x_j-aj\|-\frac{|u-v|}{50N}\\
& \overset{\eqref{x1},\eqref{x2}}{\ge} &
\frac{|u-v|}{N}\left(1-\frac{1}{20} -\frac{1}{50} \right)\\
& = & \frac{93}{100} \cdot \frac{|u-v|}{N} > \frac{9}{10} \cdot
\frac{|u-v|}{N}.
\end{eqnarray*}

\bigskip

}

Note also that, cf. Figure 7~and \eqref{lip-f},
\[
v_j = \frac{|u-v|}N e_j + \sum_{i=1}^{3} a_{ji}e_i, \qquad
|a_{ji}|\le \frac{|u-v|\sqrt{2}}{20N},
\]
\heikodetail{

\bigskip

\begin{eqnarray*}
v_j& = & F(x_j)-F(x_0)=\left(\begin{array}{c}
x_j-x_0\\
f(x_j)-f(x_0)\end{array}\right)= \left(\begin{array}{c}
a_j-x_0\\
0\end{array}\right)+ \left(\begin{array}{c}
x_j-a_j\\
f(x_j)-f(x_0)\end{array}\right)\\
& = & \left(\begin{array}{c}
\frac{|u-v|}{N}e_j\\
0\end{array}\right)+ \left(\begin{array}{c}
x_j-a_j\\
f(x_j)-f(x_0)\end{array}\right),
\end{eqnarray*}
and one uses now that
$|x_j-a_j|\le\sqrt{2}\|x_j-a_j\|\le\sqrt{2}\frac{|u-v|}{20N}$ and
\begin{eqnarray*}
|f(x_j)-f(x_0)|&\le
&\eps_0|x_j-x_0|\le\sqrt{2}\eps_0\|x_j-x_0\|\le
\sqrt{2}\eps_0\Big[ \|x_j-a_j\|+\|a_j-x_0\|\Big]\\
& \overset{\eqref{x1},\eqref{x2}}{\le} &
\sqrt{2}\eps_0\Big[\frac{1}{N}+
\frac{1}{20N}\Big]|u-v|\overset{\eqref{small-lip}}{<}\frac{\sqrt{2}}{20N}|u-v|
\end{eqnarray*}

\bigskip

} which yields
\[
|v_1\cdot v_2 | = \left|\frac{|u-v|}N(a_{12}+a_{21}) +
\sum_{i=1}^3 a_{1i}a_{2i} \right| \le \frac{|u-v|^2}{6N^2}\, .
\]
\heikodetail{

\bigskip

this last inequality may be verified as follows:
$$
2\frac{\sqrt{2}}{20}+3\frac{2}{400}=\frac{\sqrt{2}}{10}+\frac{3}{200}
<
\frac{\sqrt{2}}{10}+\frac{3}{180}=\frac{\sqrt{2}}{10}+\frac{1}{60}
=\frac{6\sqrt{2}+1}{60}<\frac{6\cdot 1.5+1}{60}=\frac 16.
$$

\bigskip

}%
Taking the estimates of $\sigma_j$ into account one more time, we
obtain $| w_1\cdot w_2 |  \le {2|u-v|^2}/({9N^2})\,$.
\heikodetail{

\bigskip

$$
|w_1\cdot w_2|=|(v_1-\sigma_1)\cdot
(v_2-\sigma_2)|\le\frac{|u-v|^2}{N^2} \left[
\frac{1}{6}+\frac{2}{50}\cdot\frac{11}{10}+\frac{1}{2500}\right]
\le  0.2111<\frac 29.
$$

\bigskip

} Combining the inequalities for these two scalar products with
the estimates of lengths of the vectors, we conclude that
\[
\max \Bigl( |\cos\ang(v_1,v_2)|, |\cos\ang(w_1,w_2)|\Bigr) \le
\frac 29 \cdot\left(\frac{10}9\right)^2 \, .
\]
Hence,
\begin{equation}
\min  \Bigl( \sin \ang(v_1,v_2), \sin\ang(w_1,w_2)\Bigr) \ge
\sqrt{1-\left[\frac 4{81}\left(\frac{10}9\right)^4\right]}>
 \frac{15}{16}\, . \label{sines}
\end{equation}
\heikodetail{

\bigskip

no fraction manipulation done here, just the desk calculator...

\bigskip

} Now,
\begin{equation}\label{A}
A:= v_1\times v_2 - w_1\times w_2 =  |v_1\times v_2|n_x -
|w_1\times w_2| n_\Sigma(x_0)\, .
\end{equation}
As $v_j=w_j+\sigma_j$ and $|w_j|\le 11|u-v|/(10N)$, we have
\begin{eqnarray*}
|A|\le  |\sigma_1|\, |w_2| + |\sigma_2|\, |w_1| + |\sigma_1|\,
|\sigma_2| & \overset{\eqref{errors}}{\le}& \left[2\cdot
\frac{11}{10}+\frac{1}{50}
\right] 2\Phi^\ast (2|u-v|/N) \, \frac{|u-v|^2}{N^2}\\
& < &6 \Phi^\ast (2|u-v|/N) \, \frac{|u-v|^2}{N^2}
\end{eqnarray*}
by \eqref{errors} and \eqref{1-100}. On the other hand, applying
the triangle inequality, using \eqref{sines}, and the estimates
$|v_j|\ge 9|u-v|/(10N)$ for $j=1,2$, we obtain first
\begin{equation}
|v_1\times v_2| \overset{\eqref{sines}}{>}
\left(\frac{9}{10}\right)^2\frac{15}{16}\frac{|u-v|^2}{N^2} >
\frac 34 \frac{|u-v|^2}{N^2},\label{cross}
\end{equation}
and then, using the second identity for $A$ in \eqref{A},
\begin{eqnarray*}
|A| & = & \bigl||v_1\times v_2| (n_x-n_\Sigma(x_0)) +
n_\Sigma(x_0)
(|v_1\times v_2|-|w_1\times w_2|)\bigr| \\
& \ge & |v_1\times v_2|\, |n_x-n_\Sigma(x_0)|
- |v_1\times v_2 - w_1\times w_2| \\
& \ge & \frac{3}{4} \frac{|u-v|^2}{N^2} |n_x-n_\Sigma(x_0)| -
|A|\, .
\end{eqnarray*}
Combining the lower and the upper estimate for  $A$  we obtain
\[
|n_x-n_\Sigma(x_0)| \le \frac 83 |A|
\left(\frac{|u-v|^2}{N^2}\right)^{-1} \le 16 \Phi^\ast
(2|u-v|/N)\, ,
\]
which yields \eqref{x-plane}.

\bigskip\noindent\textbf{Step 6:
proof of \eqref{planes}.} If $P_x$ is parallel to $P_y$, there is
nothing to prove. Let us then  suppose that these planes intersect
and denote their angle by $\gamma_0$. To show that
$\gamma_0\lesssim |u-v|^\lambda$, we use again the definition of
bad sets. Note that for
\begin{equation}
G\colon= Q(u,v)\setminus \Bigl( \Sigma_3(x_0,x_1,x_2)\,
\cup\,\Sigma_3(y_0,y_1,y_2)\Bigr)
\end{equation}
we have by \eqref{sig2}
\begin{equation}
\H^2(G) > |Q(u,v)|-2m |u-v|^2 = (2|u-v|)^2-2m|u-v|^2> |u-v|^2
\label{meas-G}
\end{equation}
by
choice of $m$. Therefore, as $\lambda-1=-8/p$, for all $z\in G$ we
have according to \eqref{sig3} the two inequalities
\begin{equation}
\label{1R-on-S3} \K(x_0,x_1,x_2,z)\le  K_0 |u-v|^{\lambda-1}\, ,
\qquad \K(y_0,y_1,y_2,z) \le K_0 |u-v|^{\lambda-1}\, ,
\end{equation}
where
\begin{eqnarray*}
K_0=K_0(p,E(u,v)) & := & (20N)^{8/p}  E(u,v)^{1/p}
\\
& \equiv & C_4(p) E(u,v)^{1/p},
\end{eqnarray*}
as we have in fact chosen $N$ depending only on
$\kappa=(p-8)/(p+16)$.

We are now going to use formula \eqref{1/R} for $\K$ to estimate
the distance from $F(z)$ to the planes $P_x$ and $P_y$. Setting
$v_j:=F(x_j)-F(x_0)$ for $j=1,2$ (as in the previous step of the
proof), and $v_3:=F(z)-F(x_0)$,  we obtain for the tetrahedron
$T:=(F(x_0),F(x_1),F(x_2),F(z))$
\[
|v_3|\le (1+\eps_0)|z-x_0|<2|u-v|, \qquad \diam T<2|u-v|
\]
by virtue of \eqref{small-lip}. Since the $|v_j|$ for $j=1,2$ have
been estimated before, this yields an estimate of the area of $T$,
\begin{eqnarray}
2 A(T) & = & |v_1\times v_2|+|v_2\times v_3|+|v_1\times
v_3|+|(v_2-v_1)\times (v_3 - v_2)| \nonumber \\
& \le & \left(\frac{11}{10}\right)^2\frac{|u-v|^2}{N^2} + 4\,
\left(\frac{11}{10}\, \frac{|u-v|}{N}\right) \diam T \nonumber \\
& \le & \frac{15|u-v|^2}{N} \qquad\mbox{as $N>1$.} \label{A-of-T}
\end{eqnarray}
Thus
\begin{eqnarray}
\K(x_0,x_1,x_2,z) & = & \frac{\dist (F(z),P_x)}{3(\diam
T)^2}\cdot\frac{|v_1\times v_2|}{ 2A(T)}\nonumber\\
& \overset{\eqref{cross}}{\ge} & \frac{\dist (F(z),P_x)}{16N^2} \,
\bigl( 2
A(T)\bigr)^{-1} \nonumber \\
& \ge & \frac{\dist (F(z),P_x)}{N^2 |u-v|^2} \, . \label{est-1R}
\end{eqnarray}
For the last inequality we have simply used \eqref{A-of-T} and the
inequality $N > (N/2)^\kappa>240$ which follows from our initial
choice of $N$.
Since the points $y_0,y_1,y_2$ have been chosen analogously to
$x_0,x_1,x_2$, it is clear that we also have
\begin{equation}
\label{est-1R-y} \K(y_0,y_1,y_2,z) > \frac{\dist
(F(z),P_y)}{N^2|u-v|^2}\, .
\end{equation}
Combining \eqref{1R-on-S3}--\eqref{est-1R-y}, we obtain
\begin{eqnarray}
\label{z-PxPy} \max\bigl(\dist(F(z),P_x),\dist(F(z),P_y)\bigr) &
< & N^2 K_0 |u-v|^{1+\lambda} \\
& \equiv &  C_5(p) E(u,v)^{1/p} |u-v|^{1+\lambda}, \qquad z\in G.
\nonumber
\end{eqnarray}
We shall show that the combination of \eqref{meas-G} and
\eqref{z-PxPy} implies that  $|n_x-n_y|\le\gamma_0=\ang(P_x,P_y)$
is estimated by a constant multiple of $|u-v|^\lambda$ thus
establishing \eqref{planes} as the only missing ingredient for the
proof of the key estimate \eqref{key}.

\medskip

Indeed, consider an affine plane $P$ which is perpendicular both
to $P_x$ and $P_y$. Let $\pi_P$ denote the orthogonal projection
onto $P$. By \eqref{z-PxPy} above, we see that $\pi_P(F(G))$ is a
subset of a rhombus $R$ contained in the plane $P$. The height of
this rhombus is equal to
\[
h=2 \cdot C_5(p) E(u,v)^{1/p} |u-v|^{1+\lambda}
\]
and the (acute) angle of $R$ is $\gamma_0$, so that the longer
diagonal of $R$ equals
\[
D=\frac{h}{\sin(\gamma_0 /2)} = \frac{2 C_5(p) E(u,v)^{1/p}
|u-v|^{1+\lambda}}{\sin (\gamma_0/2)}
\]
\heikodetail{

\bigskip

Here, a picture as on my old p. 37 should help!

\bigskip

} Therefore, the set $F(G)$ is contained in a cylinder $C_0$ with
axis $l:=P_x\cap P_y$ and radius $D/2$,
\begin{equation}
F(G) \subset C_0\colon= \{ w\, \colon \, \dist (w, l)\le D/2\}\, .
\end{equation}
The orthogonal projection of $C_0$ onto the plane containing the
domain of $f$ (recall that $F(x)=(x,f(x))$ para\-metrizes a
portion of $\Sigma$ that we consider) gives us a strip $S$ of
width $D$. This strip must contain all good parameters $z\in G$,
so that, taking \eqref{meas-G} into account, we have
\begin{eqnarray*}
3D |u-v| >2\sqrt{2} D |u-v| & = & D \, \diam Q(u,v) \\
& > & \mbox{ area of $S\cap Q(u,v)$ } \ge \H^2(G)
 > |u-v|^2\, .
\end{eqnarray*}
Hence, $D> |u-v|/3$, so that
\[
\frac{2}{\pi}\frac{\gamma_0}{2}\le\sin \frac{\gamma_0}2
=\frac{2C_5(p)E(u,v)^{1/p}|u-v|^{1+\lambda}}{D}< 6C_5(p)
E(u,v)^{1,p} |u-v|^\lambda,
\]
and hence
\[
|n_x-n_y|\le\gamma_0 <   6\pi C_5(p) E(u,v)^{1,p} |u-v|^\lambda
\]
which establishes \eqref{planes} and therefore concludes the whole
proof. Note that we have obtained the key estimate \eqref{key}
with $C(p)=6\pi C_5(p)$ depending only on $p$, as desired. \hfill
$\Box$

\smallskip

Applying the above result, one can sharpen
Corollary~\ref{cor:5.4new} as follows.

\begin{corollary}\label{cor:6.2}
Assume that $p>8$, $\M_p(\Sigma)< E < \infty.$ Then there exist
two constants, $0<\tilde{a}(p)<1<\tilde{A}(p)<\infty$, such that
for each $x\in\Sigma$ there is a function
$$
f:T_x\Sigma\to \Big( T_x\Sigma\Big)^\perp \simeq\R
$$
with the following properties:
\begin{enumerate}
\item[\rm (i)]
$f(0)=0,$ $\nabla f(0)=(0,0)$,
\item[\rm (ii)] For $\tilde{R}_1\equiv \tilde{R}_1(E,p):=\tilde
a(p)E^{-1/(p-8)}$ we have the estimate
\[
|\nabla f(y_1)-\nabla f(y_2)|\le \tilde A(p)\M_p\bigl(\Sigma\cap
B(x,{10\tilde{R}_1})\bigr)^{\frac{1}{p}}|y_1-y_2|^{1-8/p}, \qquad
y_1,y_2\in B(x,{\tilde{R}_1})
\]
\item[\rm (iii)]
\heikodetail{

\bigskip

The radius factor 10 instead of 5 comes from the standard transition
from oscillation estimates to H\"older estimates

\bigskip

}
The map
$$
\Phi(y):=x+(y,f(y)),\quad y\in T_x\Sigma\simeq\R^2,
$$
satisfies
\begin{equation}\label{graph-patch2}
\Phi(D_{\frac{3}{4}\tilde{R}_1})\subset
B(x,{\tilde{R}_1})\cap\Sigma\subset\Phi(D_{\tilde{R}_1}),
\end{equation}
where $D_{\tilde{R}_1}=B(0,\tilde{R}_1)\cap T_x\Sigma$ is a disk
in $T_x\Sigma$ around $0\in T_x\Sigma,$ and
\begin{equation}\label{graph-est2}
|D\Phi(y_1)-D\Phi(y_2)|\le A(p) \M_p\bigl(\Sigma\cap
B(x,{10\tilde{R}_1})\bigr) |y_1-y_2|^{1-8/p}, \qquad y_1,y_2\in
B(x,{\tilde{R}_1}).
\end{equation}
\end{enumerate}
\end{corollary}
Of course, in (ii) and (iii) one can replace $\M_p\bigl(\Sigma
\cap\ldots{}\bigr)$ by the total energy of the surface thus
providing clear-cut a priori estimates to be used in the next
section.


\section{Sequences of equibounded $\M_p$-energy}\label{sec:7}

\setnumbers

The main issue of this final Section is the proof of the following
compactness theorem for admissible surfaces of equibounded energy
with a uniform area bound. Notice that such an additional area
bound is necessary as the example of larger and larger spheres
shows. Let $S_\rho:=\partial B(0,\rho)$. For any tetrahedron $T$
(with non-co\-planar vertices) we estimate
\begin{equation}
\label{eq:h7.1} \K(T)\ge \frac{1}{6R(T)},
\end{equation}
where $R(T)$ denotes the radius of the circumsphere of
$T=(x_0,x_1,x_2,x_3)$. There is an explicit formula,
\[
\frac{1}{2R(T)} = \frac{\bigl|\langle z_3, z_1\times z_2
\rangle\bigr|}{\bigl| \, |z_1|^2 z_2 \times z_3 + |z_2|^2 z_3
\times z_1 + |z_3|^2 z_1 \times z_2\, \bigr| }\, ,
\]
where we have set $z_i=x_i-x_0$ for $i=1,2,3$; this formula can be
compared to \eqref{1/R} in order to obtain \eqref{eq:h7.1}. Hence,
\[
\M_p(S_\rho) \gtrsim \rho^{8-p} \to 0 \qquad\mbox{as $\rho\to
\infty$.}
\]

\begin{theorem}\label{Thm:7.1}
Let $\Sigma_j\in \A$ be a sequence of admissible surfaces. Assume
$0\in \Sigma_j$ for each $j\in \bbbn$ and let $E>0$, $p<8$ be
constants such that $\M_p(\Sigma_j)\le E$ for all $j\in \bbbn$. In
addition, assume that
\[
\sup \H^2(\Sigma_j) \le H < \infty\, .
\]
Then there is a compact $C^{1,1-\frac{8}{p}}$-manifold $\Sigma$
and a subsequence $(\Sigma_{j'})\subset (\Sigma_j)$ such that
$\Sigma_{j'}\to \Sigma$ in Hausdorff distance as $j'\to\infty$ and
moreover
\[
\M_p(\Sigma) \le \liminf_{j'\to\infty} \M_p(\Sigma_{j'})\, ,
\qquad \H^2(\Sigma) = \lim_{j'\to\infty} \H^2 (\Sigma_{j'})\, .
\]
\end{theorem}

\noindent\textbf{Remark.} The proof of this result will reveal
that the limit surface $\Sigma$ is equipped with a nice graph
representation as described in Corollary~\ref{cor:6.2}, with norms
and patch sizes uniformly controlled solely in terms of $E$ and
$p$.

\medskip\noindent\textbf{Proof of Theorem~\ref{Thm:7.1}.}
\textbf{Step 1.} We fix $j\in \bbbn$ and look at the covering
\[
\Sigma_j \subset \bigcup_{x\in \Sigma_j} B(x,R_1),
\]
where now $R_1:=\tilde{R}_1(E,p)\le R_0(E,p)$ is the radius
defined in Corollary~\ref{cor:6.2}, and $R_0(E,p)$ appeared in
\eqref{R0Ep} of Theorem~\ref{Thm:3.1}. By means of Vitali's
covering Lemma we extract a sub\-family of pair\-wise disjoint
balls $B(x_k, R_1)$, $x_k\in \Sigma_j$, such that
\begin{equation}
\label{eq:h7.2} \Sigma_j \subset \bigcup_{k} B(x,5R_1)\, .
\end{equation}
Using Theorem~\ref{Thm:3.1} for any number $N$ of
these disjoint balls (appropriately numbered) and summing with
respect to $k$, we infer
\[
N \cdot \frac \pi 2 R_1^2 \le \sum_{k=1}^N \H^2(B(x_k,R_1)\cap
\Sigma_j) \le \H^2(\Sigma_j) \le H,
\]
which means that there can be at most $\lfloor 2H/(\pi
R_1^2)\rfloor$ such disjoint balls. Therefore, \eqref{eq:h7.2}
leads to the estimate\footnote{Notice that $\tilde R_0$ depends on
$H$ and (via $R_1$) also on $E$ and $p$.}
\begin{equation}
\label{eq:h7.3} \diam \Sigma_j \le N \diam B(0,5R_1) \le
\frac{2H}{\pi R_1^2} \cdot 10 R_1 =:\tilde R_0.
\end{equation}

Since $0\in \Sigma_j$ for all $j\in\bbbn$, we find that the family
$\{\Sigma_j\}$ is contained in the closed ball $B(0,\tilde R_0)$.

\smallskip\noindent\textbf{Step 2.} Apply Blaschke's selection
theorem \cite{price} to find a compact set $\Sigma\subset
B(0,\tilde{R}_0)$ and a subsequence (still labeled with $j$) such
that
\begin{equation}
\label{eq:h7.4} \Sigma_j\to \Sigma \qquad \mbox{as $j\to \infty$}
\end{equation}
in the Hausdorff distance. Fix $\eps>0$ small (to be specified
later) and assume now that (for a further subsequence)
\begin{equation}
\label{eq:h7.4a} \dist_{\H} (\Sigma_j, \Sigma) < \frac 12 \eps R_1
\qquad \mbox{for all $j\in\bbbn$},
\end{equation}
where $\dist_{\H}(\cdot, \cdot)$ denotes the Hausdorff distance.
Next, we form an open neighbour\-hood of the limit set,
\[
\Sigma\ \subset \ B_{99\eps R_1} (\Sigma) \ \subset \
\bigcup_{y\in \Sigma} B(y, 100\eps R_1),
\]
and use Vitali's lemma again to extract a subfamily\footnote{Since
$\Sigma$ is compact, we can assume w.l.o.g. that this sub\-family
is \emph{finite\/}.} of disjoint balls $B(y_l, 100\eps R_1)$,
$y_l\in \Sigma$ for $l=1,2,\ldots, N$ such that
\begin{equation}
\label{eq:h7.5} \Sigma\ \subset \ B_{99\eps R_1} (\Sigma) \
\subset \ \bigcup_{l=1}^N B(y_l, 500\eps R_1)\, .
\end{equation}
Now, each $y_l\in \Sigma$ is a limit of some $y_l^j\in \Sigma_j$,
and according to \eqref{eq:h7.4a} we have $|y_l-y_l^j|< \frac 12
\eps R_1$ for all $l=1,\ldots, N$ and all $j\in\bbbn$. Therefore
for each fixed $j\in\bbbn$ the balls $B(y_l^j, 99\eps R_1)$ are
pairwise disjoint, since $|y_l^j-y_m^j\ge |y_l^j-y_m^j| -
|y_l^j-y_l| - |y_m-y_m^j|> 200\eps R_1 - 2 \cdot 12 \eps R_1=
199\eps R_1$. Moreover, we have
\begin{equation}
\Sigma_j \stackrel{\eqref{eq:h7.4a}}{\subset } B_{{\eps
R_1}/2}(\Sigma) \subset B_{99 \eps R_1}(\Sigma) \subset
\bigcup_{l=1}^N B(y_l, 500\eps R_1)\
\stackrel{\eqref{eq:h7.4a}}{\subset } \bigcup_{l=1}^N B(y_l^j,
501\eps R_1) \label{long-incl}
\end{equation}
for each fixed $j\in \bbbn$, since $|y-y_l^j|\le
|y-y_l|+|y_l-y_l^j|\le 501\eps R_1$ by \eqref{eq:h7.4a} for every
$y\in B(y_l, 500 \eps R_1)$. Using again Theorem~\ref{Thm:3.1} for
a fixed $j\in\bbbn$ and summing w.r.t. to $l$, we deduce
\[
N\cdot \frac \pi 2 \bigl(99\eps R_1\bigr)^2 \le \sum_{l=1}^N \H^2
(B(y_l^j, 99\eps R_1)\cap \Sigma_j) \le \H^2(\Sigma_j)\le H,
\]
whence the bound $N\le \lfloor 2H\pi^{-1} (99\eps
R_1)^{-2}\rfloor$ for the number of disjoint balls $B(y_l^j,
99\eps R_1)$ for each fixed $\eps >0$.

\smallskip\noindent\textbf{Step 3.} We consider the unit normals
$n_l^j:= n_{\Sigma_j}(y_l^j)\in \Sphere^2$ and select subsequences
finitely many times so that for all $l=1,\ldots, N$
\[
n_l^j \to n_l\in \Sphere^2 \qquad\mbox{as $j\to\infty$,}
\]
and for given small $\delta>0$  (to be specified below)
\begin{equation}
\label{eq:h7.7} |n^j_l-n_l| <\delta \qquad\mbox{for all
$j\in\bbbn$ and all $l=1,2,\ldots, N$.}
\end{equation}
Now fix $\eps>0$ so small that $ 2000\eps R_1 \le R_1$ and
\[
B(y^j_l, 2000\eps R_1) \cap \Sigma_j \subset \Phi^j_l \bigl(
D^{j,l}_{R_1} \bigr),
\]
where $\Phi^j_l(y):=y^j_l + (y,f_l^j(y))$, $y\in D^{j,l}_{R_1}
\subset T_{y^j_l} \Sigma_j \approx\bbbr^2$ is the local graph
representation of $\Sigma_j$ near $y^j_l$ on the two-dimensional
disk $D^{j,l}_{R_1} = B(0,R_1) \cap T_{y^j_l} \Sigma_j$, whose
existence is established in Corollary~\ref{cor:6.2}. If we choose
now $\delta>0$ sufficiently small (depending on $\eps$) then we
can arrange that
\[
B(y_l, 1000\eps R_1) \cap \Sigma_j \ \subset \  \tilde \Phi^j_l
\bigl( D^{l}_{\frac 56 R_1} \bigr),
\]
where $\tilde \Phi^j_l(y):= y^j_l + (y,\tilde f_l^j(y))$ for $y\in
D^l_{5R_1/6}:= B(0,\frac 56 R_1)\cap (n_l)^\perp$, and $\tilde
f_l^j$ on the \emph{fixed\/} disk $D^l_{5R_1/6}$ is obtained from
$f^j_l$ by slightly tilting the domain of $f^j_l$, i.e. by tilting
the plane $T_{y^j_l} \Sigma_j$ towards the plane $(n_l)^\perp
\approx \bbbr^2$. (That this is indeed possible is a
straightforward but a bit tedious exercise.)

The new graph functions
\[
\tilde f_l^j\colon (n_l)^\perp\supset D^l_{5R_1/6}\longrightarrow
\bbbr
\]
continue to be of class $C^{1,\lambda}$ for $\lambda=1-8/p$ with
uniform estimates for the oscillation of their gradients as in
Corollary~\ref{cor:6.2} (we use the assumption $\sup
\M_p(\Sigma_j)\le E$)  so that we may apply the theorem of
Arzela--Ascoli for each $l=1,2,\ldots, N$ to obtain subsequences
$\tilde f^{j'}_l\to f_l$ in $C^1$ as $j'\to \infty$. The limit
functions $f_l$ satisfy the same uniform $C^{1,\lambda}$
estimates. Thus, $\Sigma$ is covered by $N$ graphs
$\Phi_l(y)=y_l+(y,f_l(y))$, $l=1,2,\ldots,N$, by virtue of the
Hausdorff convergence \eqref{eq:h7.4} and the $C^1$-convergence of
the $\tilde \Phi^{j'}_l$ as $j'\to\infty$. Moreover,
\[
B(y_l,1000\eps R_1) \cap \Sigma = \Phi_l (D^l_{5R_1/6}) \cap
B(y_l,1000 \eps R_1).
\]
Now \eqref{eq:h7.5} implies that for each $y\in \Sigma$ there
exists an $l\in\{1,2,\ldots,N\}$ such that the set
\[
\Sigma \cap B(y,500\eps R_1) \stackrel{\eqref{eq:h7.5}}{\subset}
B(y_l,1000\eps R_1) \cap \Sigma
\]
so that
\[
\Sigma \cap B(y,500\eps R_1) = \Phi_l (D^l_{5R_1/6}) \cap B(y,500
\eps R_1).
\]
In particular, the limit surface $\Sigma$ is also a
$C^{1,\lambda}$ manifold for $\lambda =1 -8/p$.

\smallskip\noindent\textbf{Step 4 (lower semicontinuity of
$\M_p$).} This follows from Fatou's lemma
combined with the following properties of the inte\-grand:
\begin{eqnarray}
\K(T) & = & \lim_{i\to \infty} \K(T_i)\qquad\mbox{whenever $T_i\to
T$ and $\K(T)>0$,} \label{Ksemi1}\\
\K(T) & \le & \liminf_{i\to \infty} \K(T_i)\qquad\mbox{whenever
$T_i\to T$ and $\K(T)=0$.}\label{Ksemi2}
\end{eqnarray}
The argument is standard and uses a partition of unity in a
neighbourhood of $\Sigma$; we sketch it briefly. Take functions
$\psi_l\in C_0^\infty (B(1000\eps R_1)$, $l=1,2,\ldots, N$, such
that such that
\begin{equation}
\label{partition} \sum_{l=1}^N \psi_l \equiv 1 \qquad\mbox{on}
\quad \subset \bigcup_{l=1}^N B(y_l, 500\eps R_1)\, .
\end{equation}
This gives $\sum \psi_l \equiv 1$ on each $\Sigma_j$
for $j$ large. Inserting
\[
1=\prod_{i=0}^3 \biggl(\sum_{l_i=1}^N \psi_{l_i}(x_i)\biggr)
\]
into the integral $\M_p(\Sigma_{j'}) = \int_{(\Sigma_{j'})^4} \K
\, d\mu$ we write this integral as a sum of $N^4$ quadruple
integrals, each of them over a product of four little patches on
$\Sigma_{j'}$. Next, we use the $\Phi$'s constructed in Step~2 to
parametrize these integrals; the parameters $z_i$ (mapped to
$x_i$) belong to \emph{fixed\/} little disks $D^{l_i}$ of radius
$5R_1/6$ contained in tangent spaces to $\Sigma$. Since $\tilde
\Phi_{l}^{j'}\to \Phi_{l}$ in $C^1$, it is easy to see that all
products of $\psi_{l_i}\circ \Phi_{l_i}^{j'}(z_i)$, and all terms
where the surface measure $d\H^2(x_i)$ is expressed by $dz_i$,
converge. Combining this with \eqref{Ksemi1}--\eqref{Ksemi2},
invoking Fatou's lemma and subadditivity of $\liminf$, we see that
\[
\liminf \M_p(\Sigma_{j'}) \ge \mbox{the sum of $\liminf{}$'s of
$N^4$ terms} \ge \M_p(\Sigma).
\]
A similar argument shows that $\H^2(\Sigma_{j'}) \to
\H^2(\Sigma)$; one just replaces $\K$ by $1$ in the above
reasoning and simply passes to the limit, using the $C^1$
convergence of parametrizations. \hfill $\Box$

\medskip
\noindent {\bf Proof of Theorems \ref{thm:1.6} and
\ref{thm:1.6a}.} This follows  easily from Theorem~\ref{Thm:7.1}.
The two classes $\mathscr{C}_E(M_g)$ and $\mathscr{C}_A(M_g)$ of
surfaces $\Sigma$ which are ambiently isotopic to a fixed closed,
compact, connected, smoothly embedded reference surface $M_g$ of
genus $g$ and satisfy $\M_p(\Sigma)\le E$, or $\H^2(\Sigma)\le A$,
respectively, are nonempty. (Just take an $M_g$ of class $C^2$ to
ensure, by Proposition~\ref{prop:A.1}, that $\M_p(M_g)$ is finite;
scaling $M_g$ if necessary we can make its energy smaller than
$E$, or its area smaller than $A$.) Thus, one can take a sequence
$\Sigma_j$ contained in $ \mathscr{C}_E(M_g)$, or in
$\mathscr{C}_A(M_g)$, respectively, which is minimizing for the
area functional, or for $\M_p$. Applying Theorem~\ref{Thm:7.1}, we
obtain a subsequence of $\Sigma_j$ which converges to some
$\Sigma$ in $C^1$. Since isotopy classes are stable under
$C^1$-convergence, see \cite{blatt}, the limiting surface $\Sigma$
belongs to $\mathscr{C}_E(M_g)$, or resp. to $\mathscr{C}_A(M_g)$.

\appendix

\section{Finiteness of energy of $C^{2}$-surfaces}

\setnumbers

As before, $T=(x_0,x_1,x_2,x_3)$ stands for a tetrahedron in
$\bbbr^3$. $V(T)$ is the volume of $T$ and $A(T)$ denotes the
total area of $T$, i.e. the sum of areas of the four triangular
faces. Recall that
\begin{equation}
\label{K-int} \K(T) = \frac{V(T)}{A(T) \, (\diam T)^2}\, .
\end{equation}

\begin{proposition}
\label{prop:A.1} If $\, \Sigma\subset \bbbr^3$ is a compact,
embedded surface of class $C^2$, then there exists a constant
$C=C(\Sigma)$ such that
\[
\K(T) \le C \quad\mbox{for each $T\in \Sigma^4$.}
\]
\end{proposition}
This obviously implies that $\M_p(\Sigma)<\infty$ whenever
$\Sigma$ is of class $C^2$.

\medskip
\noindent \textbf{Proof.} Comparing $A(T)$ with the maximum of
areas of the faces, we obtain
\[
\frac 1{12} \frac{h_{\min}(T)}{(\diam T)^2}\le \K(T) \le \frac 13
\frac{h_{\min}(T)}{(\diam T)^2},
\]
where $h_{\min}(T)$ stands for the minimal height of $T$, i.e. for
the minimal distance of $x_i$ to the affine plane spanned by the
other three $x_j$'s, $i=0,1,2,3$. Since $h_{\min}(T)\le \diam T$,
it is enough to show that $\K(T)$ is bounded when $\diam T\le d_0$
for some $d_0=d_0(\Sigma)$ sufficiently small.

Thus, from now on we fix a $d_0>0$ such that for each $x\in
\Sigma$ the intersection $\Sigma \cap B(x,{2d_0})$ coincides with
a graph of a $C^2$-function defined of $x+T_x\Sigma$, and
\begin{equation}
\label{quadratic} \dist(y,x+T_x\Sigma) \le A |y-x|^2, \qquad y\in
\Sigma \cap B(x,{2d_0}).
\end{equation}
\textbf{Remark.} \eqref{quadratic} is the only thing we need from
the $C^2$-property. Such an estimate holds for $C^{1,1}$-surfaces,
too. If one represents such a surface locally by a function $g\in
C^{1,1}$ normalized to $g(0)=0$ and $\nabla g(0)=0$ then the
Lipschitz continuity of $\nabla g$ implies a quadratic height
excess as in \eqref{quadratic}.

\medskip

W.l.o.g. we can assume that $Ad_0\ll 1$.

\begin{lemma} Let $T=(x_0,x_1,x_2,x_3)$ be an arbitrary tetrahedron, with angles of
the faces denoted by $\alpha_{ij}$, $i,j=0,1,2,3$, $i\not= j$ so
that $\alpha_{ij}$ is the angle at $x_j$ on the face which is
opposite to $x_i$. Then, two cases are possible:
\begin{enumerate}
\renewcommand{\labelenumi}{{\rm (\roman{enumi})}}
\item At least one of the $\alpha_{ij}\in [\frac \pi 9, \frac{8\pi}9]$;
\item All $\alpha_{ij}\in (0,\frac \pi 9) \cup (\frac{8\pi}9,
\pi)$.
\end{enumerate}
In the latter case, eight of the $\alpha_{ij}$ are small, i.e.
belong to $(0, \frac \pi 9)$ and the remaining four are large,
i.e., belong to $(\frac{8\pi}9, \pi)$. Moreover, there is one
large angle on each face and either 0 or 2 such angles at each
vertex of $T$.
\end{lemma}

\noindent\textbf{Proof of the lemma.} We have
\begin{gather}
\sum_{0\le j \le 3, j\not=i} \alpha_{ij} = \pi \qquad\mbox{for
each $i=0,1,2,3$,} \label{180}\\
\sum_{0\le i \le 3, i\not=j} \alpha_{ij} \in (0,2\pi)
\qquad\mbox{for each $j =0,1,2,3$,} \\
\alpha_{ij} + \alpha_{lj} > \alpha_{kj} \qquad\mbox{for each
permutation $(i,j,k,l)$ of $(0,1,2,3)$.} \label{tri}
\end{gather}
(The last condition amounts to the triangle inequality for the
spherical metric.)

Now, suppose that Case (i) does not hold. If there were at most
three large angles, then the sum of all $\alpha_{ij}$ would be
strictly smaller than
\[
3\pi + 9 \cdot \frac \pi 9 = 4\pi,
\]
a contradiction. Similarly, if there were at least 5 large angles,
the sum of all angles of $T$ would be strictly larger than $4\pi$.
\heikodetail{

\bigskip

{\tt\xx By the way, I tried to find an accessible
reference for the fact that the total sum of interior angles
equals $4\pi$, without big success -- in contrast to facts
\eqref{180}--\eqref{tri}, which you find on numerous web-pages.
Finally I looked into Alexandrov's "Convex Polyeder" and combined
some results in the section where he introduces the spherical
image of a convex polyeder. So for general convex polyeders with
$n$ vertices this number $4\pi$ should change into $(n-2)2\pi$,
would you agree?}

\bigskip

}
Thus, if (i) fails, $T$ must have precisely 4 large angles. By
\eqref{180} and the pigeon-hole principle, there is precisely one
such angle on each face. Furthermore, if there is a large angle at
some vertex, then by \eqref{tri} at least one of the remaining
angles at this vertex must also be large. Since the sum of all
angles at each vertex is smaller than $2\pi$, we have precisely
either 0 or 2 large angles at each vertex. \hfill $\Box$

\bigskip Now, fix $T\in \Sigma^4$ with $d=\diam T<
d_0=d_0(\Sigma)$.

\medskip

\noindent\textbf{1.\quad} If Case (i) of the lemma holds for $T$,
we can assume w.l.o.g. that $x_0=0$, the tangent plane
$T_{x_0}\Sigma =\{(a,b,0)\mid a,b\in \bbbr\}$ is horizontal, and
$\ang(x_1,x_2)\in [\frac \pi 9, \frac{8\pi}9]$. Let $P:=\langle
x_0,x_1,x_2\rangle$. A computation shows that there is an absolute
constant $c_1$ such that
\[
\ang (P, T_{x_0}\Sigma) \le c_1 A d
\]
(which is a small angle if $d_0$ is chosen sufficiently small).
Therefore, since $\dist(x_3, T_{x_0}\Sigma) \le Ad^2$, we have
\[
\dist(x_3,P) \le c_2 A d^2,
\]
which yields $\K(T)\le c_2 A$.

\bigskip

\noindent\textbf{2.\quad} Suppose now Case (ii) holds for $T$.
W.l.o.g. we can assume that all angles at $x_0$ belong to
$(0,\pi/9)$. We can also assume that all these angles exceed $c_3
A d$ for some constant $c_3$, since otherwise there exists a
vertex and an edge of $T$ with mutual distance $\lesssim d^2$ and
we are done.

As before, we choose coordinates so that $x_0=0$ and
$T_{x_0}\Sigma =\{(a,b,0)\mid a,b\in \bbbr\}$ is horizontal.
Let $\pi_T$ stands for the orthogonal projection onto
$T_{x_0}\Sigma$.

For $i=1,2,3$, let $l_i$ be the straight line through $x_0$ and
$x_i$. Set also $x_i':=\pi_T(x_i)$, $d_i=|x_i-x_0|$,
$d_i'=|x_i'-x_0|$ and $l_i'=\pi_T (l_i)$ ($i=1,2,3$). Finally, set
$h_i=|x_i-x_i'|=\dist(x_i,x_0+T_{x_0}\Sigma)$. We have
\[
h_i \le Ad_i^2, \qquad d_i'\le d_i \le 2d_i'\, .
\]
Permuting the numbering of $x_1,x_2,x_3$, we can moreover assume
that $l_1'\not=l_3'$ (if all projections of edges meeting at $x_0$
onto the tangent plane coincide, then $V(T)=\K(T)=0$), and that the
angle $\gamma:=\ang(x_3'-x_0,x_1'-x_0)$ is the largest of all the
angles $\ang(x_j'-x_0,x_k'-x_0)$, where $j,k=1,2,3$. Set $
P:=\langle x_0,x_1,x_3\rangle.$ Note that if $\beta_i$ denotes the
angle between $l_i$ and $l_i'$, then $\sin\beta_i \le Ad_i^2/d_i
=Ad_i\le Ad \ll 1$.

Let $l\subset P$ be the straight line such that
$\pi_T(l)=l_2'=\pi_T(l_2)$. The crucial observation is that the
angle between $l$ and $l_2'$ is at most $c_4A d$ for some absolute
constant $c_4$ (here we use the piece of information that all
angles of $T$ at $x_0$ are small). Using this, we estimate
\begin{eqnarray*}
\dist(x_2,P) & \le & |x_2-x_2'| + \dist (x_2',l) \qquad\mbox{as
$l\subset P$} \\
& \le & Ad_2^2 + d_2' \sin\ang(l_2',l) \\
& \le & c_5 A d^2\, .
\end{eqnarray*}
Thus, $h_{\min}(T) \le c_5 A d^2$. This yields the desired
estimate of $\K(T)$. \hfill $\Box$

\bigskip
\noindent\textbf{Remark.} For $\Sigma$ in $C^2$, the bound that we
obtain for $\K(T)$ is of the form
\[
\K(T)\le C \cdot A,
\]
where $A$ is the maximum of the $C^2$-norms of functions that give
a graph description of $\Sigma$ in finitely many small patches.

\section{Other inte\-grands}

\setnumbers

In \cite{leger}, J.C.~L\'{e}ger suggests an integrand that could
serve as a counterpart for integral Menger curvature of
one-dimensional sets, to obtain rectifiability criteria in higher
dimensions. For $d=2$, his suggestion is to use the cube of
\begin{equation}
\label{K-Leger} \K_L (x_0,x_1,x_2,x_3) = \frac{\dist(x_3,\langle
x_0,x_1,x_2\rangle )}{\prod_{j=0}^2 |x_3-x_j|}
\end{equation}
We are going to show that $\K_L$ and some of its relatives are not
suitable for our purposes for a simple reason: even for a round
sphere, the energy given by the $L^p$-norm of such an integrand
would be infinite for all sufficiently large $p$!  This surprising
effect is due to the fact that $\K_L$ is not a symmetric function
of its variables.

To be more precise, let
\begin{equation}
\F (x,y,z,\xi) :=\frac{\dist(\xi,\langle x,y,z\rangle)
}{M(|\xi-x|, |\xi-y|, |\xi-z|)^\alpha}\, \label{bad-F}
\end{equation}
where $\alpha> 1$ is a parameter and $M\colon \bbbr_+\times
\bbbr_+\times \bbbr_+\to \bbbr_+$ is homogeneous of degree 1,
monotone nondecreasing w.r.t. each of the three variables, and
satisfies
\begin{equation}
\label{minmax} \min(t,r,s)\le M(t,r,s)\le \max(t,r,s)
\qquad\mbox{for $t,r,s\ge 0$.}
\end{equation}
Note that such $\F$ coincides with J.C.~L\'{e}ger's $\K_L$ if
$M(t,r,s)=\sqrt[3]{trs}$ is the geometric mean and $\alpha=3$.
\heikodetail{

\bigskip

Leger has $p=d+1$ which equals $3$ here so that $(\alpha-1)p=6$,
so we are kind of far away in the following proposition
from his concrete suggestion...

a sharper estimate (possibly using more of the surface of $\S^2$)
might get us closer to Leger's choice of exponents

\bigskip

}

\begin{proposition} Whenever $(\alpha-1)p\ge 12$, then
\[
\int_{\Sphere^2}\int_{\Sphere^2}\int_{\Sphere^2}\int_{\Sphere^2}
\F(x,y,z,\xi)^p \, d\H^2(x)\, d\H^2(y)\, d\H^2(z)\, d\H^2(\xi) =
+\infty\, .
\]
\label{prop-B1}
\end{proposition}

\medskip\noindent\textbf{Proof.}
\heikodetail{

\bigskip

look at my drawing on old p. 70 to understand the
strategy of patch choices to get large amounts of
energy...

\bigskip

}
We follow a suggestion of
K. Oleszkiewicz (to whom we are grateful for a brief sketch of
this proof) and consider the behaviour of $\F$ on such quadruples
of nearby points $(x,y,z,\xi)$ for which the plane $\langle x,y,z
\rangle$ is very different from the tangent plane at $\xi$. It
turns out that
\[
\int_{\Sphere^2}\int_{\Sphere^2}\int_{\Sphere^2} \F(x,y,z,\xi)^p \,
d\H^2(x)\, d\H^2(y)\, d\H^2(z)\, = +\infty\, \qquad\mbox{for each
$\xi \in \Sphere^2$}.
\]
To check this, suppose without loss of generality that
$\xi=(0,0,1)$. Fix a small $\eps\in (0,1)$ and $r_n=2^{-2n}$ for
$n=1,2,3,\ldots$. Consider the sets $\Delta_n\subset
\Sphere^2\times \Sphere^2 \times \Sphere^2$,
\begin{equation} \label{D-n}
\Delta_n := \bigl(B(a_n,{\eps r_n^2})\cap \Sphere^2\bigr)\times
\bigl(B(b_n,{\eps r_n^2})\cap \Sphere^2\bigr)\times \bigl(B(c_n,{\eps
r_n^2})\cap \Sphere^2\bigr)\, ,
\end{equation}
where
\begin{eqnarray}
a_n &:=& (r_n, 0, \sqrt{1-r_n^2}), \\
b_n &:=& (r_n, 2 r_n, \sqrt{1-5r_n^2}),\\
c_n &:=& (r_n, -2 r_n, \sqrt{1-5r_n^2}).
\end{eqnarray}
Note that for $\eps\in (0,1)$ all $\Delta_n$ are pairwise
disjoint. We shall show that whenever a triple of points
$(x,y,z)\in \Delta_n$, then the plane $P=\langle x,y,z \rangle$ is
almost perpendicular to $T_\xi \Sphere^2$ (the angle differs from
$\pi/2$ at most by a fixed constant multiple of $\eps$)
\heikodetail{

\bigskip

$|\ang (T_\xi\S^2,\langle x,y,z\rangle)-\pi/2|\le C\eps$
see also drawing in my Oct. 21st version on page 67
note also: $r_n-r_{n+1}=3r_n/4>3r_n^2$ for all $n\in\N$ to get
disjointness
.
\bigskip

}
and
\[
\dist(\xi, P)\ge r_n/2, \qquad \F^p(x,y,z,\xi)\ge A \cdot
r_n^{p(1-\alpha)}
\]
for some constant $A$ depending on $\eps$, $p$ and $\alpha$ but
\emph{not\/} on $n$. Let $v_n:=b_n-a_n$, $w_n:=c_n-a_n$
($n=1,2,\ldots$). Since $\sqrt{1-x}=1-x/2 + O(x^2)$ as $x\to 0$,
we have
\[
v_n =  (0, 2 r_n, -2r_n^2 + O(r_n^4)), \qquad w_n = (0, -2 r_n,
-2r_n^2 + O(r_n^4))\, .
\]
A computation shows that
\[
u_n:=v_n \times w_n = (-8r_n^3, 0, 0) + e_n, \qquad |e_n|\le C_1
r_n^5\, ,
\]
where $C_1$ is an absolute constant. Therefore,
\[
\sigma_n:=\frac{u_n}{|u_n|} = (-1,0,0) + f_n, \qquad |f_n|\le C_2
r_n^2\, ,
\]
again with some absolute constant $C_2$. Now, let $(x,y,z)\in
\Delta_n$ and let $v_n':=y-x$, $w_n':=z-x$. By triangle
inequality, we have
\[
\max(|v_n-v_n'|,|w_n-w_n'|)\le 2\eps r_n^2,
\]
so that another elementary computation shows that
$\sigma_n':=(v_n'\times w_n')/|v_n'\times w_n'|$ satisfies
\[
|\sigma_n-\sigma_n'|\le C_3 \eps, \qquad n=1,2,3,\ldots
\]
for $\epsilon $ sufficiently small.
\heikodetail{

\bigskip

This is the crucial estimate one would have to work with when
one wants to improve to get $\mathcal{F}$ infinite, this here
is based on
\begin{eqnarray*}
v_n'\times w_n' & = & [v_n+(v_n'-v_n)]\times
 [w_n+(w_n'-w_n)]\\
 & = & u_n+2\eps r_n^2|r_n|+O(r_n^4)
 \end{eqnarray*}
 as $n\to\infty$ so that
 $$
 6r_n^3\le |v_n'\times w_n'|\le 10r_n^3\quad\Fo \eps \ll 1.
 $$

 \bigskip

 }
Moreover,
\begin{eqnarray}
\dist(\xi, \langle x,y,z \rangle) & = &|(\xi-x)\cdot \sigma_n'| \nonumber \\
& = & \bigl|\bigl((\xi-a_n) + (a_n-x)\bigr)\cdot (\sigma_n +
(\sigma_n'-\sigma_n)\bigr| \nonumber \\
& \ge & r_n - C_4 \eps r_n \ = \ \frac{r_n}2 \label{largedist}
\end{eqnarray}
if we choose $\eps=1/2C_4$. By \eqref{minmax}, we also have
\begin{equation}
\label{M=rn} M(|\xi-x|,|\xi-y|,|\xi-z|) \approx r_n, \qquad
(x,y,z)\in \Delta_n, \quad n=1,2,3,\ldots
\end{equation}
Combining \eqref{largedist} and \eqref{M=rn}, we estimate
\begin{eqnarray*}
\lefteqn{ \int_{\Sphere^2}\int_{\Sphere^2}\int_{\Sphere^2}
\F(x,y,z,\xi)^p \,
d\H^2(x)\, d\H^2(y)\, d\H^2(z)\,} \\
& \gtrsim & \sum_{n=1}^\infty \int \!\!\int
\!\!\int_{\{(x,y,z)\in
\Delta_n\}} \F(x,y,z,\xi)^p \, d\H^2(x)\, d\H^2(y)\, d\H^2(z) \\
& \gtrsim & \sum_{n=1}^\infty (\pi \eps r_n^2)^6
\frac{r_n^p}{r_n^{\alpha p}} \\
& \approx & \sum_{n=1}^\infty (r_n)^{12+(1-\alpha)p} \\
& = & +\infty \qquad\mbox{for $(\alpha-1)p\ge 12$.}
\end{eqnarray*}
This completes the proof. \qed

\medskip\noindent\textbf{Remark.} One can check that a similar argument
shows that
\[
\int_U \int_U \int_U \int_U \F^p = +\infty \qquad \mbox{if
$(\alpha-1)p\ge 12$}
\]
whenever $U$ is a patch of of a $C^2$ surface $\Sigma\subset
\bbbr^3$ such that the Gaussian curvature of $\Sigma$ is strictly
positive on $U$.

\smallskip

The phenomenon described in Proposition~\ref{prop-B1} does not
appear for the integrand
\[
\K_R(x,y,z,\xi)= 1/R(x,y,z,\xi)\, ,
\]
where $R(x,y,z,\xi)$ denotes the radius of a circumsphere of four
points of the surface --- we simply have $1/R=\mathrm{const}$ for
all quadruples of pairwise distinct points of a round sphere.
However, one can easily find examples of smooth surfaces for which
$1/R\to \infty$ at some points: take e.g. the graph of $f(x,y)=xy$
near 0. It contains two straight lines and for every $\delta>0$
there are lots of triangles with all vertices on these lines, all
angles (say) $\ge \pi/6$ and diameter $\le \delta$. For each such
triangle $\Delta$ one can take a sphere $S$ which has the
circum\-circle of $\Delta$ as the equatorial circle. The radius of
$S$ is $\lesssim \delta$ and $S$ intersects the graph of $f$ at
infinitely many points that are not co\-planar with vertices of
$\Delta$.



\small \vspace{1cm}
\begin{minipage}{56mm}
{\sc Pawe\l{} Strzelecki}\\
Instytut Matematyki\\
Uniwersytet Warszawski\\
ul. Banacha 2\\
PL-02-097 Warsaw \\
POLAND\\
E-mail: {\tt pawelst@mimuw.edu.pl}
\end{minipage}
\hfill
\begin{minipage}{56mm}
{\sc Heiko von der Mosel}\\
Institut f\"ur Mathematik\\
RWTH Aachen\\
Templergraben 55\\
D-52062 Aachen\\
GERMANY\\
Email: {\tt heiko@}\\{\tt instmath.rwth-aachen.de}
\end{minipage}


\begin{thebibliography}{9999}
\footnotesize

\bibitem{allard}
 Allard, W. K.
 On the first variation of a varifold. \emph{Ann. of Math. (2)\/}
 \textbf{95} (1972), 417--491.


\bibitem{banavargmm}
Banavar, J.R.; Gonzalez, O.; Maddocks, J.H.; Maritan, A.
Self-interactions of strands and sheets. \emph{J. Statist.
Phys.\/} {\bf 110} (2003), {35--50}.
%

\bibitem{blatt}
Blatt, S. Note on continuously differentiable isotopies. Preprint
Nr. 34, Institut f. Mathematik, RWTH Aachen (2009).

\bibitem{BlM70}
Blumenthal, L.M.; Menger, K. {\it Studies in geometry.} Freeman
and co., San Francisco, CA,  1970.



\bibitem{daviddepauwtoro}
David, G.; De Pauw, T.; Toro, T. A generalization of Reifenberg's
theorem in $\R\sp 3$.  \emph{Geom. Funct. Anal.\/}  \textbf{18}
(2008), no. 4, 1168--1235.

\bibitem{davidkenigtoro}
David, G.; Kenig, C.; Toro, T. Asymptotically optimally doubling
measures and Reifenberg flat sets with vanishing constant.
\emph{Comm. Pure Appl. Math.}  \textbf{54}  (2001),  no. 4,
385--449.

\bibitem{davidsemmes-QM}
David, G.; Semmes, S.  Quasiminimal surfaces of codimension $1$
and John domains.  \emph{Pacific J. Math.\/} \textbf{183}  (1998),
no. 2, 213--277.

\bibitem{davidsemmes}
David, G.; Semmes, S. {\it Analysis of and on uniformly
rectifiable sets.} Math. Surveys \& Monographs {\bf 38}, AMS,
Providence, Rhode Island, 1993.

\bibitem{dnf}
Dubrovin, B.A.; Fomenko, A.T.; Novikov, S.P. \emph{Modern
geometry---methods and applications.\/} Part II. \emph{The
geometry and topology of manifolds.\/} Graduate Texts in
Mathematics, 104. Springer-Verlag, New York, 1985.

\bibitem{fenchel}
Fenchel,W. On total curvatures of Riemannian manifolds: I.
\emph{J. London Math. Soc.\/} \textbf{15}, (1940). 15--22

\bibitem{GT98}
Gilbarg, D.; Trudinger, N.S.
\emph{Elliptic partial differential equations of second order.\/}
 Reprint of the 1998-ed., Springer, Berlin, New York, 2001.


\bibitem{GM99}
Gonzalez, O.;   Maddocks, J.H. Global  curvature, thickness, and
the ideal shape of knots. \emph{Proc. Natl. Acad.  Sci. USA\/}
{\bf 96} (1999), 4769--4773.

\bibitem{GMSvdM02}
Gonzalez, O.; Maddocks, J.H.; Schuricht, F.; von der Mosel, H.
Global curvature and self-contact of nonlinearly elastic curves
and rods. \emph{Calc. Var. Partial Differential Equations\/} {\bf
14} (2002), 29--68.



\bibitem{leger}
L{\'e}ger, J.C.
 {Menger curvature and rectifiability}.
 {\emph{Ann. of Math.\/} (2)}
 {\bf 149}
 (1999),
 {831--869}.

\bibitem{LW08a}
Lerman, G.; Whitehouse, J.T. High-Dimensional Menger-Type
Curvatures -- Part~I: Geo\-met\-ric Multipoles and Multiscale
Inequalities. arXiv:0805.1425v1 (2008).


\bibitem{LW08b}
Lerman, G.; Whitehouse, J.T. High-Dimensional Menger-Type
Curvatures -- Part II: $d$-Se\-pa\-ra\-tion and a Menagerie of
Curvatures.  \emph{J. Constructive Approximation\/} \textbf{30}
(2009), 325--360.

\bibitem{lima}
Lima, E.L. Orientability of smooth hypersurfaces and the
Jordan--Brouwer separation theorem. \emph{Expo. Math.\/}
\textbf{5} (1987), 283--286.

\bibitem{maly}
Mal\'{y}, J. Absolutely continuous functions of several variables.
\emph{J. Math. Anal. Appl.\/} {\bf 231} (1999), no. 2, 492--508.



\bibitem{Ma98}
Mattila, P. Rectifiability, analytic capacity, and singular
integrals. Proc. ICM, Vol. II (Berlin 1998),  \emph{Doc. Math.\/}
{\bf 1998}, Extra Vol. II, 657--664 (electronic).

\bibitem{Ma04}
Mattila, P. Search for geometric criteria for removable sets of
bounded analytic functions.  \emph{Cubo} {\bf 6} (2004), 113--132.


\bibitem{Mel95}
Melnikov, M. Analytic capacity: discrete approach and curvature of
measure. \emph{Sb. Mat.\/} {\bf 186} (1995), 827--846.

\bibitem{MelV95}
Melnikov, M.; Verdera, J. A geometric proof of the $L^2$
boundedness of the Cauchy integral on Lipschitz curves.
\emph{Intern. Math. Research Notices\/} {\bf 7} (1995), 325--331.

\bibitem{menger}
Menger, K. Untersuchungen \"uber allgemeine Metrik. Vierte
Untersuchung. Zur Metrik der Kurven. \emph{Math. Ann.\/} {\bf 103}
(1930), 466--501.




\bibitem{koles}
Oleszkiewicz, K. Private communication.


\bibitem{P02}
Pajot, H. {\it Analytic capacity, rectifiability, Menger curvature
and the Cauchy integral.} Springer Lecture Notes 1799, Springer
Berlin, Heidelberg, New York, 2002.


\bibitem{ptt}
Preiss, D.; Tolsa, X.; Toro, T. On the smoothness of H\"{o}lder
doubling measures. \emph{Calc. Var. Partial Differential
Equations\/} \textbf{35} (2009), no. 3, 339--363.

\bibitem{price}
Price, G.B. On the completeness of a certain metric space with an
application to Blaschke's selection theorem. \emph{Bull. Amer.
Math. Soc.\/} \textbf{46} (1940), 278--280.

\bibitem{reifenberg}
Reifenberg, E.R. Solution of the Plateau Problem for
$m$-dimensional surfaces of varying topological type.  {\it Acta
Math.} \textbf{104}  (1960), 1--92.



\bibitem{SvdM03a}  Schuricht, F.; von der Mosel, H.
Global curvature for rectifiable loops.  \emph{Math. Z.\/} {\bf
243} (2003), 37--77.

     %
\bibitem{SvdM03b}
Schuricht, F.; von der Mosel, H. Euler-Lagrange equations for
nonlinearly elastic rods with self-contact. \emph{Arch. Rat. Mech.
Anal.\/} {\bf 168} (2003), 35--82.
      %
\bibitem{SvdM04}
Schuricht, F.; von der Mosel, H. Characterization of ideal knots.
\emph{Calc. Var. Partial Differential Equations\/} {\bf 19}
(2004), 281--305.

\bibitem{Si93}
Simon, L.: Existence of surfaces minimizing the Willmore
functional. \emph{Comm. Anal. Geom.\/} \textbf{1} (1993), no. 2,
281-326.



\bibitem{ssvdm-triple}
Strzelecki, P.; Szuma\'{n}ska, M.; von der Mosel, H. Regularizing
and self-avoidance effects of integral Menger curvature. Preprint
Nr. 29, Institut f. Mathematik, RWTH Aachen (2008). To appear in
Annali Scuola Norm. Sup. di Pisa.

\bibitem{StvdM05}
Strzelecki, P.; von der Mosel, H. On a mathematical model for
thick surfaces. In: Calvo, Millett, Rawdon, Stasiak (eds.) {\it
Physical and Numerical Models in Knot Theory}, pp. 547--564. Ser.
on Knots and Everything 36, World Scientific, Singapore, 2005.


\bibitem{svdm-global}
Strzelecki, P.; von der Mosel, H. Global curvature for surfaces
and area minimization under a thickness constraint. \emph{Calc.
Var. Partial Differential Equations\/} \textbf{25} (2006),
431--467.

\bibitem{marta}
Szuma\'{n}ska, M. Integral versions of Menger curvature: smoothing
potentials for rectifiable curves. (In Polish.) Ph.D. thesis,
Technical University of Warsaw, 2009.


\bibitem{Toro-w22}
Toro, T. Surfaces with generalized second fundamental form in
$L\sp 2$ are Lipschitz manifolds.  \emph{J. Differential Geom.\/}
39 (1994), no. 1, 65--101.


\bibitem{V01} Verdera, J. The $L^2$ boundedness of the Cauchy
integral and Menger curvature. In {\it Harmonic analysis and
boundary value problems}, pp. 139--158. Contemp. Math. {\bf 277},
 AMS, Providence, RI, 2001.


\end{thebibliography}
\end{document}